\documentclass[final]{elsarticle}
\usepackage{latexsym, graphicx, epsfig, amsmath, amsfonts,amssymb,bm}
\usepackage{subfigure,cases}
\usepackage[outdir=./Figures/]{epstopdf}
\usepackage{xcolor}
\usepackage{colortbl,booktabs}
\usepackage{float}
\usepackage{cals}

\usepackage{geometry}
\newtheorem{lemma}{Lemma}
\newtheorem{theorem}{Theorem}
\newtheorem{remark}{Remark}
\newtheorem{example}{Example}
\newenvironment{proof}{\paragraph{Proof}}{\hfill$\square$}

\def\u{{\bf u}}
\def\f{{\bf f}}

\def\c{\mathbf{c}}

\def\dfr#1#2{\displaystyle{\frac{#1}{#2}}}

\def\be{\begin{equation}}
\def\ee{\end{equation}}

\def\w{\omega}

\def\x{\mathbf{x}}
\def\c{\mathbf{c}}

\def\R{{\mathbb R}}

\newcommand{\argmin}{\operatornamewithlimits{argmin}}

\begin{document}

\begin{frontmatter}
\title{Numerical Aspects for Approximating Governing Equations Using Data} 


\author{Kailiang Wu}
\ead{wu.3423@osu.edu}

\author{Dongbin Xiu\corref{mycorrespondingauthor}}
\cortext[mycorrespondingauthor]{Corresponding author}
\ead{xiu.16@osu.edu}

\address{Department of Mathematics, The Ohio State University, Columbus, OH 43210, USA.}

\begin{abstract}
We present effective numerical algorithms for locally 
recovering unknown governing differential equations 
from measurement data. We employ a set of standard basis functions, e.g.,
polynomials, to approximate the governing equation with high accuracy.
Upon recasting the problem into a function approximation problem,  we discuss
several important aspects for accurate approximation. Most notably, we discuss the 
importance of  using a large number of short bursts of trajectory data, rather than using data from a single
long trajectory. Several options for the numerical algorithms to perform accurate approximation are then presented,
along with an error estimate of the final equation approximation.
We then present an extensive set of numerical examples of both linear and nonlinear systems to demonstrate the properties and effectiveness
of our equation recovery algorithms.

\end{abstract}

\end{frontmatter}

\section{Introduction} \label{sec:intro}

Recently there has been a growing interest in discovering governing
equations of certain physical problems using observational data. It is
a result of the wide recognition that all physical laws, e.g., Newton's
law, Kepler's law, etc., were established based on empirical observational
data. Since all models and physical laws are approximations to the true
physics, a numerically constructed law or governing equation, which albeit
may not possess a concise analytical form, could serve as a good
model, so long as it is reasonably accurate.

Early work in this direction includes
\cite{bongard2007automated,schmidt2009distilling}, where symbolic
regression was used to select physical laws and 
determine the structure of the underlying dynamical system.
More recently, 
a sparsity-promoting method \cite{brunton2016discovering} was proposed to discover
governing equations from noisy measurement data. 
It was based on a sparse regression method (\cite{tibshirani1996regression})
and a Pareto analysis to construct parsimonious 
governing equations from a combinatorially large set of candidate models. 
This methodology was further applied to recovering  
partial differential equations \cite{rudy2017data,schaeffer2017learning}, as well as inferring structures of network
models \cite{mangan2016inferring}.  
Conditions were provided in \cite{tran2017exact} for recovering 
the governing equations from possibly highly corrupted measurement data, 
along with a sparse convex optimization method to recover 
the coefficients of the underlying equations in the space of multivariate polynomials.  
A nonconvex sparse regression approach  
was proposed in \cite{schaeffer2017sparse} 
for selection and
 identification of a dynamical system directly from noisy data. 
A model selection approach was developed in 
\cite{Mangan20170009} for dynamical systems via sparse regression and information criteria. 
Combining delay embedding and Koopman theory, 
a data-driven method was designed in \cite{brunton2017chaos} to decompose chaotic systems as an intermittently forced linear system. 
Very recently, 
three sampling strategies 
were developed in \cite{schaeffer2017extracting} 
for learning quadratic high-dimensional differential equations from under-sampled data. 
There are other techniques that address various aspects of dynamical 
system discovery problem. 
These include, but are not limited to,
methods to 
reconstruct the equations of motion from a data series \cite{crutchfield1987equations}, 
artificial neural networks \cite{gonzalez1998identification},
nonlinear regression \cite{voss1999amplitude}, 
equation-free modeling \cite{kevrekidis2003equation}, 
normal form identification \cite{majda2009normal}, 
empirical dynamic modeling \cite{sugihara2012detecting,ye2015equation}, 
nonlinear Laplacian spectral analysis \cite{giannakis2012nonlinear}, 
modeling emergent behavior \cite{roberts2014model},  
and automated inference of dynamics \cite{schmidt2011automated,daniels2015automated,daniels2015efficient}.  
The readers are also referred to the introduction in 
\cite{BruntonSIAMNews2017}. 

The focus of this paper is on the study of local recovery/approximation of 
 general first-order nonlinear ordinary differential/algebraic equations 
 from the time-varying measurement data.
The major new features of the proposed numerical framework include the following.
First, our method seeks ``accurate approximation'' to the underlying
governing equations.  This is different from many of the existing
studies, where ``exact recovery'' of the equations is the goal. The
justification for ``approximation'', rather than ``recovery'', is that
all governing equations will eventually be discretized and
approximated by certain numerical schemes during implementation. As
far as computations are concerned, all numerical simulations
are conducted using ``approximated'' equations. Therefore, as long as
the errors in the equation approximation/recovery step is
small, comparable to the discretization errors, the
approximated equations can be used as sufficient model for
simulations.
By relaxing the requirement from ``exact recovery'' to ``accurate
approximation'', we are able to utilize many standard approximation
methods using standard basis functions such as
polynomials, which are used in this paper. Consequently, one does not need to choose a
(combinatorially) large dictionary in order to include all possible
terms in the equations. Since most of the governing
equations we encounter in practice have relatively simple and smooth
forms, polynomial approximation can be highly accurate at
moderately low degrees.
Secondly, 
we propose the use of multiple, in fact
a large number of, trajectory data of very short span, instead of using a
single long trajectory data as in many existing studies. The length of
each burst of the trajectory can be exceptionally short -- it 
needs to be merely long enough to enable accurate estimate of the time
derivatives. For noiseless trajectories, the length can as little as 2,
for first-order estimation of the time derivative, or 3 for
second-order estimation, and so on. We do require that the number of
such short bursts of trajectories to be large so that they provide
``good'' samples in the domain of interest.
Upon doing so, we convert the equation recovery problem into a
function approximation problem with over-sampled data. A large
number of well established approximation algorithms can then be employed.
(The readers are referred to earlier books such as \cite{Cheney66,
 	Davis75, Powell81, Rivlin69} for the classical results and
      methods in function approximation theory.)
In this paper, we discuss the results obtained by using least squares
method, which is one of the most established methods; $\ell_1$
regression method, which is capable of removing data corruption
(\cite{candes2005decoding, ShinXiu}); and a matrix-free sequential
algorithm (\cite{ShinXiu_SISC17,
 	WuShinXiu_SISC17, WuXiu_JCP18}), which is suitable for
      extraordinarily large data sets or stream data.
We remark that the use of a very lagre number of trajectory bursts
does not necessarily induce large simulation cost, as each burst is of
very short length. In fact, in many cases the total amount of data in
the large number of bursts is less than that in a single long trajectory.
We also remark that the idea of using multiple bursts of trajectories
has also be discussed in \cite{schaeffer2017extracting}. The major
difference between \cite{schaeffer2017extracting} and this work is
that the goal in \cite{schaeffer2017extracting} is for exact recovery
of the equation of quadratic form. The bursts are random selected
based on the theory of sparse recovery (compressed sensing). In this
paper, the bursts are selected using approximation theory, as we only
seek to approximate the equations and not recover them. And the
equations can be of arbitrary form, beyond polynomials and including
algebraic constraints. The choices of
the bursts can be random or deterministic, as long as they enable one
to obtain accurate approximation.

This paper is organized as follows. After the basis problem setup in
Section \ref{sec:setup}, we present the main method in Section
\ref{sec:method}, which includes the general formulation, basis
selections, trajectory selections, choices of approximation
algorithms, and an error bound. We then present, in Section
\ref{sec:examples}, an extensive set of numerical examples, covering a
wide variety of linear/nonlinear differential/algebraic equations, to
demonstrate the effectiveness of the proposed algorithms.

\section{Problem Setup} \label{sec:setup}

We consider an 
autonomous system in the form of
\begin{equation}\label{eq:ODE}
\dfrac{d {\bf u}(t)}{dt}={\bf f} ({\bf u}), \qquad \u(t_0) = \u_0,
\end{equation}
where ${\bf u}=(u_1,u_2,\cdots,u_d)^\top$, $d\geq 1$, are state
variables, $t_0$ the initial time, and
${\bf f}$ defines the evolution of the states. Throughout this
paper we call $\u(t; \u_0)$ a trajectory, whose
dependence on the initial state $\u_0$ is explicitly shown. The
right-hand-side ${\bf f}$ is the ``true'' and unknown dynamics.
Our goal is to construct an accurate representation of ${\bf f}$ using
measurement data of the state variables ${\bf u}$.

\subsection{Domain-of-Interest}

While the state variables $\u$ may reside in the entire space
$\R^d$, we restrict our equation recovery to a bounded domain
$D\subset \R^d$. We refer to $D$ as the domain-of-interest, which stands for
the range of the state variables where one focuses on the underlying
problem. Depending on the (unknown) dynamical system and one's
interest in studying it, the definition of the domain-of-interest
varies.

Let
\begin{equation}
U_i = \left[u_{\min,i}, u_{\max, i}\right], \qquad 1\leq i\leq d,
\end{equation}
be a closed interval and represent a range of the state variable
$u_i$, $1\leq i\leq d$,
where one wishes to understand the behavior of the system. We then
define $D$ as the product of these intervals, i.e.,
\begin{equation}
D = \times_{i=1}^d U_i.
\end{equation}
Effectively this creates a hypercube in $d$ dimensions. We then equip
it with a measure $\w$, which is 
absolutely continuous and (without loss of generality) satisfies $\int_D
d \omega = 1$. We then invoke the standard definition of inner product
\begin{align*}
\left(  g,  h \right)_{L^2_\omega} := \int_D  g({\bf x})  h({\bf x})
  d\omega({\bf x}),
\end{align*}
and its induced norm $\|\cdot \|_{L^2_\omega}$. 

\begin{remark}
Upon defining the domain-of-interest, we have restricted our equation
recovery effort to a ``local'' manner. This is because recovering the
equation in the entire space $\R^d$, albeit desirable, may be highly challenging, when the
underlying dynamical system exhibits complex behaviors that are different in different
parts of the space. The use of the domain-of-interest allows one to
focus on a region of the domain where the problem can be well behaved and
${\bf f}$ possesses a regular/smooth form, which would be highly advantageous
for the recovery effort. The size of $D$ depends on the underlying
problem and is not necessarily small.
\end{remark}

\subsection{Data}

Let $t_0<t_1<\cdots $ be a sequence of time instances.
We assume that the true state variables $\u$ can be observed on these
time instances and denote
$$
{\bf x}(t_j; \u_0) = {\bf u}(t_j; \u_0) + {\bm \epsilon}_j, \quad
j=0,1,\dots,
$$
the observed data for the state trajectory, where
${\bm \epsilon}_j$ stands for observation errors. The
trajectory data can be from different sequences of trajectories, which
are originated by different
initial states.
Let $\u_0^{(1)},\dots, \u_0^{(M)}$, $M\geq 1$, be a set of different initial
states. We use the following shorthanded notation to denote the
trajectory data from the $m$-th initial condition,
\begin{equation} \label{x}
{\bf x}^{(m)}(t_{j}) = {\bf u}(t_{j};\u_0^{(m)}) + {\bm \epsilon}_{j}^{(m)}\in D,
\qquad 0\leq  j \le J_m, \quad 1\leq m \leq M,
\end{equation}
where $J_m$ is the number of intervals in the $m$-th trajectory data.
The collection of all these data shall be used to estimate the true
state equation. 

\subsection{Extension to Differential-Algebraic Equations}

It is worth noting that the method discussed in this paper is
applicable to semi-explicit
	differential algebraic equations (DAEs), i.e., DAEs with
        index 1, in the following form,
	\begin{subnumcases}{\label{eq:DAE}}
	\dfrac{d {\bf u}(t)}{dt}={ {\bf F}}({\bf u},{\bf v}),\label{eq:DAEa} \\
	{ {\bf G}} ({\bf u},{\bf v})=0, \label{eq:DAEb}
	\end{subnumcases} 
	where $\u$ and ${\bf v}$ are the state variables and ${\bf G}$
        represents algebraic constraints.
	Let us assume that the implicit function theorem holds such
        that the constraint \eqref{eq:DAEb} induces the following
        explicit relation
	\begin{equation}\label{eq:DAEnewb}
	{\bf v} = {\bf g} ({\bf u}).
	\end{equation}
Then, by letting
${\bf f}({\bf u}) = { {\bf F}} ({\bf u},{\bf g}({\bf u}))$, the DAEs
can be written in the same form as the autonomous system \eqref{eq:ODE}.

\section{Main Numerical Algorithm} \label{sec:method}

In this section we discuss numerical aspects of our equation
approximation/recovery method. We first present the general formulation, which is rather
similar to many of the existing work (e.g., \cite{bongard2007automated,brunton2016discovering,schaeffer2017sparse}). The  details of the numerical algorithm,
which is different and new, are then presented, along with an error analysis.

\subsection{General Formulation}

Let $\{\x\}$ be a collection of trajectory data, in the form of
\eqref{x}. We
assume that the time derivatives of the trajectory are also available at some time instances $t_j, j \in \Theta_m \subseteq \{ 0,1,\cdots,J_m \}$,
either produced by the measurement procedure or approximated
numerically using $\{\x\}$. 
Let
$$
{\bf {\dot x}} (t) =  {\bf {\dot u}}(t) + \bm {\tau}(t),
$$
be the derivative data, where $\bm{\tau}(t)$ stands for the errors
or noises. 
This creates a set of data pairings at the time instances $\{t_j,j\in \Theta_m \}$, where both
the state variables and their derivatives are available. The
collection of such pairings are denoted as
\begin{equation} \label{pairing}
\left\{ {\bf x}^{(m)}(t_{j}), \dot{\bf x}^{(m)}(t_{j}) \right\}, \qquad j \in \Theta_m,\quad 1\le m \le M,
\end{equation}
where the total number of data pairings is 
$$
M_{tot} = \sum_{m=1}^M \# \Theta_m.
$$

We then choose a (large) set of dictionary $\Phi=\{\phi_1,\dots,\phi_N\}$ to represent $\f$ in the
governing equation \eqref{eq:ODE}. For the vector function $\f =
(f_1,\dots f_d)$, its
approximation is typically conducted component-wise. That is, for each
$\ell=1,\dots,d$, we seek
\begin{equation}
f_\ell(\x) \approx \tilde{f}_\ell (\x) := \sum_{j=1}^N c_{j,\ell} \phi_j(\x).
\end{equation}
Here the functions
$\{\phi_j (\x) \}$ in the dictionary should contain all the
possible functions one wishes to incorporate, for the accurate
approximation of $\f$. Some examples are
\begin{equation} \label{Phi}
\Phi = \{{\bf 1}, \x, \x^2,\dots, \x^n, \sin\x, \cos\x, \dots, \sin m\x, \cos m\x,
\exp(\pm\x), \dots, \exp(\pm k\x), \dots\}.
\end{equation}
Let $\mathbf{A}$ be a $(M_{tot}\times N)$ model matrix
$$
\mathbf{A} = (a_{ij}), \qquad a_{ij} = \phi_j(\x_i), \qquad 1\leq i
\leq M_{tot}, \quad 1\leq j\leq N,
$$
where $(\x_i,{\dot \x}_i)$ represents the $i$-th data pairing in \eqref{pairing}. 
The equation recovery problem is then cast into an approximation
problem 
\begin{equation} \label{Ax}
\mathbf{A c}_\ell \approx \dot{\mathbf{X}}_\ell , \qquad \ell=1,\dots, d,
\end{equation}
where $\dot{\mathbf{X}}_\ell$ is the vector of $\ell$-th component of all available derivatives in
the pairings \eqref{pairing}, and $\mathbf{c}_\ell=(c_{1,\ell}\dots, c_{N,\ell})^T$ is the coefficient
vector to be solved.

Once an accurate approximation is obtained, we consider
\begin{equation}{\label{eq:appODE}}
\dfrac{d {\bf x}(t)}{dt}=\tilde{\bf f}({\bf x}),
\end{equation}
as an approximation of the true system \eqref{eq:ODE}.

It should be noted that converting the equation recovery problem into
an approximation problem in the form of \eqref{Ax} is a rather
standard framework in the existing studies. In most of the existing
studies, one typically chooses the
dictionary \eqref{Phi} to be very large (often combinatorially large) and then use a certain
sparsity-promoting method to construct a parsimonious solution, e.g.,
$\ell_1$ minimization for underdetermined system. By doing so, if the
terms in the true equations are already included in the dictionary
\eqref{Phi}, it is then possible to exactly recover the form of the
governing equations.

Hereafter, we will present the detail of our numerical approach, which
is different from most of the existing approaches.

\subsection{Equation Approximation using Standard Basis}

Contrary to most of the existing studies,
our method here does not utilize a
combinatorially large dictionary to include all possible terms in the
equations. Instead, we propose to use  a ``standard'' set of
basis functions, particularly polynomials in term of $\x$, to
approximate $\f$ directly. 
One certainly is not restricted to polynomials and can use other basis
such as radial basis functions. Here we choose polynomials merely because
it is one of the most widely used basis for function
approximation. 
By doing so, we
give up the idea of ``exact recovery'' of $\f$, as exact recovery is
not possible unless the true governing equations are strictly of
polynomial form (which in fact is not uncommon). Instead,  our approach focuses on
obtaining an approximation of $\f$ with sufficiently high accuracy. We
remark that this is satisfactory for practical computations, as the
exact governing equations (even if known) are often discretized and approximated by a
certain numerical scheme after all. So long as the errors in the
equation recovery/approximation step is sufficiently small, compared
to the discretization errors, the approximated governing equations can
be as effective as the true model. Even though our method will not
recover the true equations exactly in most cases, we shall continue to
loosely use the term ``equation recovery'' hereafter, interchangeable
with ``equation approximation''.

To construct the approximation
$\tilde{\bf f}$, we confine ourself to a finite dimensional
polynomial space $ V \subset L^2_\omega(D)$ with $\dim V=N\geq 1$.
Without loss of generality, we take $d\omega = \frac{1}{\int_D d{\bf x}}d{\bf x}$, and  
let the subspace $V$ be $\Pi_n^d$, the linear subspace of
polynomials of degree up to $n\geq 1$. That is,
\begin{equation} \label{Poly-space}
\Pi_n^d = \text{span} \{ {\bf x}^{{\bf i}} = x_1^{i_1}\cdots x_d^{i_d}, |{\bf i}| \le n\},
\end{equation}
where ${\bf i} = (i_1,\dots,i_d)$ is multi-index with $|{\bf i}|=i_1 +
\cdots + i_d$. Immediately, we have
\begin{equation} \label{Poly-space-dim}
N = \dim \Pi_n^d = {n+d \choose d} =\frac{(n+d)!}{n!d!}.
\end{equation}

Let $\{\phi_j( {\bf x} )\}_{j=1}^N$ be a basis of $V$. 
We then seek an approximation $\tilde{\bf f} \approx \f$ in a
component-wise manner. That is, for the $\ell$-th component,
we seek
\begin{equation} \label{f_approx}
\tilde  f_\ell ( {\bf x} ) =\sum_{j=1}^N  c_{j,\ell} \phi_j( {\bf x} ),\quad \ell=1,\dots,d.
\end{equation}
Upon adopting the following vector notations
\begin{equation} \label{vector1}
{\bf \Phi}( {\bf x} ) = \big( \phi_1( {\bf x} ), \cdots, \phi_N( {\bf x} ) \big)^{\top},
\qquad
{\bf c}_\ell = \big( c_{1,\ell}, c_{2,\ell},\cdots,c_{N,\ell} \big)^\top, 
\end{equation}
we write
\begin{equation}
\tilde  f_\ell ( {\bf x} ) = \langle {\bf c}_\ell, {\bf \Phi}( {\bf x}
) \rangle, \qquad \ell = 1,\dots, d.
\end{equation}
Our goal is to compute the expansion coefficients in ${\bf c}_\ell$, $\ell=1,\dots,d$.
Later in Section \ref{sec:strategy} we will introduce several algorithms to
compute the coefficients ${\bf c}_\ell$.


Unless the form of the true system \eqref{eq:ODE} is of
polynomial, we shall not obtain exact recovery of the
equation. For most physical systems, their governing equations
\eqref{eq:ODE} consist of smooth functions, we then expect 
polynomials to produce highly accurate approximations with low or 
modest degrees. 
On the other hand, if one possesses sufficient prior knowledge of the
underlying system and is certain
that some particular terms should be in the governing equations,
then those terms can be incorporated in the dictionary, in addition to
the standard polynomial basis.

\subsection{Selection of Trajectories}

Upon setting the basis for the approximation, it is important to
``select'' the data. Since the problem is now casted into an approximation
problem in $D\subset \R^d$, many approximation theories can be relied
upon. One of the key facts from approximation
theory is that, in order to construct accurate function approximation,
data should fill up the underlying domain in a systematic manner.
This implies that we should enforce the data pairings
\eqref{pairing} to ``spread'' out the domain-of-interest $D$.
Since the trajectories of many dynamical systems often evolve on
certain low-dimensional manifolds, using a single trajectory, no matter
how long it is, will unlikely to produce a set of data pairings for
accurate approximation, except for some special
underlying systems such as the chaotic systems
\cite{tran2017exact}. 
 This is the reason why we propose to use multiple number of
short trajectories. Specifically, the trajectories are
started from different initial conditions. The length of each
trajectories may be very short --- just long enough to allow one to
estimate the time derivatives. More specifically, we conduct the
following steps.

\begin{itemize}
\item Let $M\geq 1$ be the number of trajectories, and
$$
\x^{(m)}_0 = \u_0^{(m)} + {\bm \epsilon}_0^{(m)}, \qquad m=1,\dots, M,
$$
are different initial data, where ${\bm \epsilon}_0^{(m)}$ is the
observation error.

\item For each initial state, collect trajectory data generated by the
  underlying system \eqref{eq:ODE}
\begin{equation} \label{eq:data}
{\x}^{(m)}(t_{j}) = {\bf u}(t_{j};\u_0^{(m)}) + {\bm \epsilon}_{j}^{(m)}\in D,
\qquad 0\leq  j \le J_m, \quad 1\leq m \leq M,
\end{equation}
where $J_m$ is the number of data points collected for the $m$-th
trajectory and ${\bm \epsilon}_{j}^{(m)}$ is the observation noise.

\item Collect or estimate the time derivatives of the trajectories by using the data in \eqref{eq:data}
\begin{equation}
 {\dot \x}^{(m)} (t_{j}) =  {\dot \u}(t_{j}; \u_0^{(m)}) + \bm
 {\tau}^{(m)}(t_{j}), \qquad 
 j \in \Theta_m, \quad 1\leq m \leq M,
\end{equation}
where $\bm {\tau}^{(m)}(t_{j})$ stands for the estimation/observation
errors, and $\Theta_m=\{ j_1^{(m)},\cdots,j_{s_m}^{(m)} \}\subseteq\{ 0,\cdots,J_m  \}$. 

\item The pairings
\begin{equation}
\left\{ \x^{(m)}(t_{j}),  {\dot \x}^{(m)} (t_{j})\right\}, \qquad
j \in \Theta_m, \quad 1\leq m \leq M,
\end{equation}
are the data we use to approximation the equation, where
${\x}^{(m)}(t_{0}) = \x^{(m)}_0$ is the initial state.

\end{itemize}

We propose to use a large number of short bursts of trajectories,
i.e.,
\begin{equation} \label{traj}
J_m \sim {\mathcal O}(1), \quad m=1,\dots,M, \qquad M\gg 1.
\end{equation}
The length of each trajectory, $(J_m+1)$, is very small -- it should be merely long
enough to allow accurate estimation of the time derivatives. In case
the time derivatives are available by other means, e.g., direct
measurement, then $J_m$ can be set to 0.

There are several approaches to choose different initial states $\{{\bf u}_0^{(m)}\}$ from the region $D$. They include but are not limited to 
\begin{itemize}
	\item Sampling the initial states $\{{\bf u}_0^{(m)}\}$ from some random distribution defined on $D$, e.g., the uniform distribution. 
	\item Sampling the initial states $\{{\bf u}_0^{(m)}\}$ by quasi-Monte Carlo sequence, 
	such as the Sobol or Halton sequence, or other lattice rules.
	\item If $D$ is a hypercube, 
	the initial state set $\{{\bf u}_0^{(m)}\}$ can be taken as the set of quadrature points (e.g., the Gaussian points) or its random subset \cite{ZhouNarayanXiu_2015}. One can also sample the quadrature points 
	according some discrete measure \cite{WuShinXiu_SISC17}. 
\end{itemize} 
 

\subsection{Numerical Strategies} \label{sec:strategy}

We now discuss several key aspects of the numerical implementation of
the method.

\subsubsection{Computation of time derivatives}

In most of the realistic situations, 
only the trajectory data $\{{\bf x}^{(m)}(t_j)\}$  are available, 
and the time derivatives ${\bf \dot x}^{(m)}(t_j)$ need to be approximated numerically. 
Given a set of trajectory data  
$\{ {\bf x}^{(m)}(t_j), 0\leq j \le J_m  \}$, we seek to
numerically estimate the approximate time derivative ${\bf \dot
  x}^{(m)}(t_{j})$ at certain time instances $t_j$, $j \in
\Theta_m\subseteq \{0,\cdots,J_m\}$.  
This topic belongs to the problem of numerical differentiation of 
discrete data.  

If the data samples are noiseless, the velocity can be computed by a high-order 
finite difference method, e.g., 
a second-order central approximation with $J_m=2$ and $t_j=j\Delta t$ is
\begin{equation}\label{eq:finitediff}
{\bf \dot x}^{(m)}(t_{1})
= \frac{ {\bf \dot x}^{(m)}(t_{2}) -  {\bf \dot x}^{(m)}(t_{0})   }{2 \Delta t}.
\end{equation}
However, when data contain noise, estimation of 
derivatives can be challenging. 
In this case, if one uses  
the standard finite difference methods, then  
the resulting noise level in the first derivatives $ {\bf \dot x}^{(m)}(t_{j}) $ scales as 
$\sim \frac{ {\bm \epsilon}_{j}^{(m)} }{\Delta t}$.  
(See \cite{wagner2015regularised} for more analysis.)
Several approaches were developed for numerical differentiation of noisy data. They include 
least squares approximation methods \cite{knowles2014methods,wagner2015regularised},  
Tikhonov regularization \cite{doi:10.1137/0708026}, total variation regularization \cite{chartrand2011numerical},  
Knowles and Wallace variational method \cite{knowles1995variational}, etc. 
For an overview, see \cite{knowles2014methods}.    

In this paper, we use the least squares approach 
to denoise the derivative ${\bf \dot x}^{(m)}(t_{j})$. 
The approach consists of constructing a least squares
approximation of the trajectory, followed by an analytical
differentiation of the approximation to obtain the derivatives.
With the $m$-th burst of trajectory data, $m=1,\dots, M$,
we first construct a polynomial ${\bf p}(t)=\sum_{k=0}^{L} {\bf a}_k
\left(t- t_j\right)^k $ as an approximation 
to the trajectory ${\bf u}(t;{\bf u}_0^{(m)})$, by using the least
squares method with $1\le L \le J_m$. We then 
estimate the derivative as ${\bf \dot x}^{(m)}(t_{j})={\bf a}_1$. 


We remark that other methods for derivative computation with denoising
effect can be certainly used. We found the least squares method to be
fairly robust and easy to implement. When the noise level in the data
become too large to obtain accurate derivatives, our method would become unreliable. This is, however, a
universal challenge for all methods for equation recovery.

We also remark that even though the number of trajectories can be very
large, i.e., $M\gg 1$. Since each trajectory is very short with $J_m
\sim {\mathcal O}(1)$, $m=1,\dots, M$, the total amount of data for
the state  variables scales as $\mathcal{O}(M)$. In many cases this is
fairly small.

\subsubsection{Recovery algorithms: Regression methods}

The regression methods seek 
the expansion coefficients ${\bf c}_\ell$ by solving the following minimization problems
\begin{equation}\label{eq:LSf}
{\bf c}_\ell = \argmin_{ {\bf c} \in \mathbb{R}^{N} } \| {\bf A} {\bf
  c} - \dot{\mathbf{X}}_\ell \|,\quad \ell = 1,\dots, d, 
\end{equation}
where the norm $\| \cdot \|$ can be taken as the vector $\ell_2$-norm
$\| \cdot \|_2$, which corresponds to the well known least squares  method and can handle noisy data. 
It can also be taken as the vector $\ell_1$-norm, which is the least
absolute deviation (LAD) method and can eliminate potential corruption
errors in data (cf. \cite{candes2005decoding,ShinXiu}). If one would like to seek the sparse coefficients ${\bf c}_\ell$, 
the LASSO (least absolute shrinkage and selection operator) method
\cite{tibshirani1996regression} may be used,
\begin{equation}\label{eq:Lasso}
{\bf c}_\ell = \argmin_{ {\bf c} \in \mathbb{R}^{N} } 
\frac12 \| {\bf A} {\bf c} - \dot{\mathbf{X}}_\ell \|^2 + \lambda \| {\bf c} \|_1,\quad \ell = 1,\dots,d, 
\end{equation}
where $\lambda$ is a parameter to be chosen.
Since our method here seeks accurate approximation to the governing
equations, rather than exact recovery, we focus on the use of least
squares and LAD.

\begin{remark}
The methods discussed here are based on component-wise
  approximation of the state variable $\u$. That is, each component of
the states is approximated independently. It is possible to seek joint
approximation of all components of $\u$ simultaneously. A joint
approximation may offer computational advantages and yet introduce new
challenges. We shall leave this to a separate study.
\end{remark}

\subsubsection{Recovery algorithms: Matrix-free method}

When the data set is exceptionally large ($M\gg 1$), or unlimited as
in stream data, matrix-free method such as the
sequential approximation (SA)
method (cf. \cite{ShinXiu_SISC17, WuShinXiu_SISC17, WuXiu_JCP18, ShinWuXiu_JCP18}), can be more effective.

The SA method is of iterative nature. It aims at conducting
approximation using only one data point at each step. By 
starting from arbitrary choices of the initial approximation, which
corresponds to an arbitrary choice of the coefficient vectors
$\c_{\ell}^{(0)}$ at $k=0$ step, the iteration of the coefficient
vector at $k\geq 1$ step is
\begin{align} \label{eq:updatec}
&{\bf c}_{\ell}^{(k)}   = {\bf c}_{\ell}^{(k-1)}  + \frac{ { {\dot{x}}}_{\ell}^{(k)} ( t_{j_k}) -  \langle  {\bf
		c}^{(k-1)}_{\ell},{\bf \Phi} ( {\bf x}^{(k)} ( t_{j_k} ) ) \rangle  } { \left\|  {\bf \Phi} ( {\bf x}^{(k)} ( t_{j_k} ) ) \right\|_2^2 + \gamma_k }   {\bf \Phi} ( {\bf x}^{(k)} ( t_{j_k} ) ),\quad k=1,\dots,
\end{align}
where $\gamma_k\geq 0$ are parameters chosen according to the noise level in the $k$-th piece of data  $\big(   {\bf x}^{(k)} ( t_{j_k} )  , {\dot { \bf x}}^{(k)} ( t_{j_k} )   \big)$, 
and $ \dot{x}_{\ell} $ denotes the $\ell$-th component of $\dot{\bf
  x}$. For noiseless data, one can set $\gamma_k = 0$.

The scheme is derived by minimizing a cost function that consists of
data mismatch at the current arriving data point and the difference in the polynomial approximation
between two consecutive steps. The method is particularly suitable for
stream data arriving sequentially in time, as it uses only one data
point at each step and does not require storage or knowledge of the
entire data set. Moroever, the scheme uses only vector operations and
does not involve any matrices. For more details of the method,
including convergence and error analysis, see \cite{ShinXiu_SISC17,
  WuShinXiu_SISC17, WuXiu_JCP18, ShinWuXiu_JCP18}.

\subsection{Error analysis}

Assuming the aforementioned numerical algorithms can construct a
sufficiently accurate approximation of ${\bf f}$ with a
error bound
$\| \tilde{\bf f} - {\bf f} \|_{2,L^2_\omega}$, we can 
estimate the error in solution of the recovered system \eqref{eq:appODE}.

Let ${\bf x}(t)$ be the exact solution of the recovered system
\eqref{eq:appODE} and ${\bf u}(t)$ 
the exact solutions of true system \eqref{eq:ODE}.
Assume that ${\bf x}(t_0)={\bf u}(t_0) \in D$,
we then have the following error bound for the difference between ${\bf x}(t)$ and ${\bf u}(t)$.

\begin{theorem}\label{thm:errorx2}
Assume $\bf f$ is Lipschitz continuous on $D$,
$$
\| {\bf f} ({\bf u}_1 ) - {\bf f} (  {\bf u}_2 ) \|_2 \le L_f \| {\bf u}_1 -  {\bf u}_2  \|_2,\quad \forall {\bf u}_1, {\bf u}_2 \in D.
$$
If ${\bf x} (t),{\bf u}(t) \in D$, for $t_0\le t \le T$, then
\begin{equation}\label{eq:L2Bx}
\big\|  {\bf x} (t)- {\bf u}(t) \big\|_2 \le (t-t_0)  \| \tilde{\bf f} - {\bf f} \|_{2,L^\infty} + L_f \int_{t_0}^t \big\|  {\bf x} (s) - {\bf u}(s) \big\|_2 ds,
\end{equation}
Furthermore
\begin{equation} \label{eq:L2BxFinal}
  \big\|  {\bf x} (t) - {\bf u}(t) \big\|_2  \le
 L_f^{-1} \big( {\rm e}^{L_f ( t-t_0) } - 1 \big) 
\Big( \|{\bf f}-{\mathcal P}_V {\bf f}\|_{2,L^\infty} +    \| \tilde {\bf f}  - {\mathcal P}_V  {\bf f}  \|_{2,L^2_\omega} \| \sqrt{K}\|_{L^\infty} \Big),
\end{equation}
where the operator ${\mathcal P}_V$ denotes the component-wise projection on to $V$, 
and $K(\x)$ denotes  the ``diagonal'' of the
reproducing kernel of $V$. With any orthogonal basis $\{\psi_j(\x)\}_{j=1}^N$ of $V$, the kernel $K$ can be expressed as $K(\x) = \sum_{j=1}^{N} \psi_j^2(\x).$
\end{theorem}

The proof can be found in \ref{app:bound}. We remark that the error
bound here is rather crude, as it addresses a general system of
equation.  More refined error bounds can be derived
when one has more knowledge of the underlying system of equations.

\section{Numerical Examples} \label{sec:examples}

In this section we present numerical examples to verify 
the performance 
of the proposed methods for recovering the ordinary differential equations (ODEs) 
and semi-explicit differential algebraic equations (DAEs).


In all the test cases, 
we assume that there is no prior knowledge about the governing equations.  
Without loss of generality, we assume that $t_j =j \Delta t$. 
Unless otherwise stated, 
$\Delta t$ is taken as $0.005$, the initial states $\{ {\bf u}_0^{(m)}\}$ 
are sampled by the uniform distribution over the region $D$, 
and the basis is taken as $\{ {\bf x}^{\bf i}, |{\bf i}| \le n \}$.  
We employ a variety of test cases of both linear and nonlinear
systems of ODEs. The examples include relatively simple and well understood ODE systems
from the textbook
\cite{boyce2009elementary}, as well as more complex nonlinear systems
of coupled ODE and DAEs. Moreover, we also introduce corruption errors
in the trajectory data, in addition to the standard random noises. We
examine not only the errors in approximating the exact equations, but
also the errors in the predictions made by the recovered equations.

Our primary focus is to demonstrate the flexibility of using
polynomial approximation for accurate equation approximation, rather
than exact recovery, as well as the effect of using a large number of
short bursts of trajectories.
For the actual computation of the approximation, we primarily use the
standard least squares method. When data corruption exists, we also
use LAD, i.e., $\ell_1$ regression, and employed
the software package $\ell_1$-Magic \cite{candes2005l1}. We also
include an example with a very large data set and use the sequential
method to demonstrate the flexibility of the matrix-free method.

\subsection{Linear ODEs}

We first present a set of numerical results for simple linear ODEs with $d=2$, to 
examine the performance of the methods in recovering the 
structure of ODEs around different types of critical points.  All
examples are textbook problems from \cite{boyce2009elementary}.

\begin{example}\label{ex:linearODE}  \rm 
Six linear ODE systems are considered here, and their structures 
cover all the elementary critical points.  
The description of the systems and the corresponding numerical results are displayed in Table \ref{tab:1}. 
In all the tests, $n=1$ is used. 
The first and second tests use $M=10$ sets of noiseless  trajectory data. The third and forth tests use  $M=10$ sets of noisy trajectory data, where the noise follows uniform distribution in $[-0.01,0.01]$. 
The last two tests considers $M=20$ sequences of
trajectory data, with two uncertain (randomly chosen) 
sequences ($10\%$ of the total sequences) have 
large corruption errors following the normal distribution $N(0.5,1)$. 
It is seen from Table \ref{tab:1} that 
the proposed methods are able to 
recover the correct structure of 
the ODE systems from few measurement data. 
The values of approximate coefficients in the learned governing equations 
are very close to the true values.

 We now discuss the need for using a large number of short bursts of trajectory data.
	Let us consider the
	fourth linear ODE system in Table \ref{tab:1}, which is the nodal sink system
	\begin{equation}\label{eq:linearODE4}
	\begin{cases}
	\dot u_1 =-2u_1+u_2-2,\\
	\dot u_2=u_1-2u_2+1.
	\end{cases}
	\end{equation}
Assume that the observed data are from a single trajectory 
	$$
	{\bf x} (t_j;{\bf u}_0) = {\bf u} (t_j;{\bf u}_0),
	$$
	originated from the initial state ${\bf u}_0$. If ${\bf u}_0$ belongs to the 
	line $\Gamma=\{ {\bf u}=(u_1,u_2): u_2-u_1-1=0 \}$, then all
        the data on the trajectory $\{{\bf x} (t_j;{\bf u}_0)\}$
        belong to $\Gamma$. It is straightforward to see that
        the solution of the following
        linear system 
	\begin{equation}\label{eq:linearODE4app}
	\begin{cases}
	\dot u_1 = -u_2 + \eta_1 ( u_2-u_1-1 ),\\
	\dot u_2= -u_2 + \eta_2 ( u_2-u_1-1 ),
	\end{cases}
	\end{equation}  
	can exactly fit the data on $\Gamma$, for arbitrary $\eta_1,
        \eta_2 \in \mathbb R $. 
       This implies that recovering \eqref{eq:linearODE4} becomes
       an ill-posed problem, as there are infinite number of linear
       systems \eqref{eq:linearODE4app} that correspond to such a
       single trajectory data.
From the point of view of function approximation, this is not
surprising, as a single trajectory of a given system often falls onto a
certain low dimensional manifold (unless the system is
chaotic). Approximating an unknown function using data from a
low-dimensional manifold is highly undesirable and can be
ill-posed. The use of a large number of short bursts of trajectory
effectively eliminates this potential difficulty. 

\end{example}

	\begin{table}[!p]\label{tab:1}
		\caption{Example \ref{ex:linearODE}: Local recovery of linear ODEs}
\begin{calstable}
	\alignC
	\makeatletter
	\colwidths{
		{0.16\textwidth}
		{0.21\textwidth}
		{0.23\textwidth}
		{0.20\textwidth}
		{0.16\textwidth}		
	}

	\def\cals@framers@width{.4pt }

	\thead{
		\brow
		\cell{\tt Phase portrait} 
		\cell{\tt True system}
		\cell{\tt Computational settings}
		\cell{\tt Learned system}
		\cell{\tt Learned portrait} 
		\erow
	}
	
	
	\brow
	\cell{\includegraphics[width=0.14\textwidth]{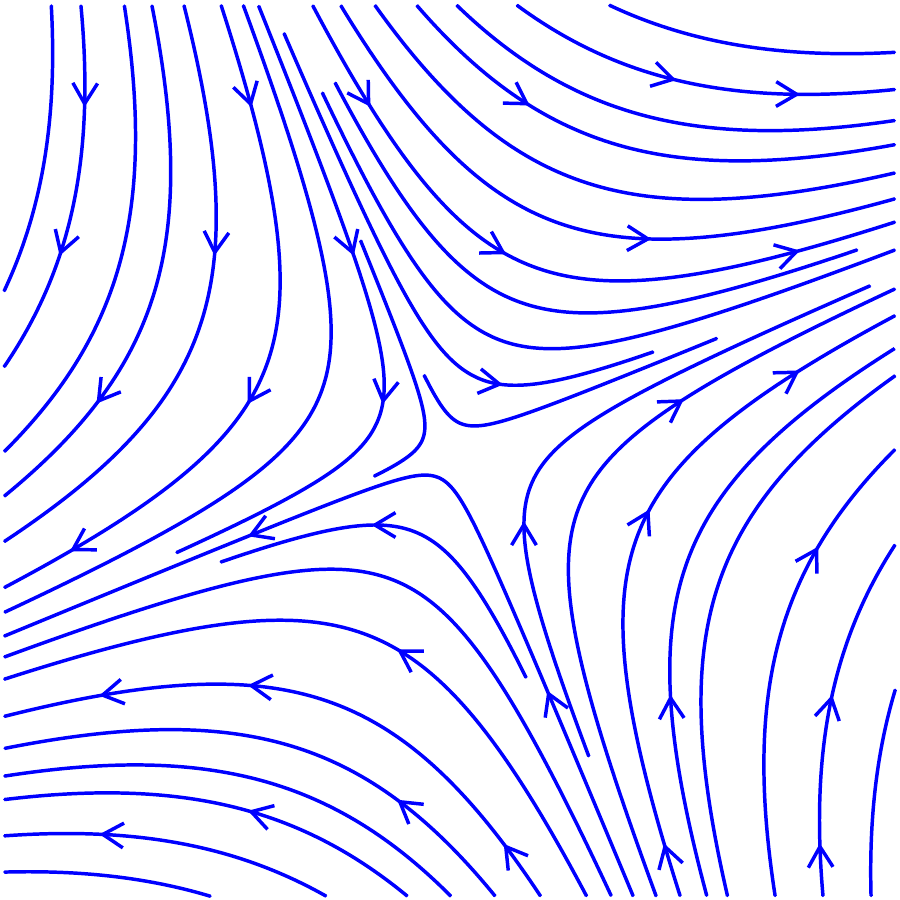}}
	\cell{ 
		\vspace{4mm}
		\scriptsize
			Saddle point
			\\[2mm]
			 $\begin{aligned}
				& \dot u_1 = u_1+u_2-2 \\ 
				& \dot u_2 = u_1 - u_2 \end{aligned}$
		}
	\cell{
		\vspace{4mm}
	\scriptsize
	 $D=[0,2]^2$
	\\[1mm]
	Noiseless data with $J_m=2$
	\\
	Least squares regression
	}
	\cell{ \scriptsize
	{  $\begin{aligned}
		 \dot x_1 & = 1.0000083x_1
		 \\
		 &\quad +1.0000003x_2
		 \\
		 &\quad 1.0000083 \\ 
		 \dot x_2 &= 1.0000083 x_1 
		 \\
		 & \quad -1.0000083 x_2 
		 \\
		 &\quad -1.65 \times 10^{-14} \end{aligned}$}
}	
	\cell{\includegraphics[width=0.14\textwidth]{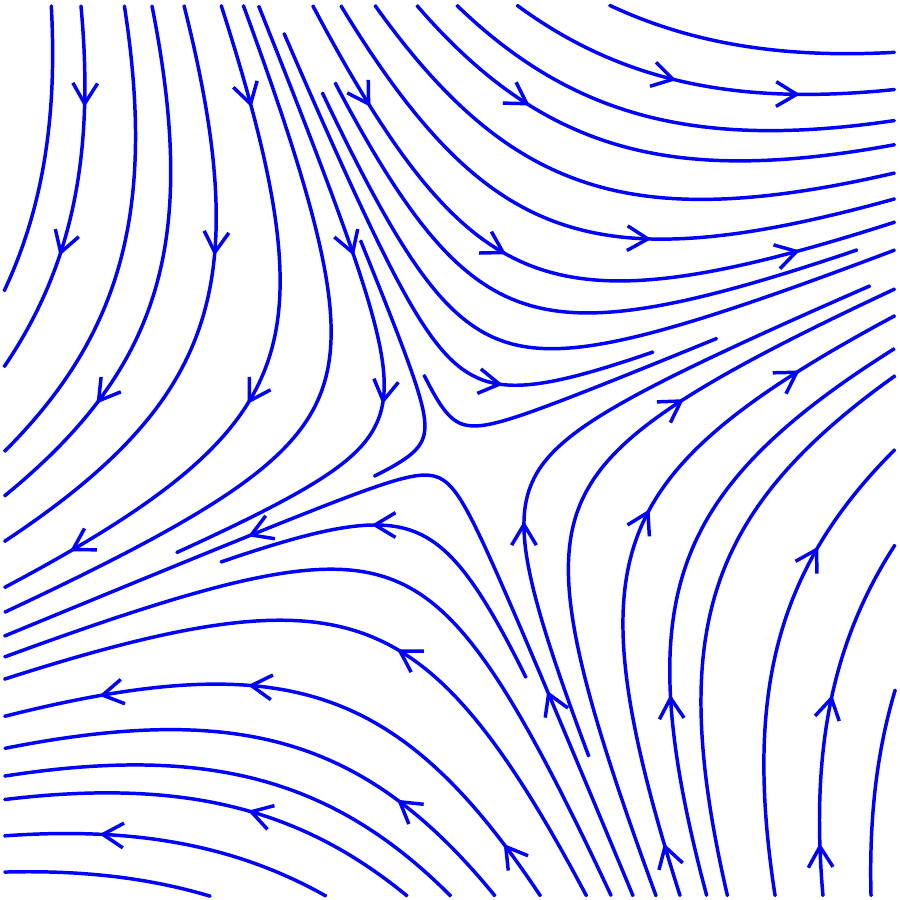}}
	\erow

\brow
\cell{\includegraphics[width=0.14\textwidth]{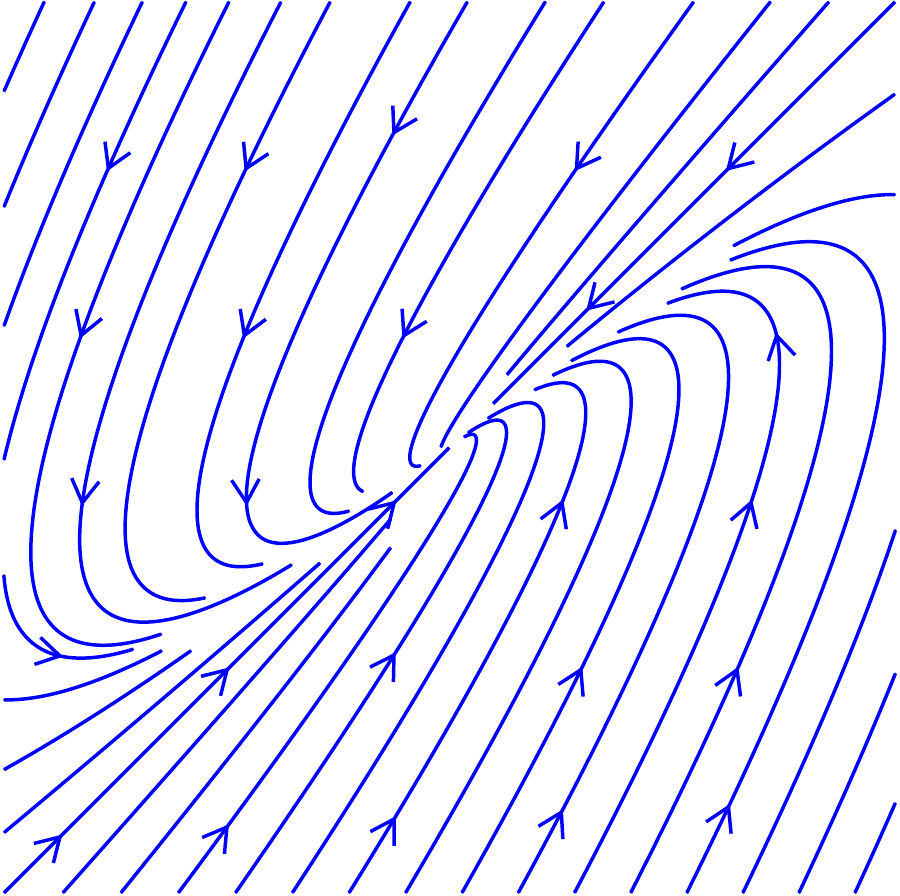}}
\cell{ 
	\vspace{4mm}
	\scriptsize
	Improper node
	\\[2mm]
	$\begin{aligned}
	&\dot u_1 = u_1-4u_2 \\ 
	& \dot u_2 = 4u_1 - 7u_2 \end{aligned}$
}
\cell{
	\vspace{4mm}
	\scriptsize
	$D=[-1,1]^2$
	\\[1mm]
	Noiseless data with $J_m=2$
	\\
	$\ell_1$-regression
}
\cell{ \scriptsize
	{  $\begin{aligned}
		\dot x_1 & = 1.00034x_1
		\\
		&\quad -4.00045x_2
		\\
		&\quad +4.05 \times 10^{-16} \\ 
		\dot x_2 &= 4.00045 x_1 
		\\
		& \quad -7.00056 x_2 
		\\
		&\quad + 1.50 \times 10^{-15} \end{aligned}$}
}	
\cell{\includegraphics[width=0.14\textwidth]{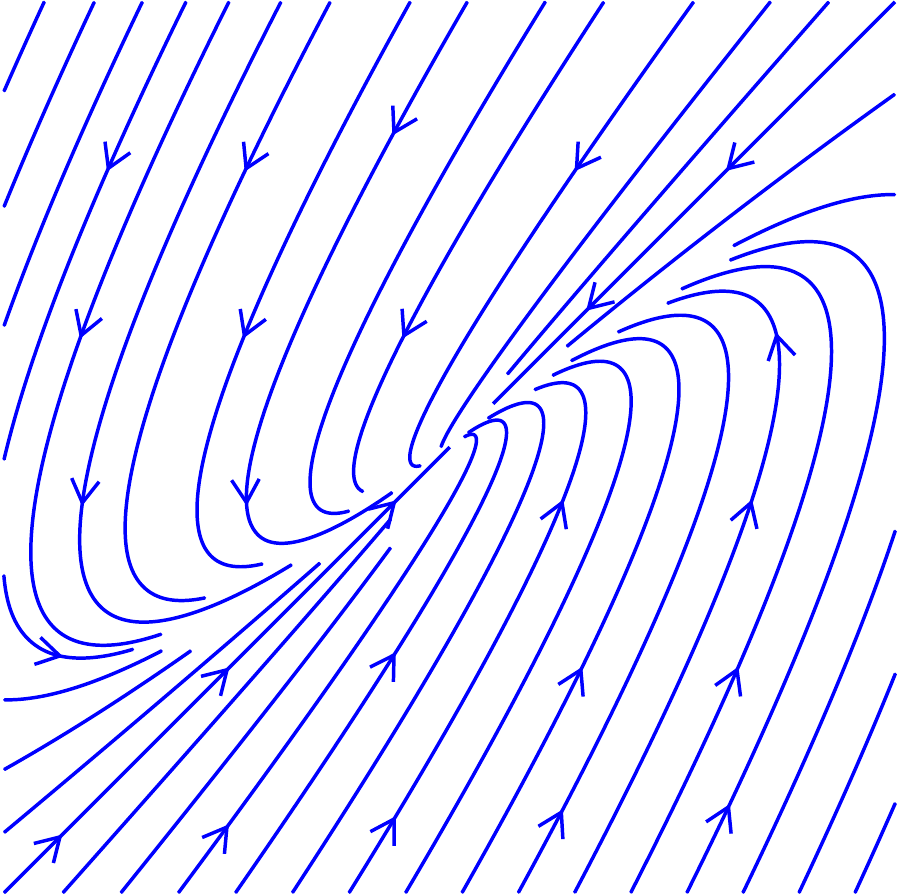}}
\erow

\brow
\cell{\includegraphics[width=0.14\textwidth]{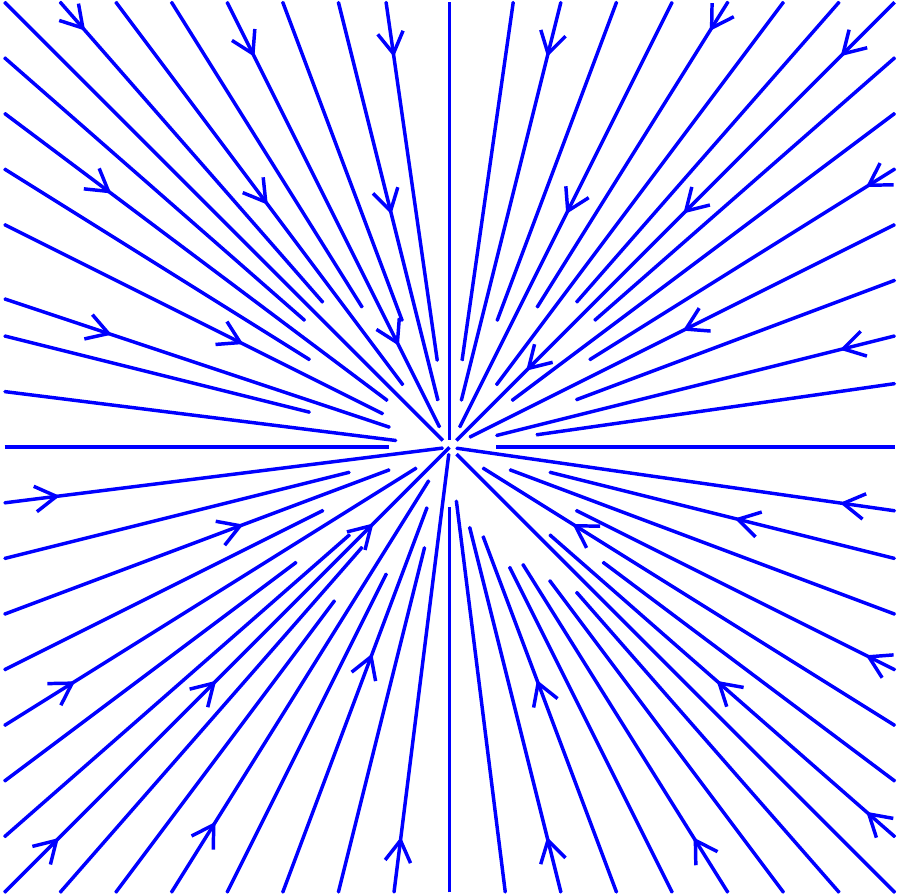}}
\cell{ 
	\vspace{4mm}
	\scriptsize
	Star point
	\\[2mm]
	$\begin{gathered}
	\dot u_1 = -u_1 \\ 
	\dot u_2 = -u_2 \end{gathered}$
}
\cell{
	\vspace{4mm}
	\scriptsize
	$D=[-1,1]^2$
	\\[1mm]
	Noisy data with $J_m=19$
	\\
	Least squares regression
}
\cell{ \scriptsize
	{  $\begin{aligned}
		\dot x_1 & = -1.0265x_1
		\\
		&\quad +0.0131x_2
		\\
		&\quad -0.0146 \\ 
		\dot x_2 &= -0.0392 x_1 
		\\
		& \quad -1.000038 x_2 
		\\
		&\quad + 0.0036 \end{aligned}$}
}	
\cell{\includegraphics[width=0.14\textwidth]{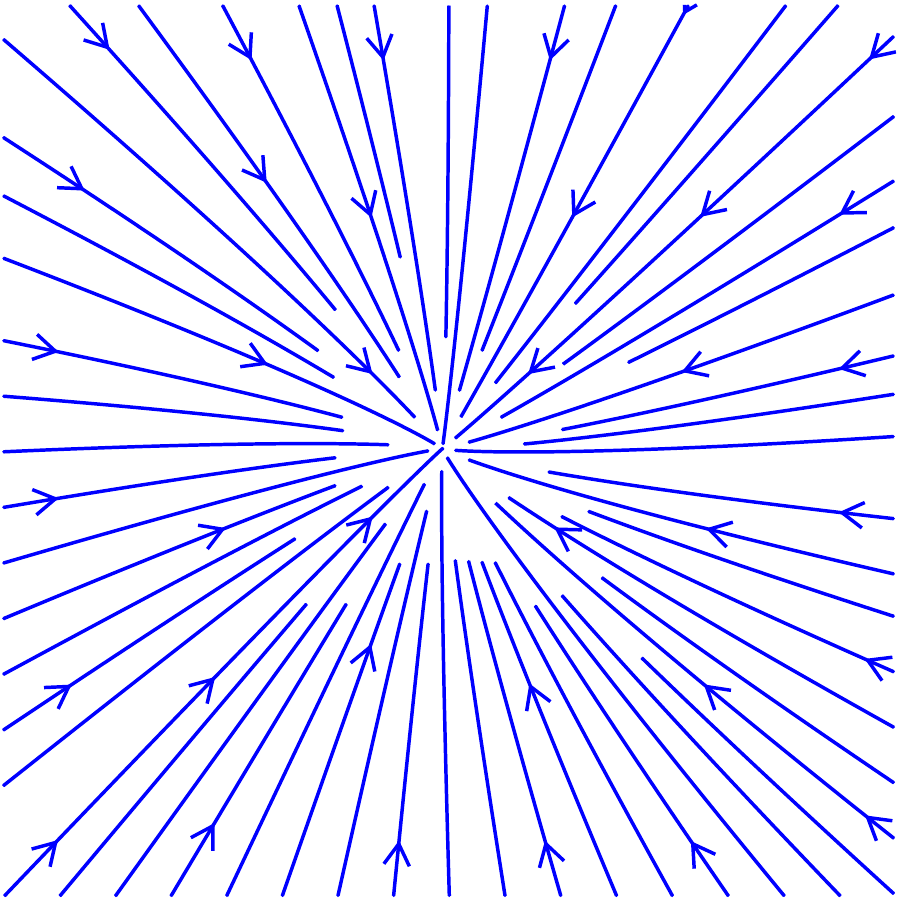}}
\erow

\brow
\cell{\includegraphics[width=0.14\textwidth]{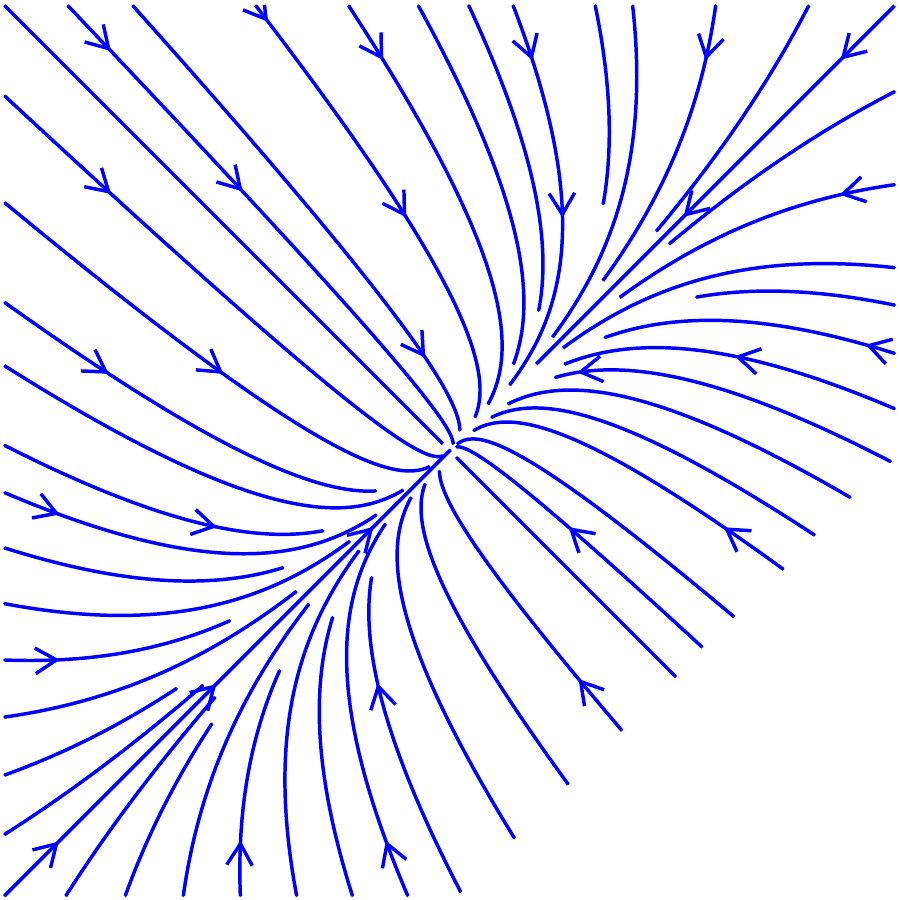}}
\cell{ 
	\scriptsize
	\vspace{4mm}
	Nodal sink
	\\[2mm]
	$\begin{gathered}
	\dot u_1= -2 u_1 +u_2-2 \\ 
	\dot u_2 = u_1 - 2u_2+1 
	 \\
	 u_2 > u_1 \end{gathered}$
}
\cell{
	\scriptsize
	\vspace{4mm}
	{\tiny$D=[-2,0]\times[-1,1] \cap  \{ {\bf x}| x_2>x_1 \}$}
	\\ 
	Noisy data with $J_m=19$
	\\
	Least squares regression
}
\cell{ \scriptsize
	{  $\begin{aligned}
		\dot x_1 & = -1.9556x_1
		\\
		&\quad 0.9693 x_2
		\\
		&\quad -1.9638 \\ 
		\dot x_2 &= 1.0061 x_1 
		\\
		& \quad -1.9920 x_2  
		\\
		&\quad + 1.0011 \end{aligned}$}
}	
\cell{\includegraphics[width=0.14\textwidth]{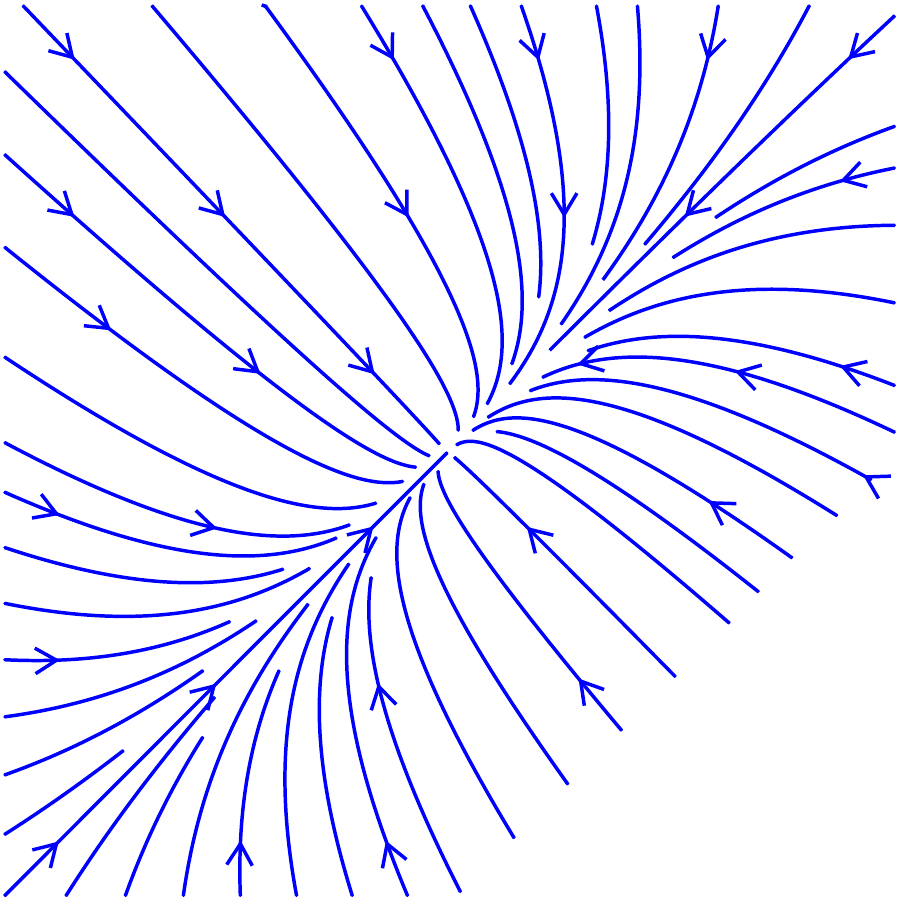}}
\erow

\brow
\cell{\includegraphics[width=0.14\textwidth]{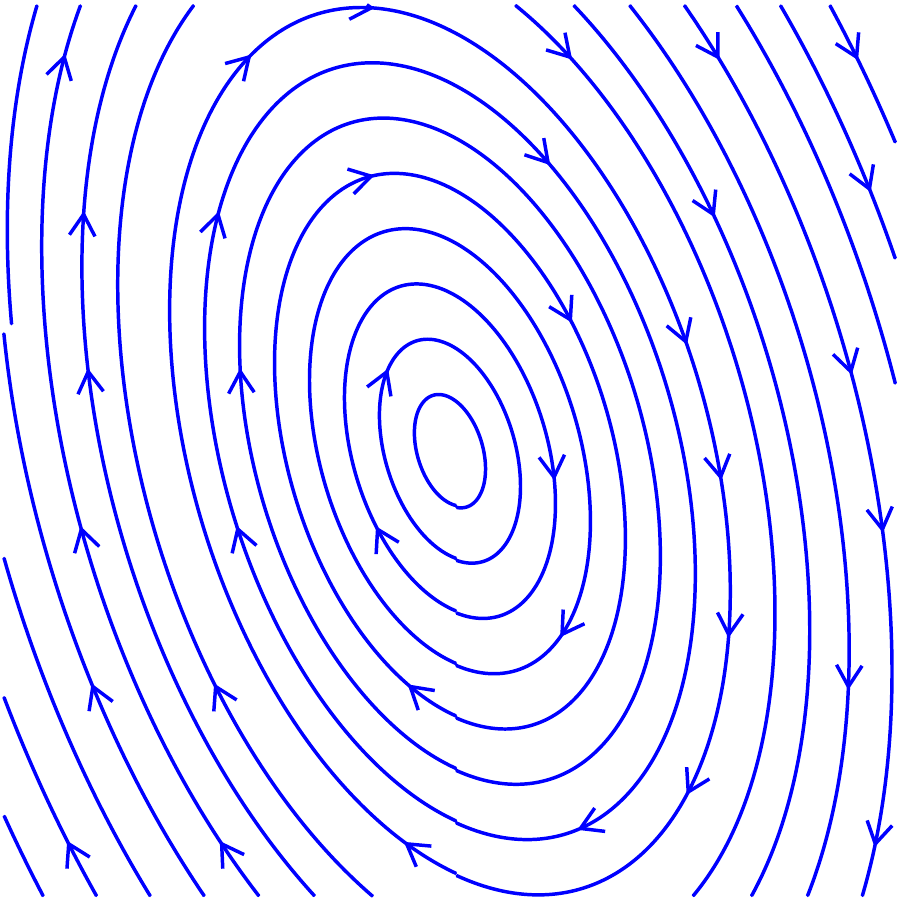}}
\cell{ 
	\vspace{4mm}
	\scriptsize
	Center point
	\\[2mm]
	$\begin{aligned}
	&\dot u_1 = u_1+2u_2 \\ 
	&\dot u_2 = -5u_1 - u_2 \end{aligned}$
}
\cell{
	\vspace{4mm}
	\scriptsize
	$D=[-1,1]^2$
	\\[1mm]
	Corrupted data with $J_m=2$
	\\
	$\ell_1$-regression
}
\cell{ \scriptsize
	{  $\begin{aligned}
		\dot x_1 & = 0.9999625x_1
		\\
		&\quad +1.999925x_2
		\\
		&\quad + 6.90 \times 10^{-14} \\ 
		\dot x_2 &=  -4.9998125 x_1 
		\\
		& \quad -0.9999625 x_2 
		\\
		&\quad -4.10 \times 10^{-13} \end{aligned}$}
}	
\cell{\includegraphics[width=0.14\textwidth]{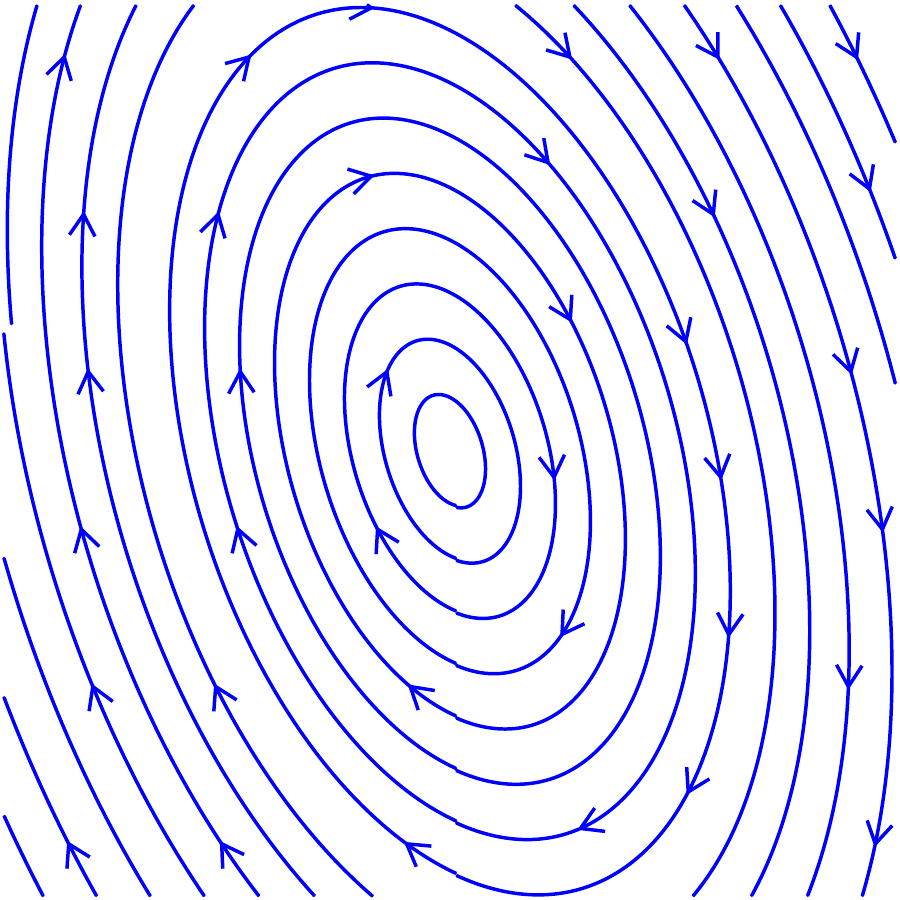}}
\erow

\brow
\cell{\includegraphics[width=0.14\textwidth]{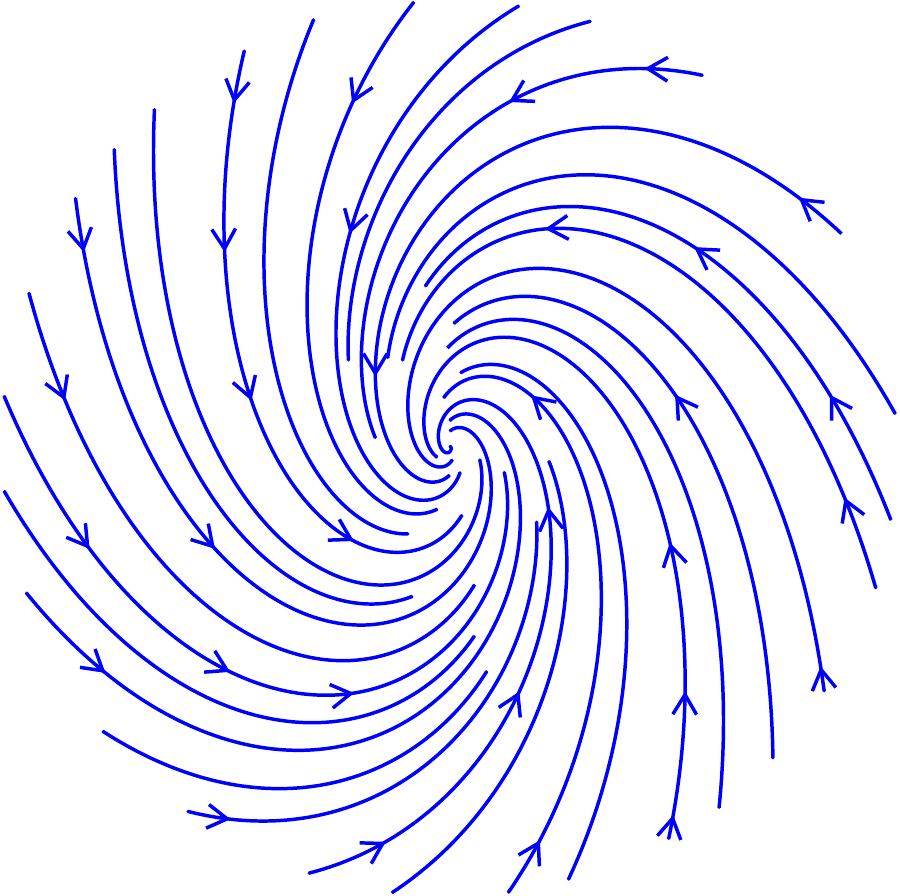}}
\cell{ 
	\scriptsize
	\vspace{4mm}
	Spiral point
	\\[2mm]
	$\begin{gathered}
	\dot u_1 = -u_1-u_2-1 \\ 
	\dot u_2 = 2u_1 - u_2+5
	\\
	 (u_1+2)^2+(u_2-1)^2 \le 1 \end{gathered}$
}
\cell{
	\vspace{4mm}
	\scriptsize
	{\tiny$D= \{ {\bf x}| (x_1+2)^2+(x_2-1)^2 \le 1 \}$}
	\\[1mm]
	Corrupted data with $J_m=2$
	\\
	$\ell_1$-regression
}
\cell{ \scriptsize
	{  $\begin{aligned}
		\dot x_1 & = -0.999979x_1
		\\
		&\quad -1.000004x_2
		\\
		&\quad -0.999954 \\ 
		\dot x_2 &= 2.0000083 x_1 
		\\
		& \quad -0.99997917 x_2 
		\\
		&\quad + 4.9999958 \end{aligned}$}
}	
\cell{\includegraphics[width=0.14\textwidth]{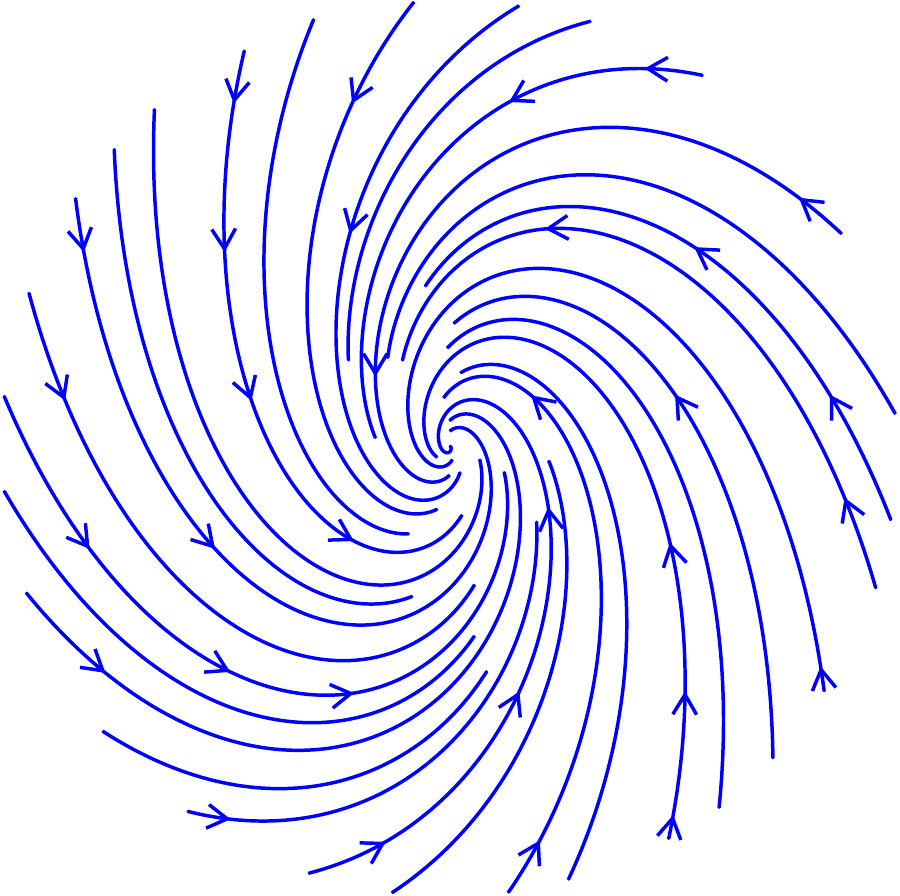}}
\erow

\end{calstable}

\end{table}


\subsection{Nonlinear ODEs}

We then test the proposed methods for recovering nonlinear ODE systems
exhibiting more complex dynamics. 
In general, the nonlinear ODEs do not have explicit analytical solutions, 
we therefore generate the measurement data by numerically solving the ODEs using the forth-order explicit  Runge-Kutta method 
with a small time step-size.

\begin{example}\label{ex:nonODE1} \rm
We first consider
the undamped Duffing equation
$$
\ddot u+u+\epsilon u^3=0,
$$
where $\epsilon$ is a small positive number and specified as $ 10^{-4}$. 
Taking $u_1=u$ and $u_2=\dot u$ gives
\begin{equation}\label{eq:nonODE1}
\begin{cases}
\dot u_1=u_2,\\
\dot u_2=-u_1-\epsilon u_1^3. 
\end{cases}
\end{equation}
The region $D$ is specified as $[0,2]\times[-1,1]$. 
From $M=30$ bursts of noiseless trajectory data with $J_m=2$, we try to learn the governing system using  
 different regression methods with $n=3$. The errors of approximate expansion coefficients in the learned systems 
 are displayed in Fig. \ref{fig:nonODE1}. 
As we can see, the coefficients are recovered with high accuracy by
$\ell_1$ and $\ell_2$ regressions. 

 \begin{figure}[htbp]
 	\centering
 	\subfigure[Errors in coefficients $\c_1$]{\includegraphics[width=0.47\textwidth]{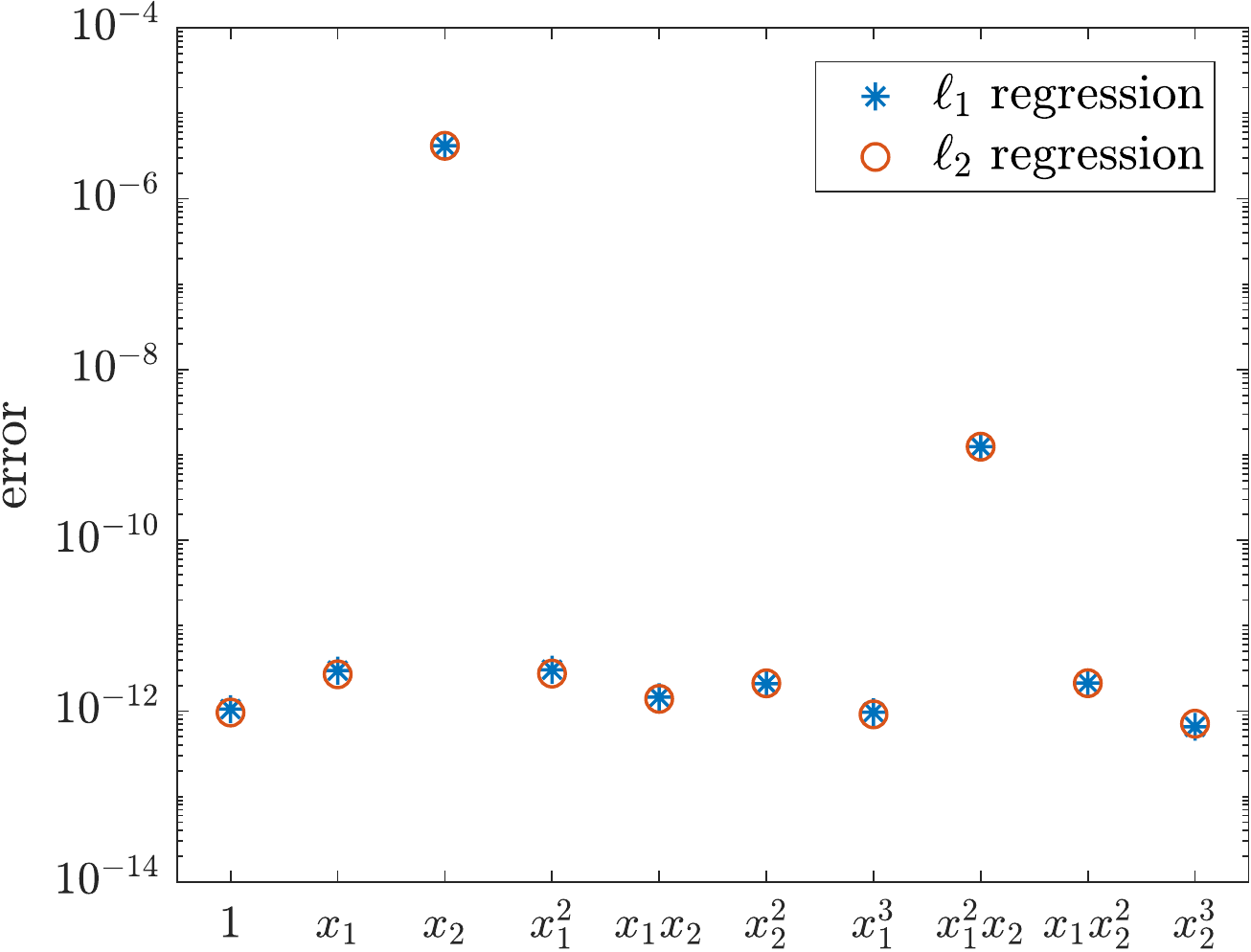}}~~~~
 	\subfigure[Errors in coefficients $\c_2$]{\includegraphics[width=0.47\textwidth]{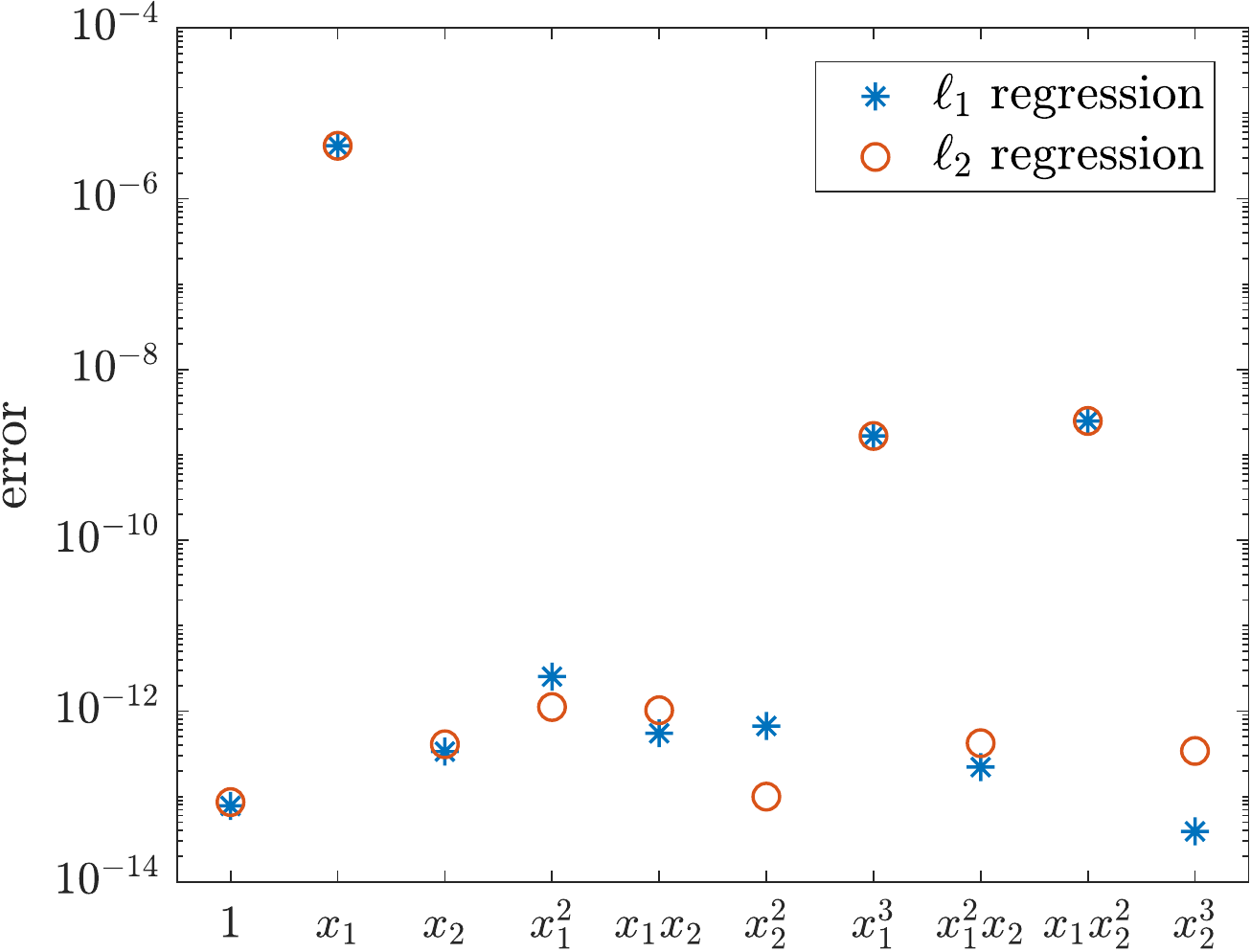}}
 	\caption{\small
 		Example \ref{ex:nonODE1}: Results for recovering the system \eqref{eq:nonODE1}. 
 		The symbols ``$\ast$'' and ``$\circ$'' represent the results of $\ell_1$ and $\ell_2$ regression methods, respectively. 
 		The coefficient errors which are not shown in the figures are zero errors.   
 	}\label{fig:nonODE1}
 \end{figure}
\end{example}


\begin{example}\label{ex:nonODE2} \rm
We consider 
\begin{equation}\label{eq:nonODE2}
\begin{cases}
\dot u_1=u_1 (1-u_1 - u_2),\\
\dot u_2=u_2 ( 0.5 - 0.25 u_2 - 0.75 u_1 ),
\end{cases}
\end{equation}
whose qualitative behavior was studied in \cite{boyce2009elementary}. 
We take $D=[-1,2]\times[-0.5,3]$, which contains the four critical points, 
an unstable node $(0,0)$, two asymptotically stable nodes $(1,0)$ and $(0,2)$, and 
a saddle point $(0.5,0.5)$. Suppose we have $M=30$ sets of trajectory data with $J_m=30$, 
and all the data contain i.i.d. random noises following uniform distribution in $[-0.01,0.01]$. 
Learning from these data, the least squares regression method with $n=2$ produced
the following approximate system 
\begin{equation}\label{eq:nonODE2app}
\begin{cases}
\dot x_1=x_1 (1.0217-1.0166x_1 -1.0164 x_2) -0.01695 + 0.0350 x_2 -0.0128 x_2^2,\\
\dot x_2=x_2 ( 0.50914 - 0.2500 x_2 -0.7544 x_1 ) -0.0160 + 0.0086 x_1 - 0.0016x_1^2.
\end{cases}
\end{equation}
Although in this case the data contain noises at a relatively large level, 
the proposed method can still recover the coefficients with errors
less than $4\times 10^{-2}$. 
In Figure \ref{fig:nonODE2}, 
the phase portraits of the original and recovered systems 
are plotted to provide a qualitative comparison between 
the true and the learned dynamics. 
Good agreements are observed in the structures, and 
the four critical points are correctly detected. 
Figure \ref{fig:nonODE2val} shows the solutions of the original and recovered systems 
at the initial state $(1.25,1.75)$, as well as the errors in 
the learned solution.

 \begin{figure}[htbp]
	\centering
	\includegraphics[width=0.47\textwidth]{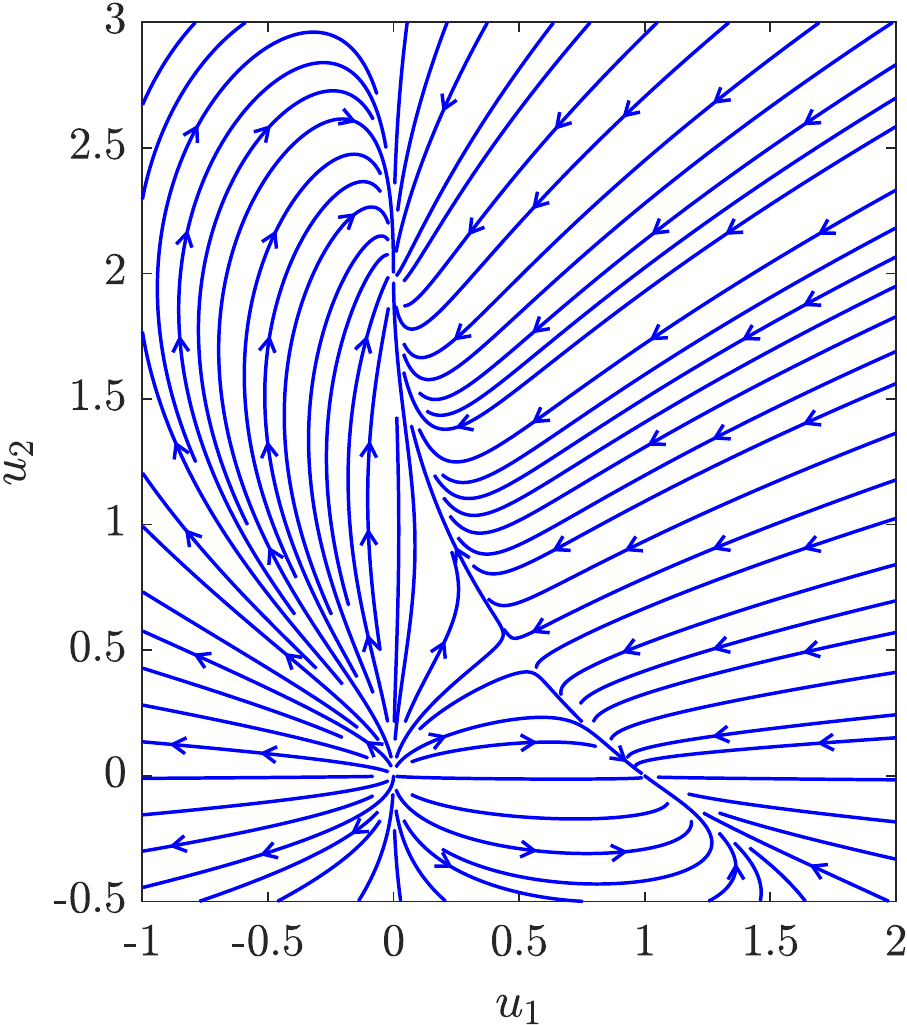}~~~~
	\includegraphics[width=0.47\textwidth]{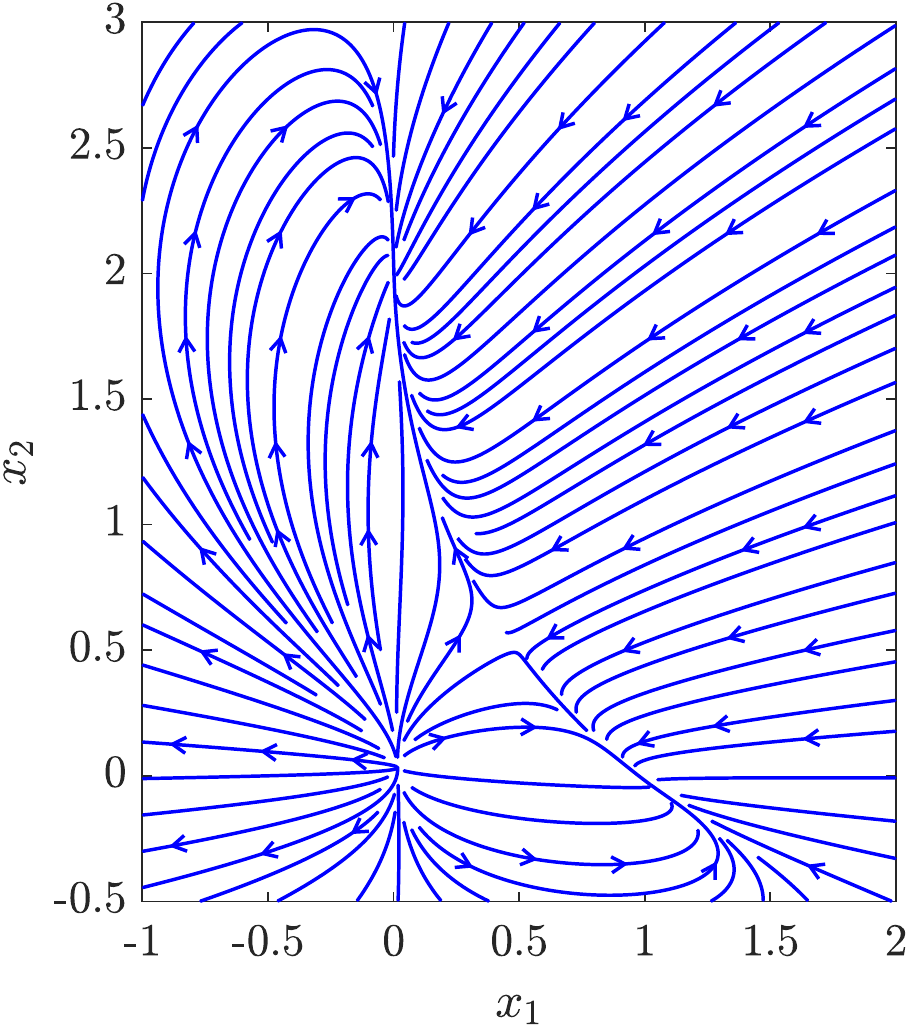}
	\caption{\small
		Example \ref{ex:nonODE2}: Phase portraits for the true  system \eqref{eq:nonODE2} (left) and the learned system \eqref{eq:nonODE2app} (right).   
	}\label{fig:nonODE2}
\end{figure}

 \begin{figure}[htbp]
	\centering
	\includegraphics[width=0.47\textwidth]{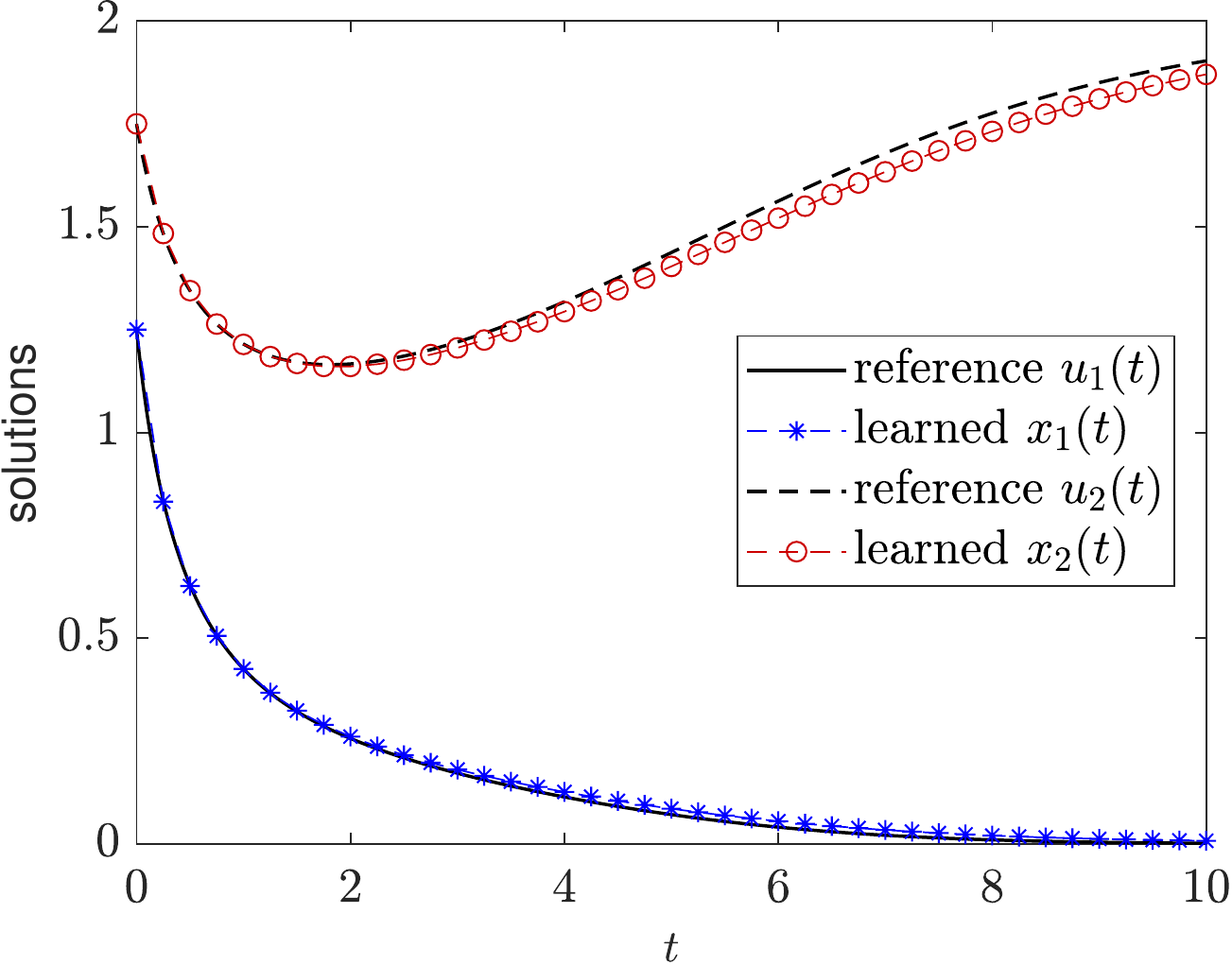}~~~~
	\includegraphics[width=0.47\textwidth]{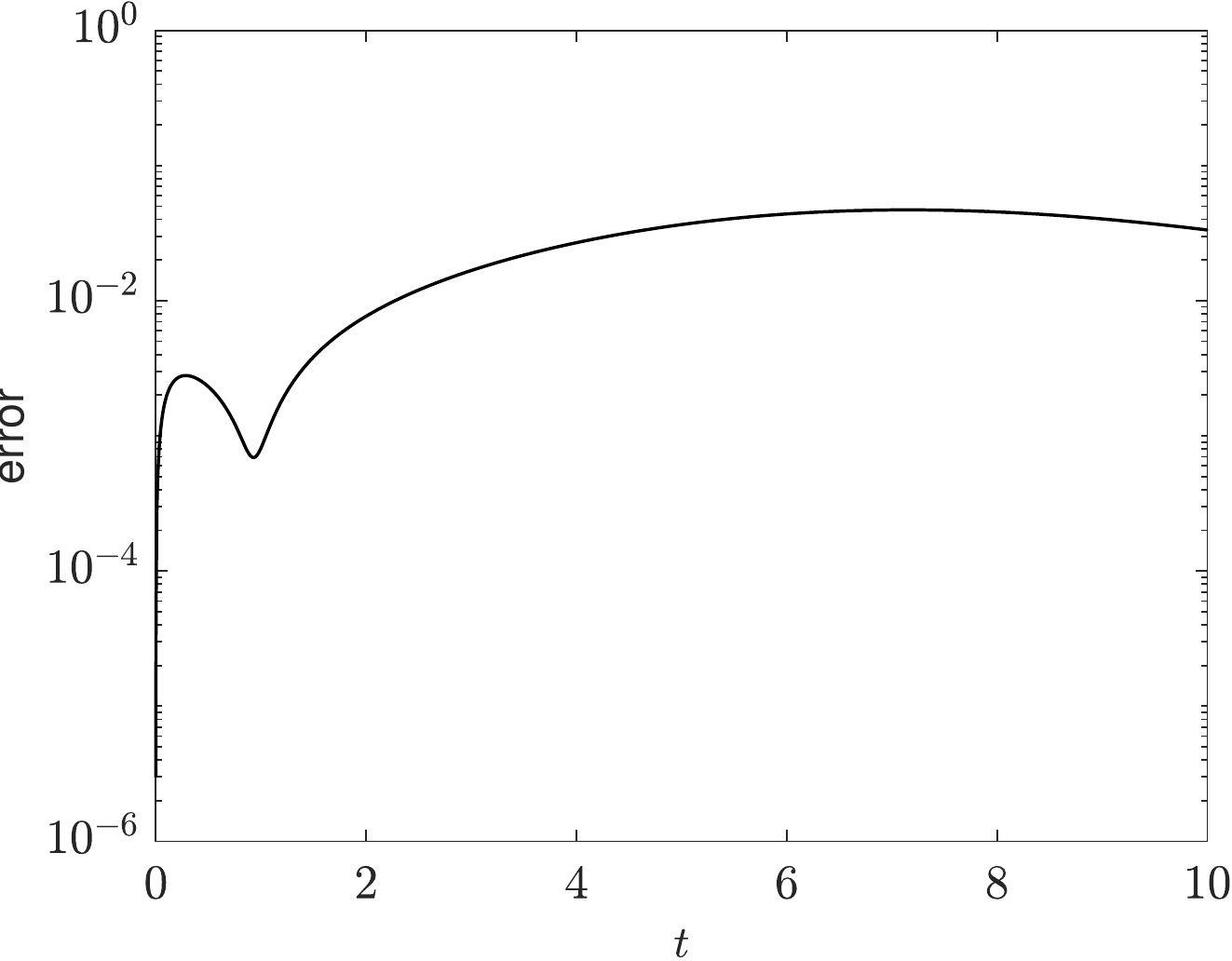}
	\caption{\small
		Example \ref{ex:nonODE2}: Validation of the recovered system \eqref{eq:nonODE2app} 
		with initial state $(1.25,1.75)$. Left: solutions. Right: error $\| {\bf x}(t) - {\bf u} (t) \|_2$. The horizontal axis denotes the time.   
	}\label{fig:nonODE2val}
\end{figure}

\end{example}

\begin{example}\label{ex:nonODE3} \rm
We now consider
\begin{equation}\label{eq:nonODE3}
\begin{cases}
\dot u_1 =u_2 -u_1 ( u_1^2 + u_2^2 - 1 ),\\
\dot u_2=-u_1 - u_2 ( u_1^2 + u_2^2 -1 ).
\end{cases}
\end{equation}
An important feature of this example is that 
we assume the state variables in the unit disk $D_b=\{ {\bf u}:
u_1^2+u_2^2 <1 \}$ are not accessible. Therefore we have data only
from outside this region. This rather arbitrary choice of
accessibility of data is to mimic the practical situation where data
collection may be subject to certain restrictions.
Also, we introduce relatively large corruption errors to the data.
Therefore, we have 
a nonstandard region $D=[-2,2]^2 \setminus D_b$, where data are
available and corrupted.

First, we use $M=40$ bursts of trajectory data with $J_m=2$
(very short burst consisting only 3 data points), 
where $10\%$ of the data contain relatively large corruption errors
following i.i.d. normal distribution $N(0.5,1)$. (Note the relatively
large bias in the mean value.)
We then use $\ell_1$-regression to learn the governing equations from these data, and obtain 
\begin{equation}\label{eq:nonODE3app}
\begin{cases}
\dot x_1 &= 0.9986 x_2 -x_1 ( 1.0020 x_1^2 + 1.0018 x_2^2 - 1.0047 )   -0.0006 \\
& \quad  + 0.0005 x_1^2 - 0.0001  x_1x_2 + 0.0002 x_2^2 
+ 0.0006 x_1^2 x_2 + 0.0005 x_2^3 ,\\
\dot x_2 &= -1.0025 x_1 -  x_2 ( 1.0015 x_1^2 + 1.0011 x_2^2 - 1.0025 ) -0.0005 \\
	& \quad - 0.00008 x_1 ^2 + 0.0004 x_1 x_2 + 0.0004 x_2^2 
	+ 0.0006 x_1^3 + 0.0005 x_1 x_2^2.
\end{cases}
\end{equation}
We see that the $\ell_1$-regression can effectively eliminate
the sparse corruption errors, 
and the coefficients 
are recovered with high accuracy. 
The phase portraits in Figure \ref{fig:nonODE3} 
further exhibit the accurate dynamics of the learned system  \eqref{eq:nonODE3app}. We also validate the solution of 
the recovered system  \eqref{eq:nonODE3app} 
with an arbitrarily chosen initial state $(-1.325,1.874)$, 
see Figure \ref{fig:nonODE3val}. 
The solution ${\bf x}(t)$ agrees well with the reference 
solution ${\bf u}(t)$ obtained from the exact equation. The increase
of the error over time is expected, as in any numerical solvers for
time dependent problems.

 \begin{figure}[htbp]
	\centering
	\includegraphics[width=0.47\textwidth]{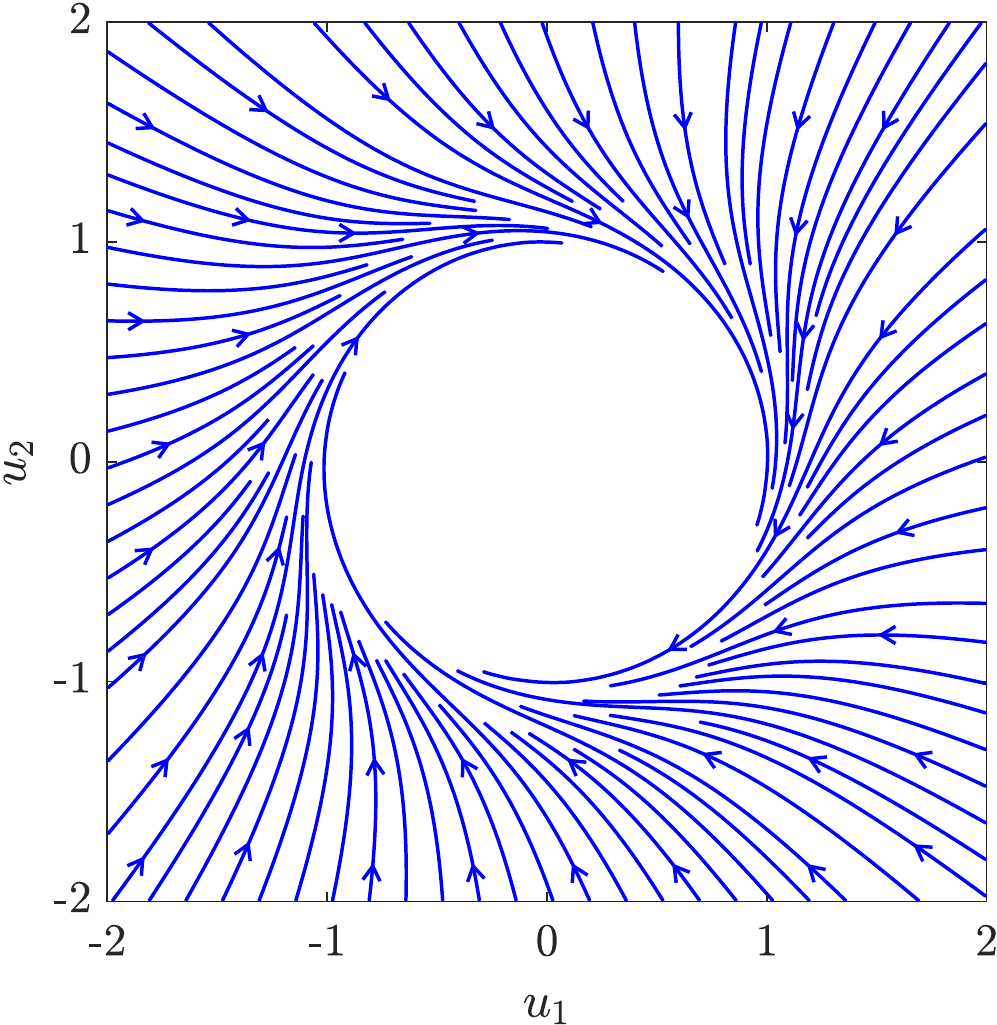}~~~~
	\includegraphics[width=0.47\textwidth]{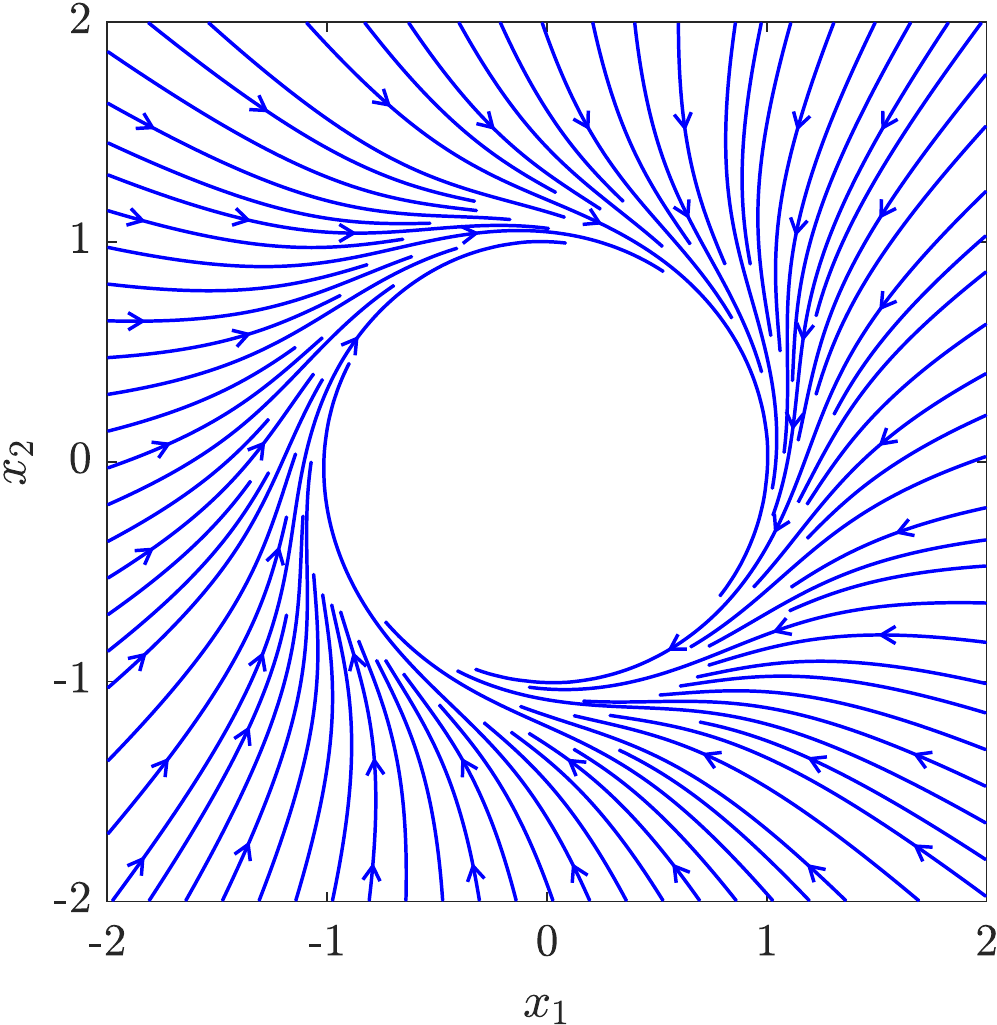}
	\caption{\small
		Example \ref{ex:nonODE3}: Phase portraits for the true  system \eqref{eq:nonODE3} (left) and the learned system \eqref{eq:nonODE3app} (right).   
	}\label{fig:nonODE3}
\end{figure}

 \begin{figure}[htbp]
	\centering
	\includegraphics[width=0.47\textwidth]{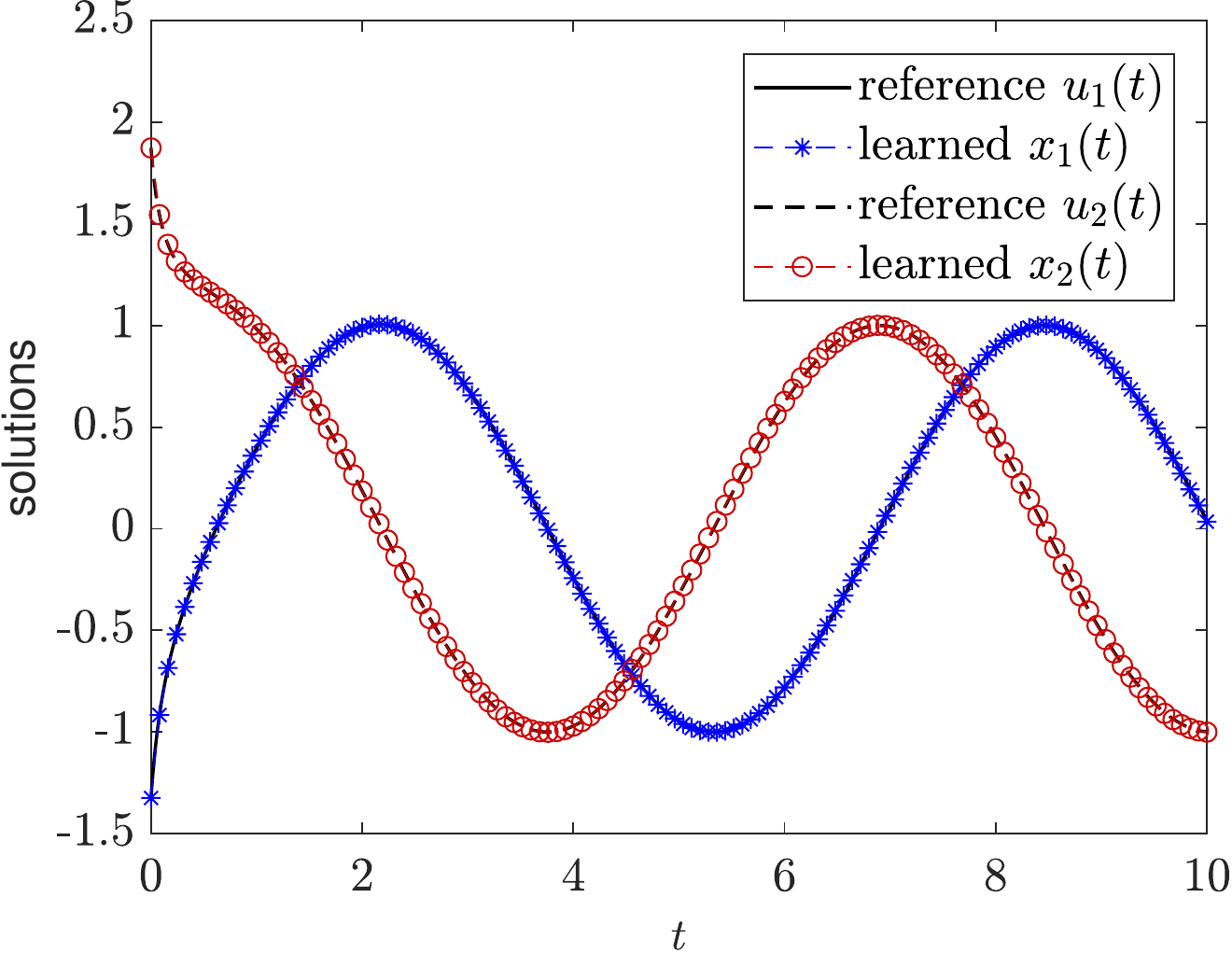}~~~~
	\includegraphics[width=0.47\textwidth]{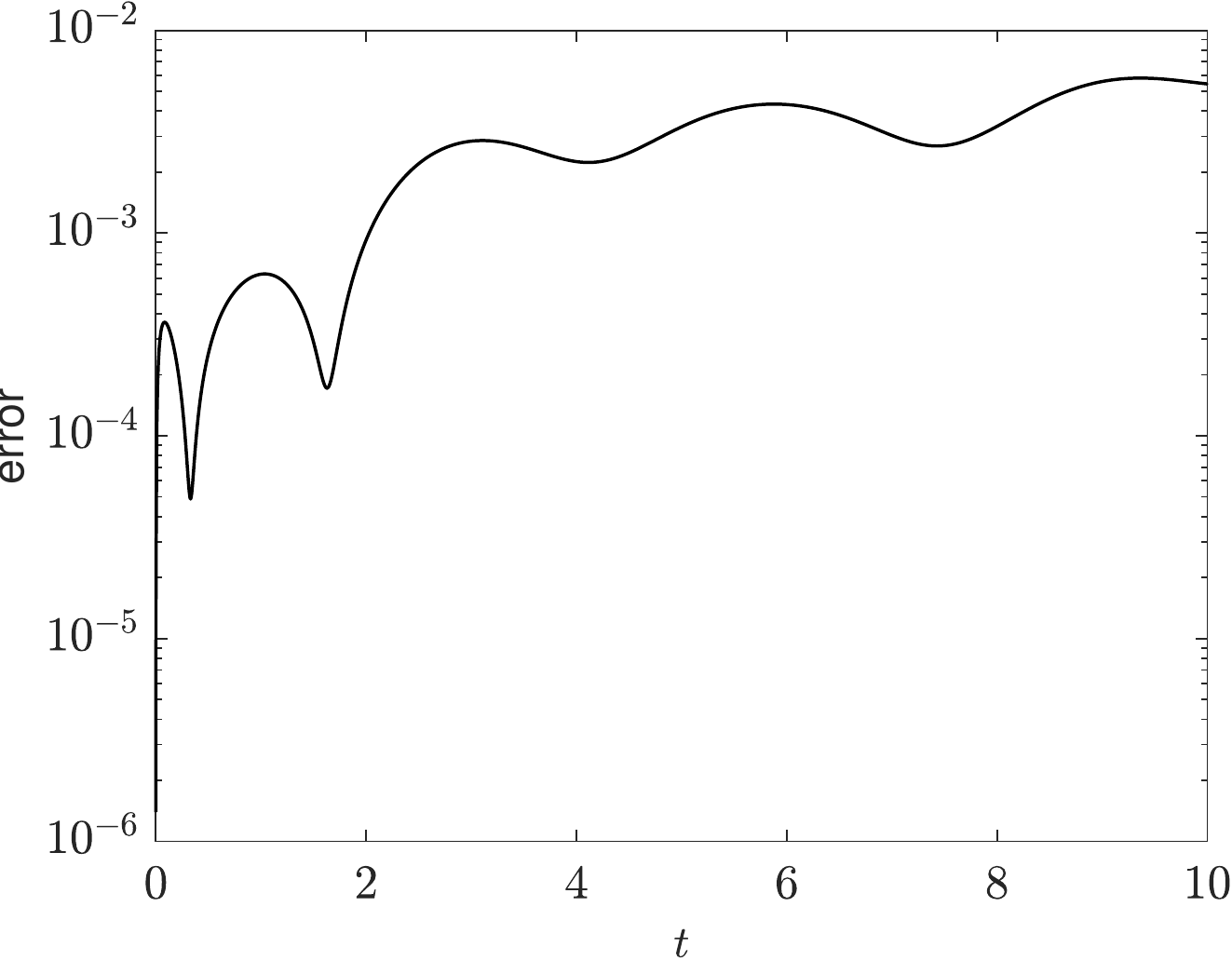}
	\caption{\small
		Example \ref{ex:nonODE3}: Validation of the recovered system \eqref{eq:nonODE3app} 
		with initial state $(-1.325,1.874)$. Left: solutions. Right: error $\| {\bf x}(t) -  {\bf u} (t) \|_2$. The horizontal axis denotes the time.   
	}\label{fig:nonODE3val}
\end{figure}

	We now consider data only with standard random noises and study the 
effect of noise level on the accuracy of recovered equations. 
Suppose we have $M=40$ sets of trajectory data with $J_m=10$, 
and all the data contain i.i.d. random noises following uniform
distribution in $[-\eta,\eta]$. Here the parameter $\eta$ represents
the noise level, and 
$\eta=0$ corresponds to noiseless data. 
	The errors of the expansion coefficients in the learned systems 
	are displayed in Fig. \ref{fig:nonODE3err}, for different
        noise levels at $\eta=0, 0.005, 0.01, 0.05, 0.1, 0.2$. 
	 The vector 2-norms of the coefficient errors are shown in
         Fig.~\ref{fig:nonODE3err2}.
It is clearly seen that the errors become larger at larger noise level
in the data, which is not surprising.
	
	 \begin{figure}[htbp]
		\centering
		\subfigure[Errors in coefficients $\c_1$]{\includegraphics[width=0.47\textwidth]{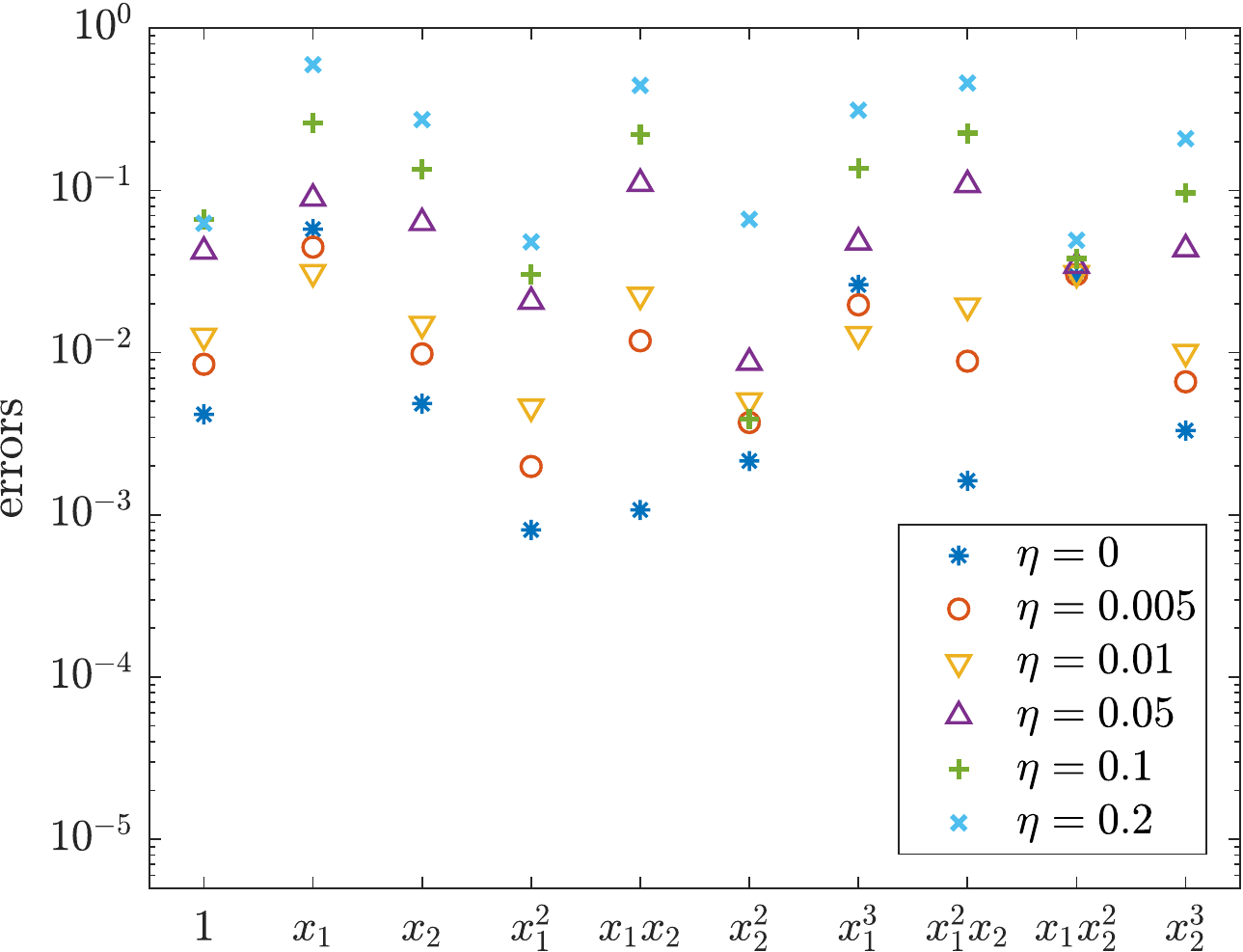}}~~~~
		\subfigure[Errors in coefficients $\c_2$]{\includegraphics[width=0.47\textwidth]{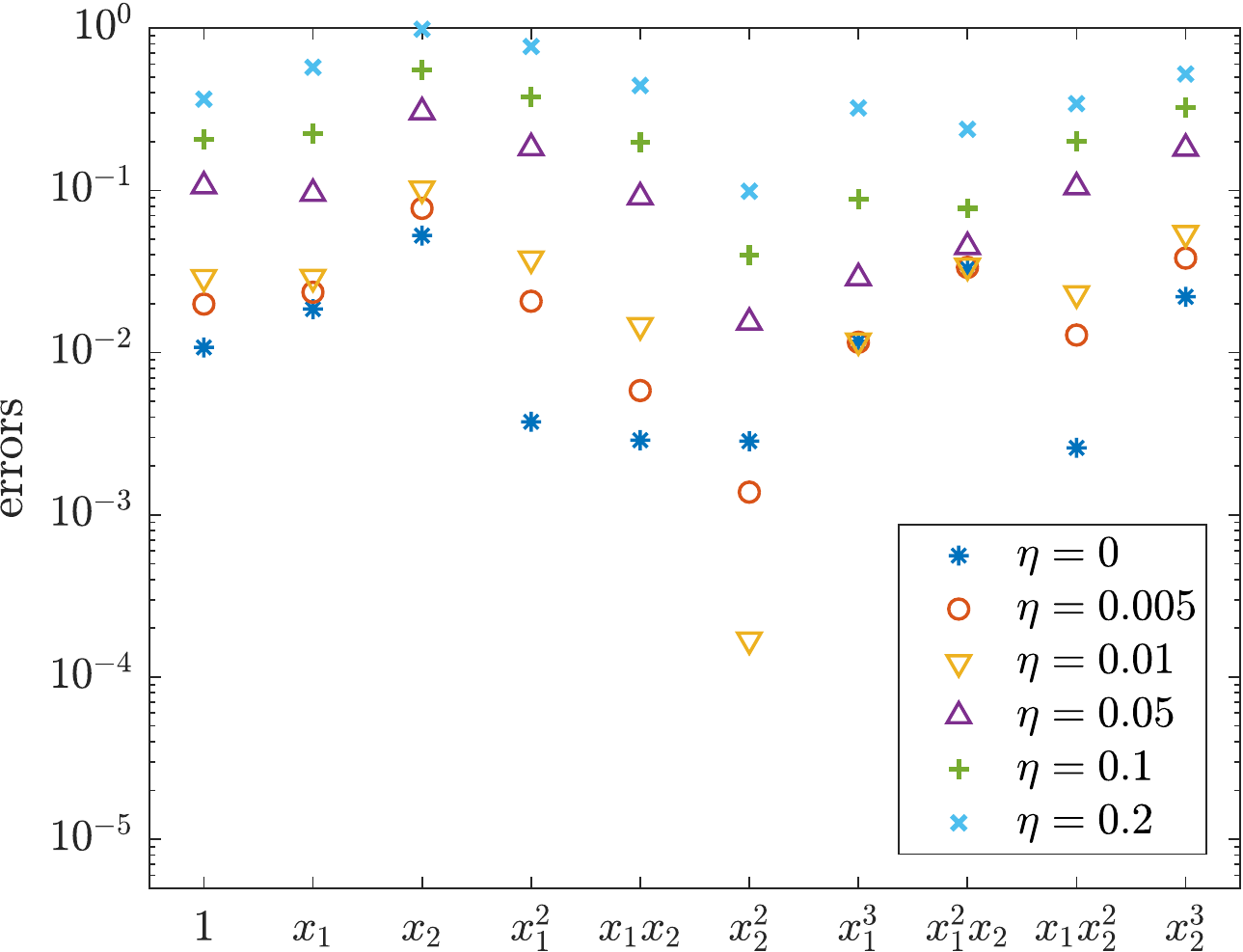}}
		\caption{\small
			Example \ref{ex:nonODE3}: Errors in each terms
                        of the recovered equation at different noise levels.   
		}\label{fig:nonODE3err}
	\end{figure}
	
	 \begin{figure}[htbp]
	\centering
{\includegraphics[width=0.47\textwidth]{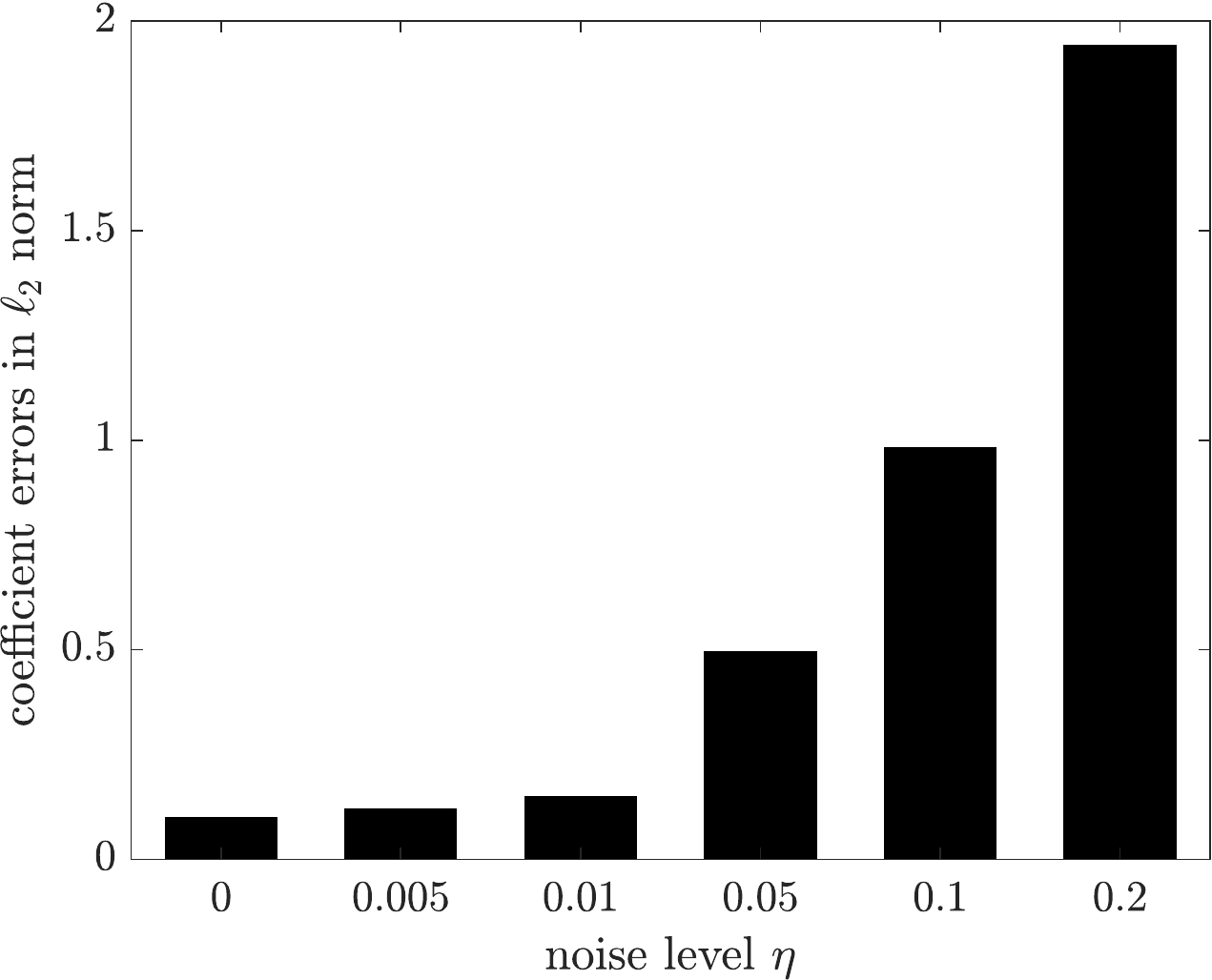}}
	\caption{\small
		Example \ref{ex:nonODE3}: The $\ell_2$ norm of the
                errors in the recovered coefficients of the equations.
	}\label{fig:nonODE3err2}
\end{figure}

\end{example}

\begin{example}\label{ex:nonODE4} \rm
We now consider a damped pendulum problem. 
Let $l$ be the length of the
pendulum, $\alpha$ be the damping constant,  and $g$ be the
gravitational constant.
 Then, for the unit mass, the governing equation is well known as
$$
\frac{d^2 \theta}{d t^2} + \frac{\alpha}{l} \frac{ d\theta }{dt} + \frac{g}{l} \sin\theta = 0,
$$
where $\theta$ is the angle from its equilibrium vertical position.
By defining $u_1=\theta$ and $u_2=\dot \theta$, we rewrite the
equation as
\begin{equation}\label{eq:nonODE4}
\begin{cases}
\dot u_1=u_2  ,\\
\dot u_2=- \frac{g}{l} \sin u_1  - \frac{\alpha}{l} u_2 .
\end{cases}
\end{equation} 
We set $l=1.1$ and $\alpha=0.22$, and 
use a large number of short trajectory data ($M=200,J_m=2$) to recover
the equation.
The initial state of each trajectory is randomly and uniformly chosen from the region $D=[-\pi,\pi]\times[-2\pi,2\pi]$. 
 The least squares $\ell_2$-regression is used, with polynomial
 approximation at 
different degrees $n$. 
Figure \ref{fig:nonODE4} shows the phase portraits for the system \eqref{eq:nonODE4} and the learned systems. 
As $n$ increases, the recovered dynamics becomes more accurate and can
capture the true structures of the solution.
More detailed examination of the learned system is
shown in Figure \ref{fig:nonODE4val}, with polynomial degree $n=6$.
We observe good agreements 
between the solution of the learned system ${\bf x}(t)$ and true
solution ${\bf u}(t)$, at an arbitrarily chosen initial state
$(-1.193,-3.876)$. The error evolution over time remains bounded at
$10^{-2}$ level.

 \begin{figure}[htbp]
	\centering
	\subfigure[$n=2$]{\includegraphics[width=0.47\textwidth]{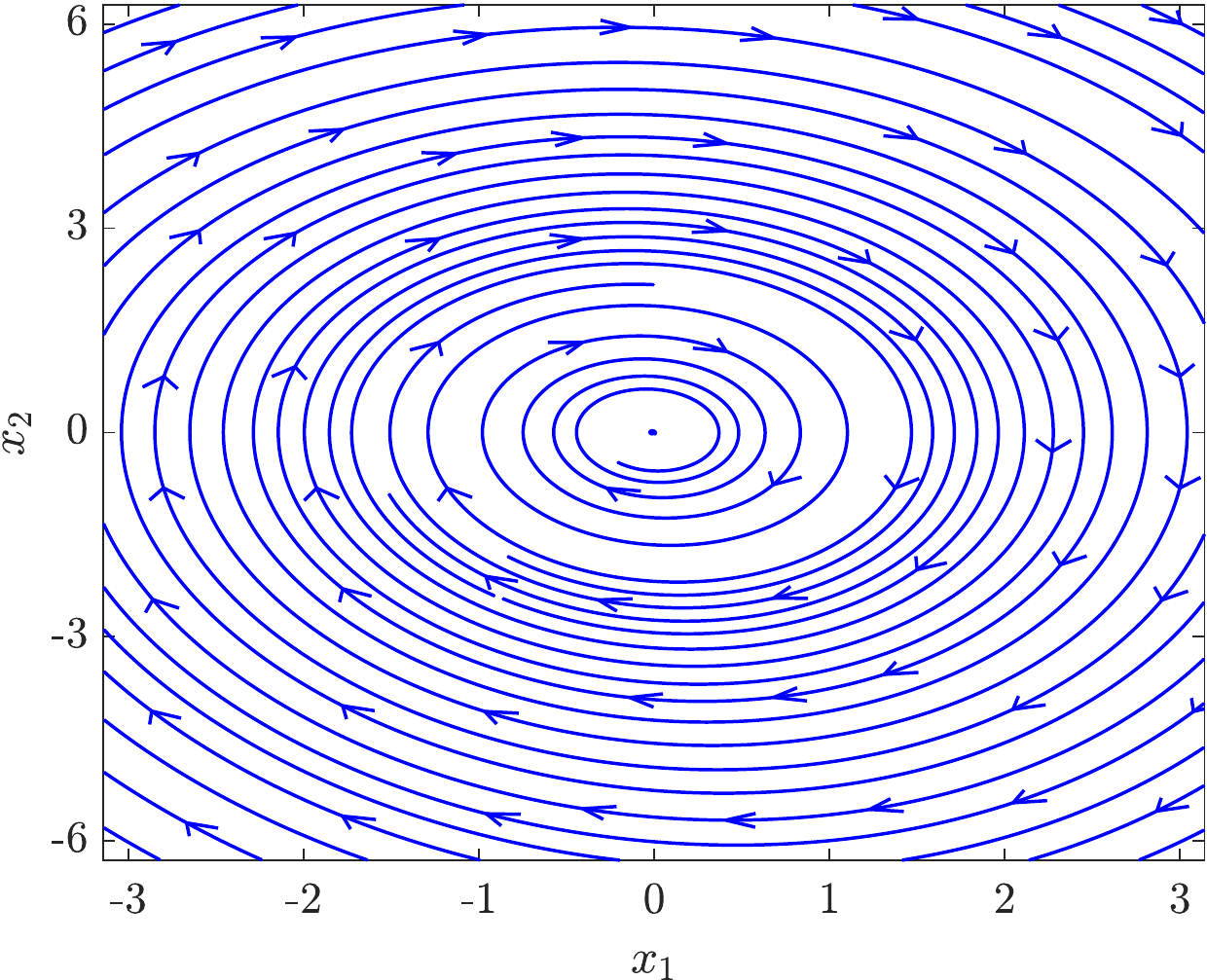}}~~~~
	\subfigure[$n=4$]{\includegraphics[width=0.47\textwidth]{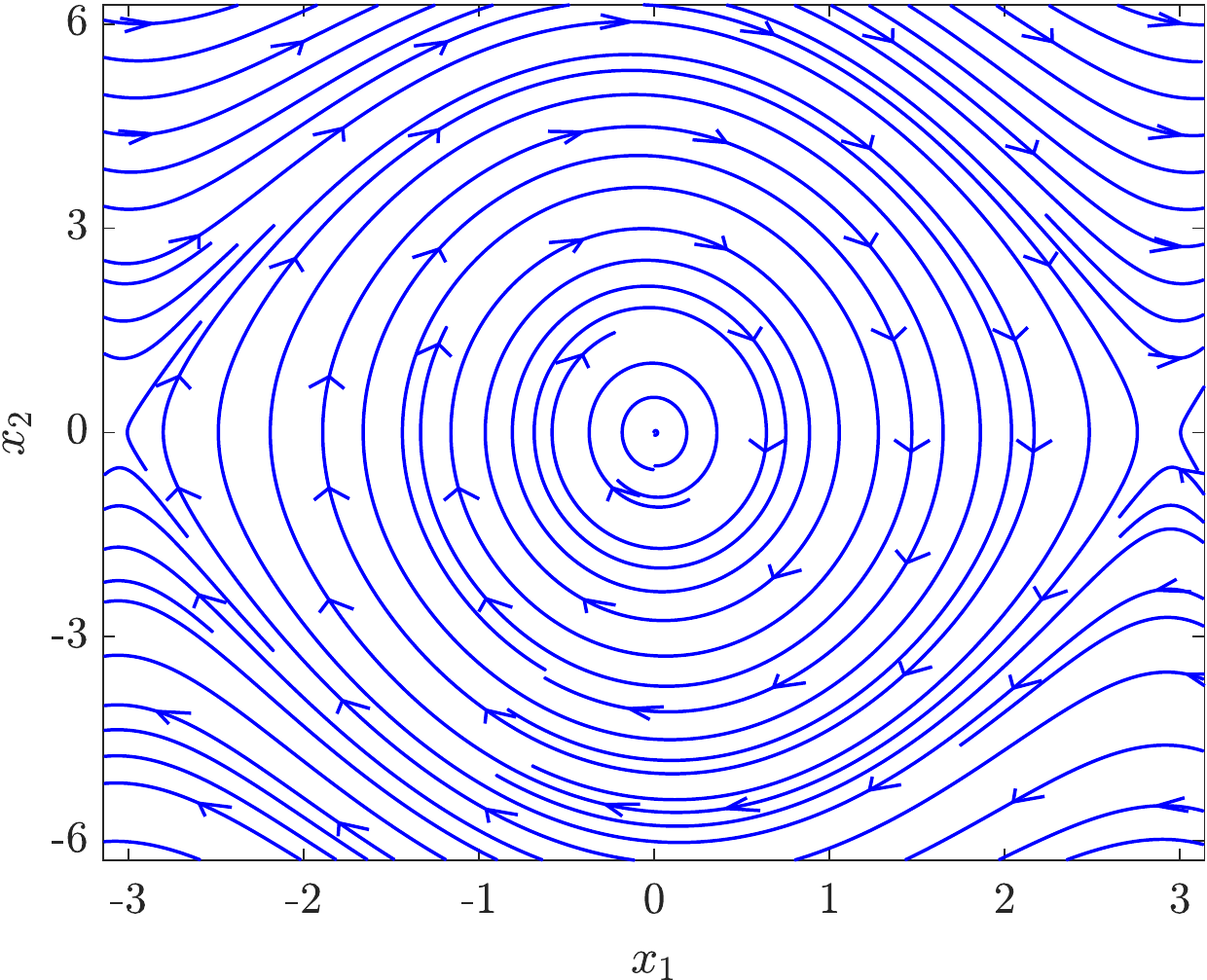}}
	\subfigure[$n=6$]{\includegraphics[width=0.47\textwidth]{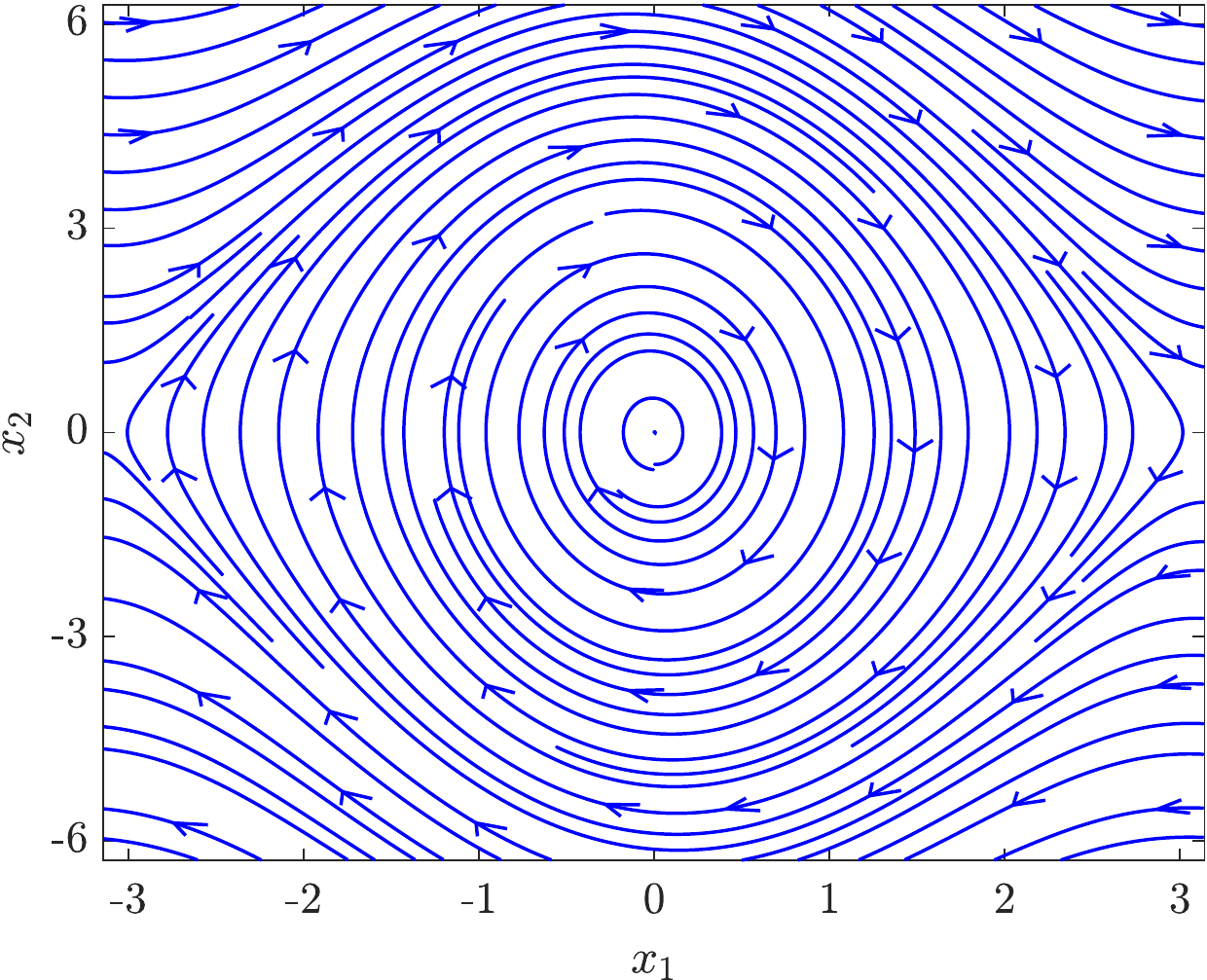}}~~~~
	\subfigure[Exact]{\includegraphics[width=0.47\textwidth]{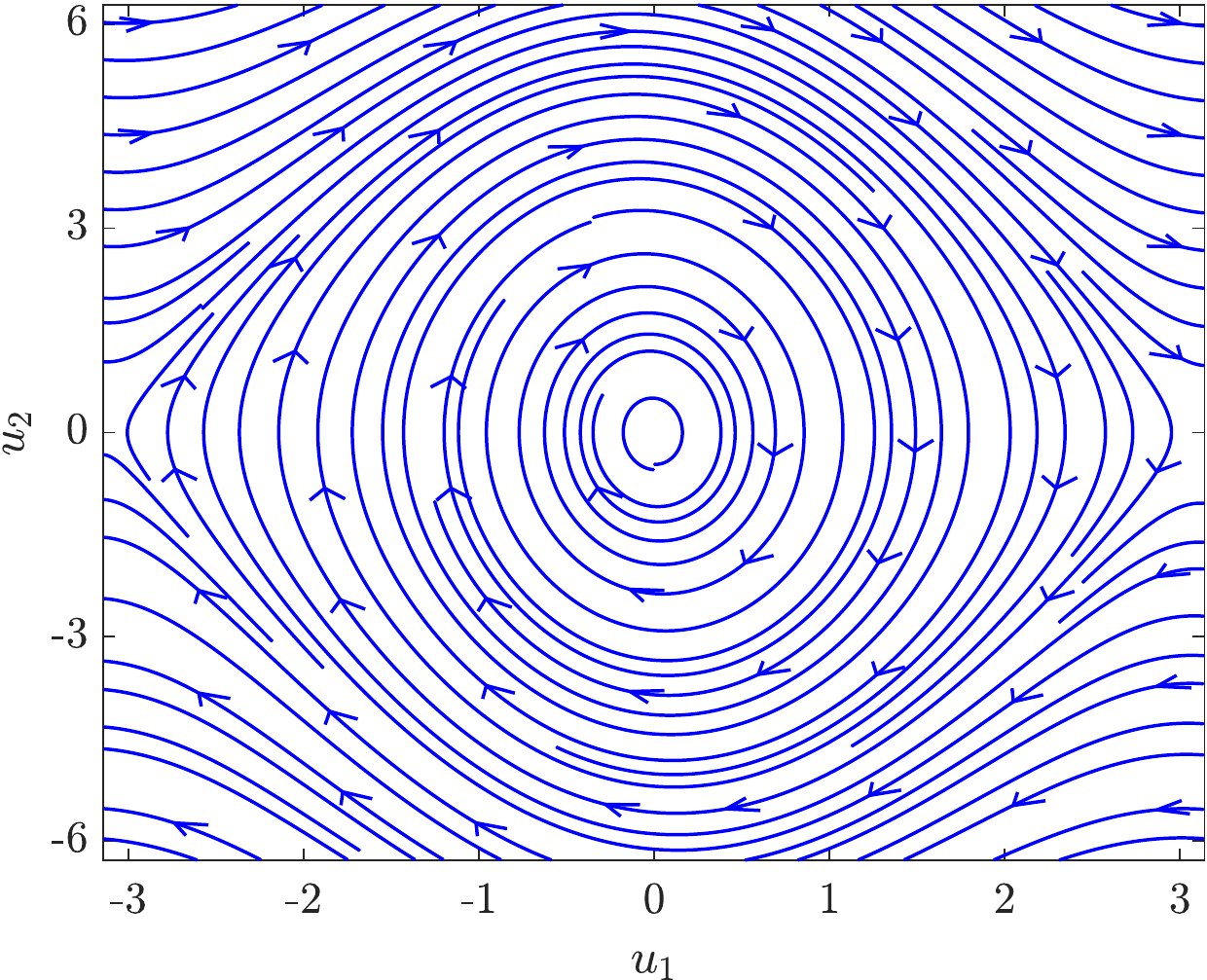}}
	\caption{\small
		Example \ref{ex:nonODE4}: Phase portraits for the system \eqref{eq:nonODE4} and its approximate systems learned from measurement data with different $n$.   
	}\label{fig:nonODE4}
\end{figure}

\begin{figure}[htbp]
	\centering
	\includegraphics[width=0.47\textwidth]{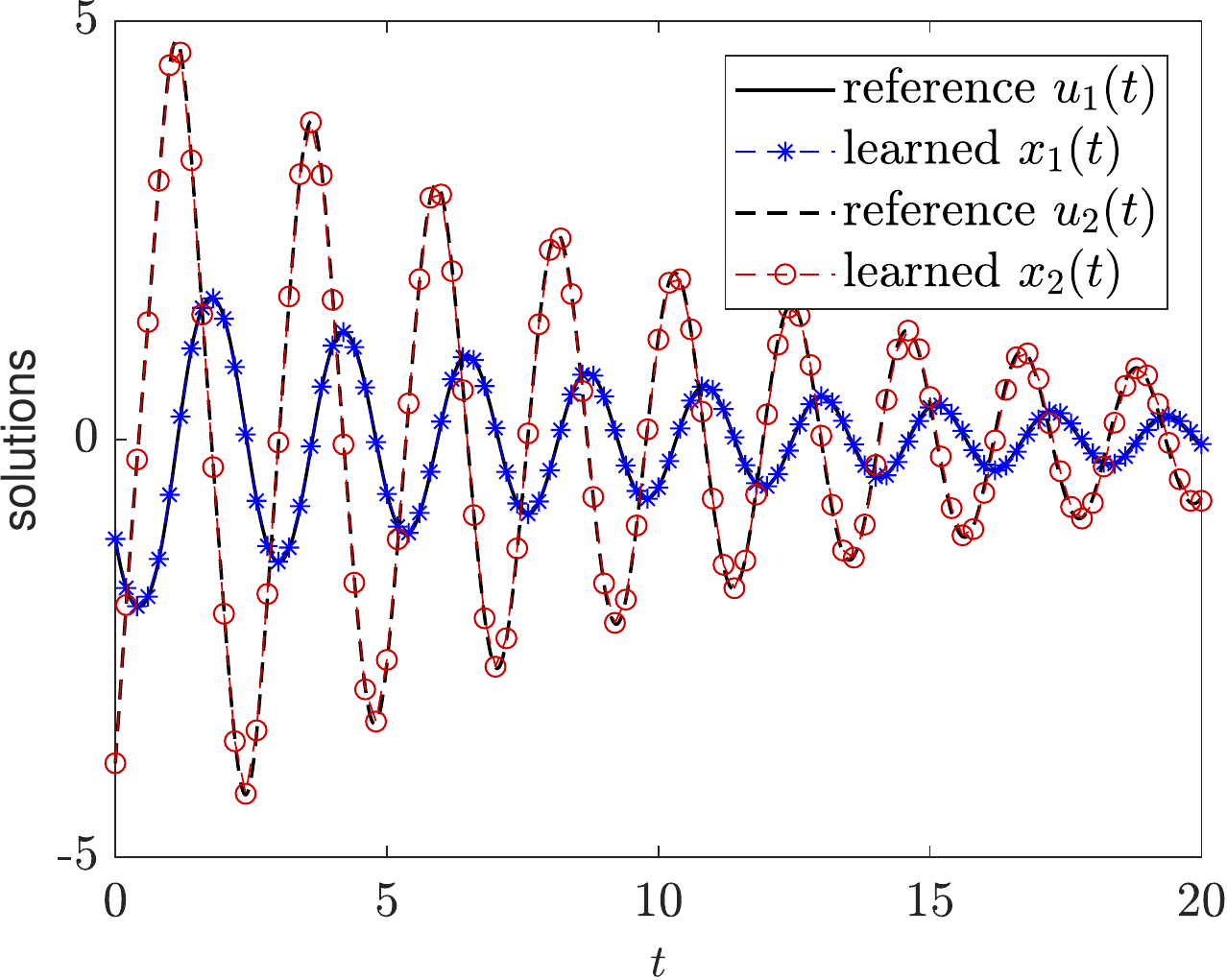}~~~~
	\includegraphics[width=0.47\textwidth]{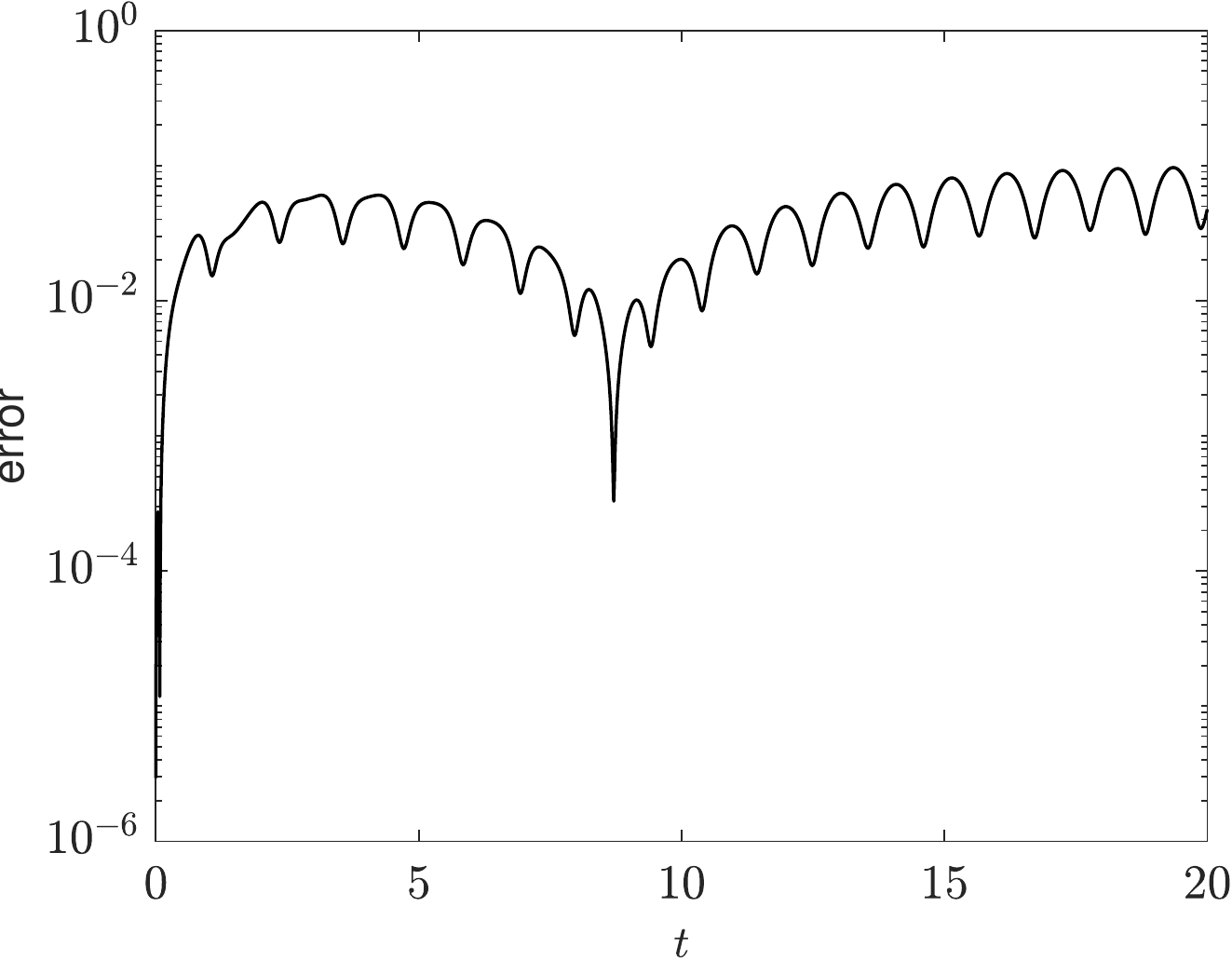}
	\caption{\small
		Example \ref{ex:nonODE4}: Validation of the recovered system ($n=6$)  
		with initial state $(-1.193,-3.876)$. Left: solutions. Right: error $\| {\bf x}(t) -  {\bf u} (t) \|_2$. The horizontal axis denotes the time.   
	}\label{fig:nonODE4val}
\end{figure}
\end{example}

\subsection{Nonlinear DAEs}

We now study 
semi-explicit DAEs \eqref{eq:DAE} 
from measurement data. The learned approximating DAEs are 
denoted by 
\begin{subnumcases}{\label{eq:app:DAE}}
\dfrac{d {\bf x}(t)}{dt}={ \tilde{\bf f}}({\bf x}(t)),\label{eq:app:DAEa} \\
{\bf y}(t)= \tilde{\bf g} ( {\bf x}(t) ), \label{eq:app:DAEb}
\end{subnumcases} 
where ${\bf x}(t)$ and ${\bf y}(t)$ are 
approximations to ${\bf u}(t)$ and ${\bf v}(t)$, respectively. 
In the tests, 
the measurement data are first generated by numerically solving the true DAEs using the forth-order explicit  Runge-Kutta method 
with a small time step-size. Possible noises 
or corruption errors are then added to the trajectories to generate
our synthetic data.

\begin{example}\label{ex:DAE1} \rm
The first DAEs model we consider is a model for
nonlinear electric network \cite{pulch2013polynomial}, 
defined as
\begin{equation}\label{eq:DAE1}
\begin{cases}
\dot u_1 = v_2/C,
\\
\dot u_2  = u_1/L,
\\
0 = v_1 - (G_0-G_\infty) U_0 {\rm tanh} (u_1/U_0) - G_\infty u_1,
\\
0=v_2 + u_2 + v_1,
\end{cases}
\end{equation}
where $u_1$ denotes the node voltage, and 
$u_2$, $v_1$ and $v_2$ are branch currents. Following \cite{pulch2013polynomial}, 
the  physical parameters are specified as  
$C=10^{-9}$, $L=10^{-6}$, $U_0 = 1$, $G_0=-0.1$ and $G_\infty=0.25$. 
In our test, $D=[-2,2]\times [-0.2,0.2]$, 
$M=200$ sequences of trajectory data with $\Delta t=10^{-9}$ 
and $J_m=19$ 
are generated by numerically solving \eqref{eq:DAE1} and then
adding random noises. 
The added noises in $u_1$ and $u_2$
follow  uniform distribution on $[-0.0001,0.0001]$, 
while the noises in $v_1$ and $v_2$ follow
normal distribution $N(0,0.001^2)$ and $N(0,0.005^2)$ respectively. We then use the $\ell_2$-regression with $n=8$  
to learn the governing equations from these data. 
The phase portraits for the system $\dot{\bf u} = {\bf f}({\bf u})$ and the learned system $\dot{\bf x}=\tilde{\bf f}({\bf x})$ 
are compared in Figure \ref{fig:DAE1}. 
It is observed that the proposed method accurately reproduces 
the dynamics and detects the correct features. 
Comparison the solutions between the recovered DAEs and the true system is shown
Figure \ref{fig:DAE1val}, using an arbitrary choice of
the initial state $(0,0.1)$. The evolution of the numerical errors in the solution of
the learned system are shown
in Figure \ref{fig:DAE1err}.

 \begin{figure}[htbp]
	\centering
	\subfigure{\includegraphics[width=0.47\textwidth]{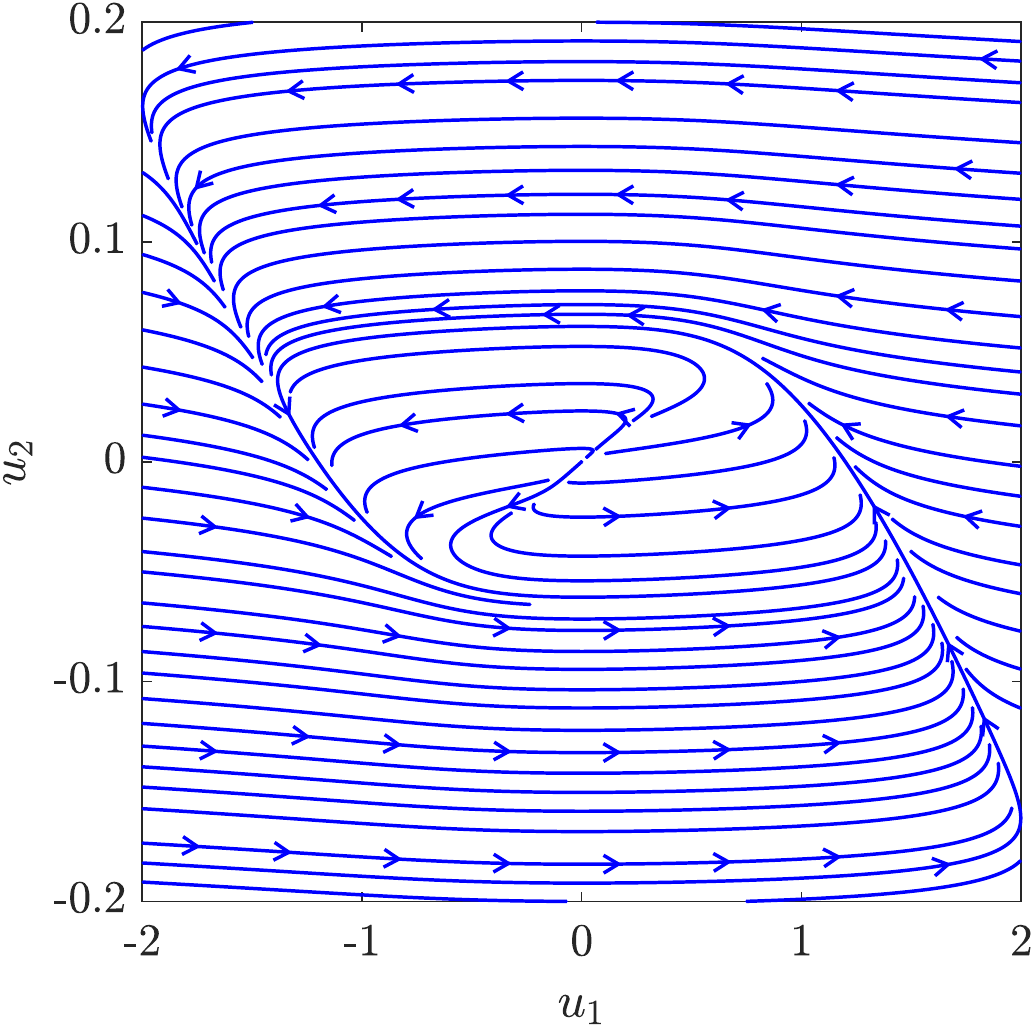}}~~~~
	\subfigure{\includegraphics[width=0.47\textwidth]{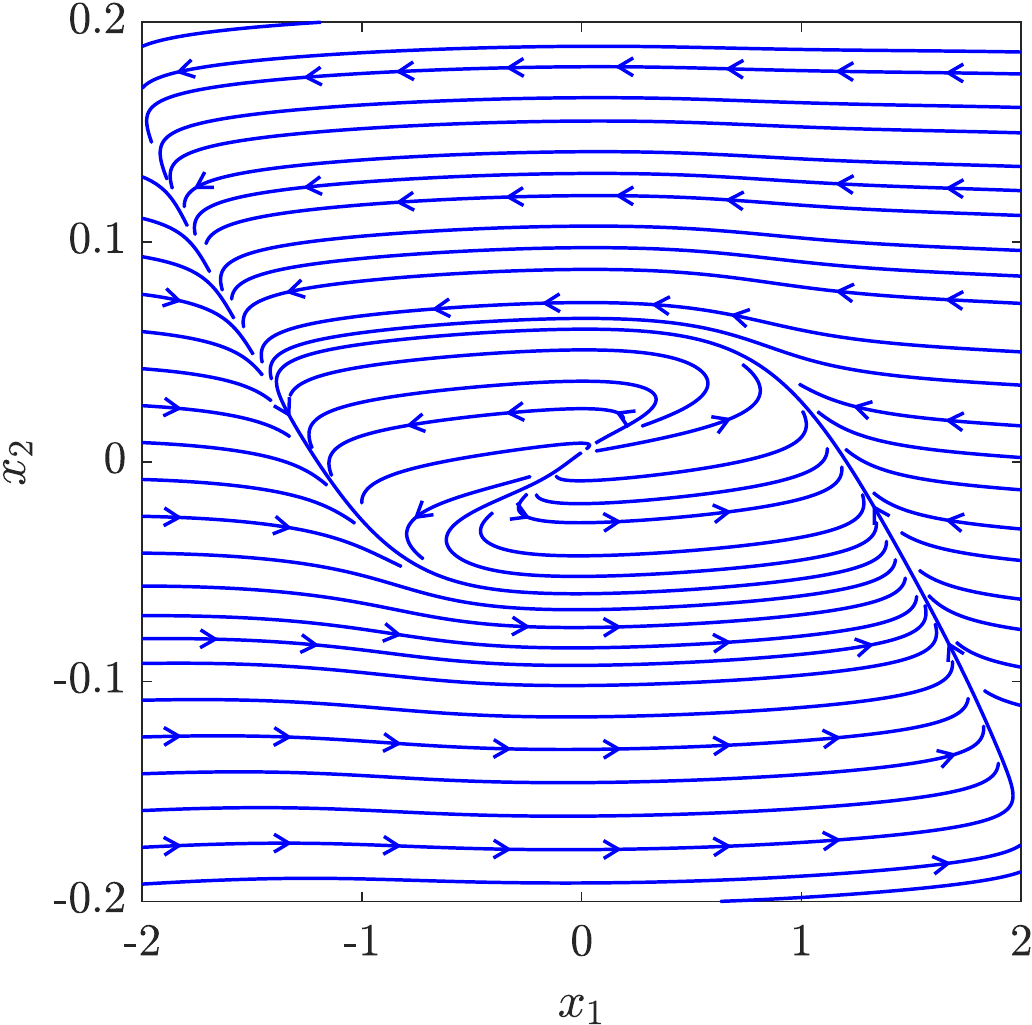}}
	\caption{\small
		Example \ref{ex:DAE1}: Phase portraits for  the system $\dot{\bf u} = {\bf f}({\bf u})$ (left) and the learned system $\dot{\bf x}=\tilde{\bf f}({\bf x})$ (right).   
	}\label{fig:DAE1}
\end{figure}

\begin{figure}[htbp]
	\centering
	\includegraphics[width=0.47\textwidth]{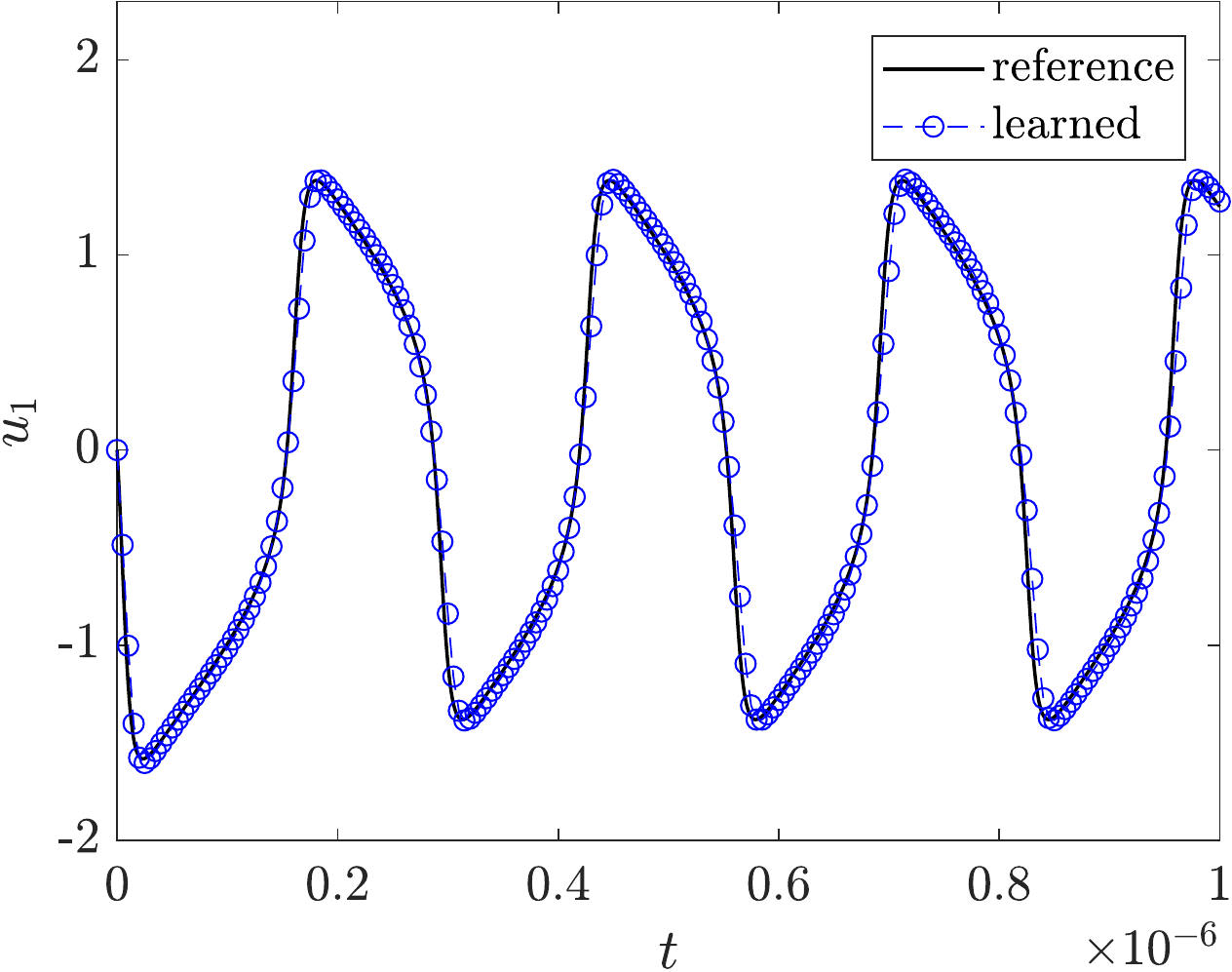}~~~~
	\includegraphics[width=0.47\textwidth]{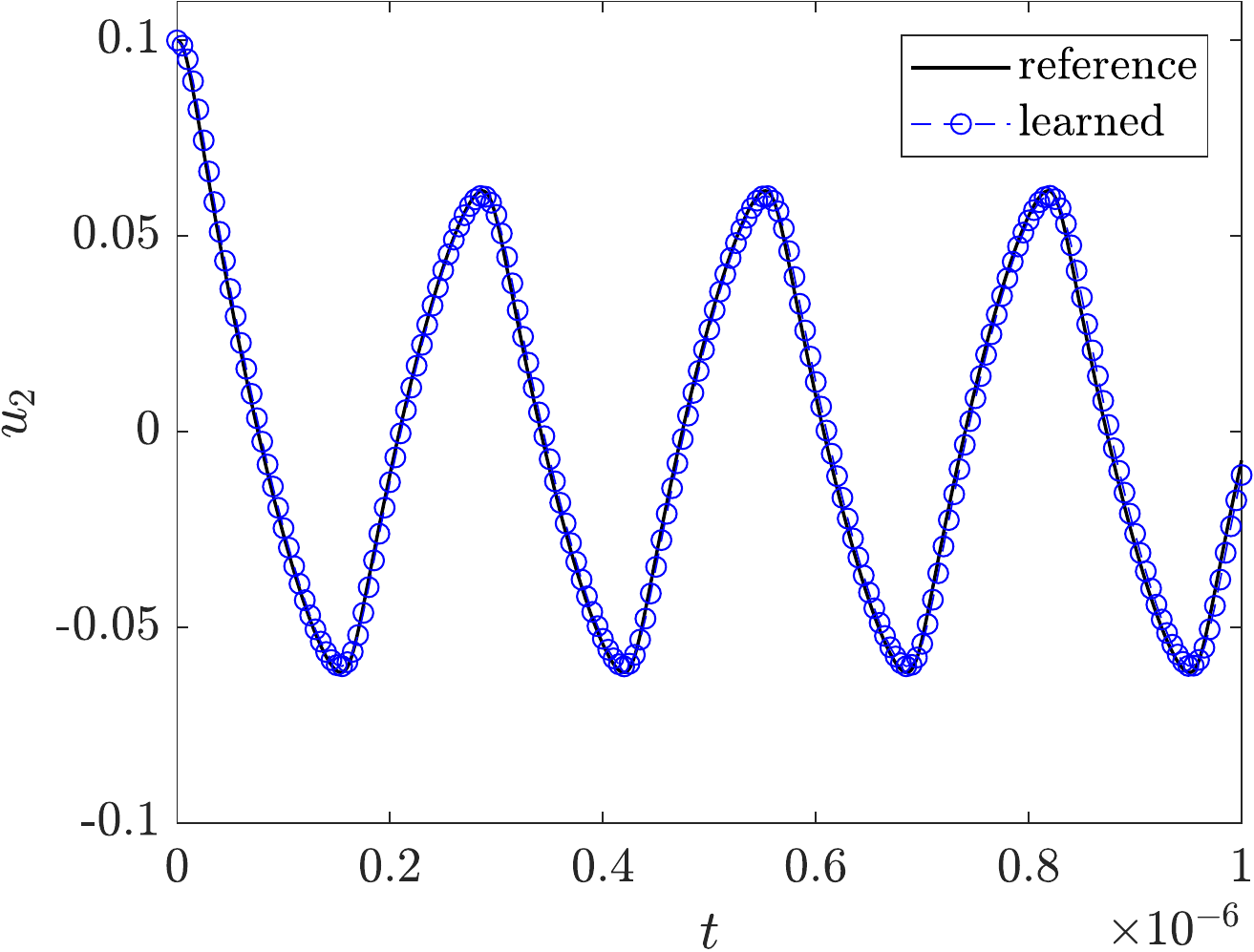}
	\includegraphics[width=0.47\textwidth]{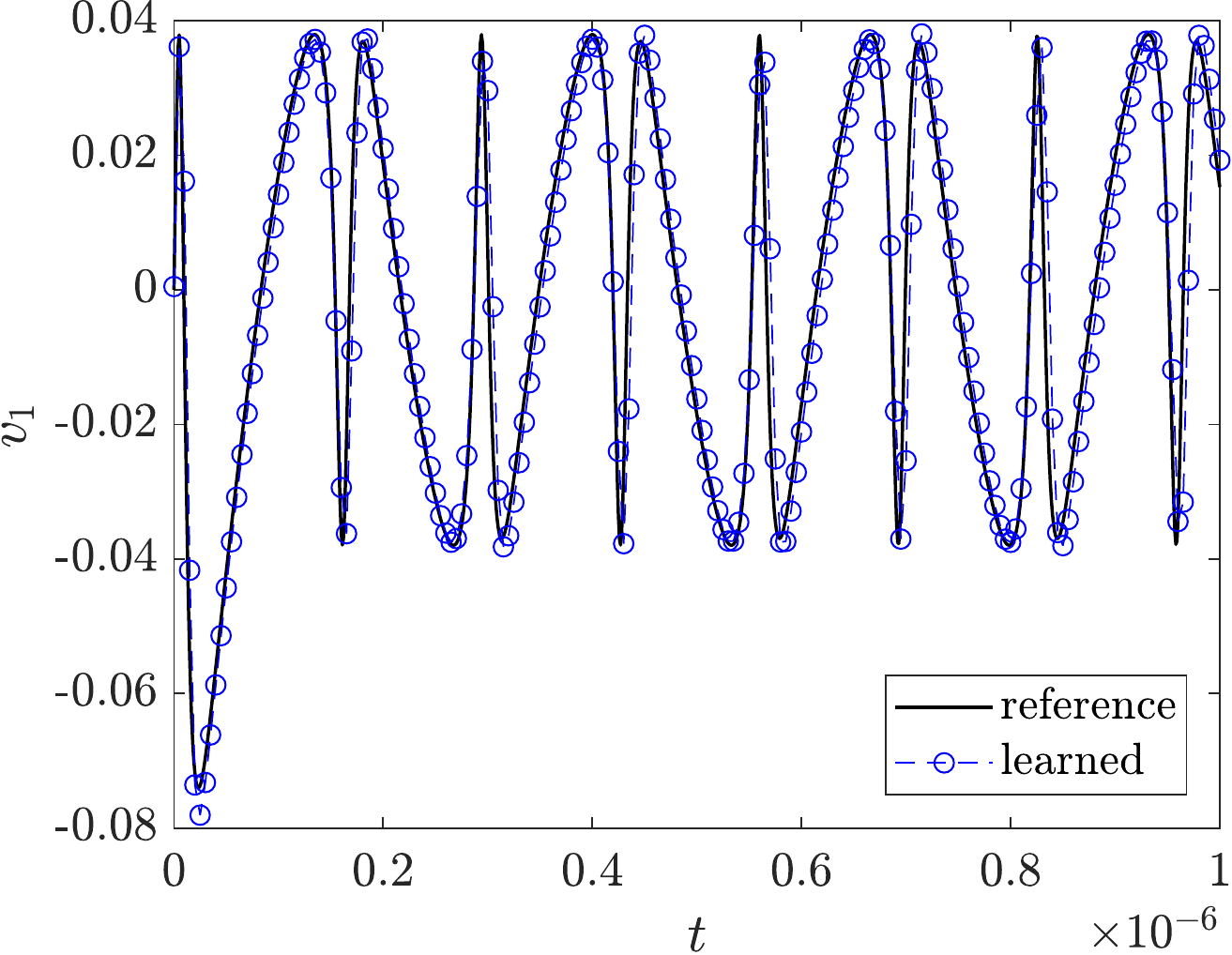}~~~~
	\includegraphics[width=0.47\textwidth]{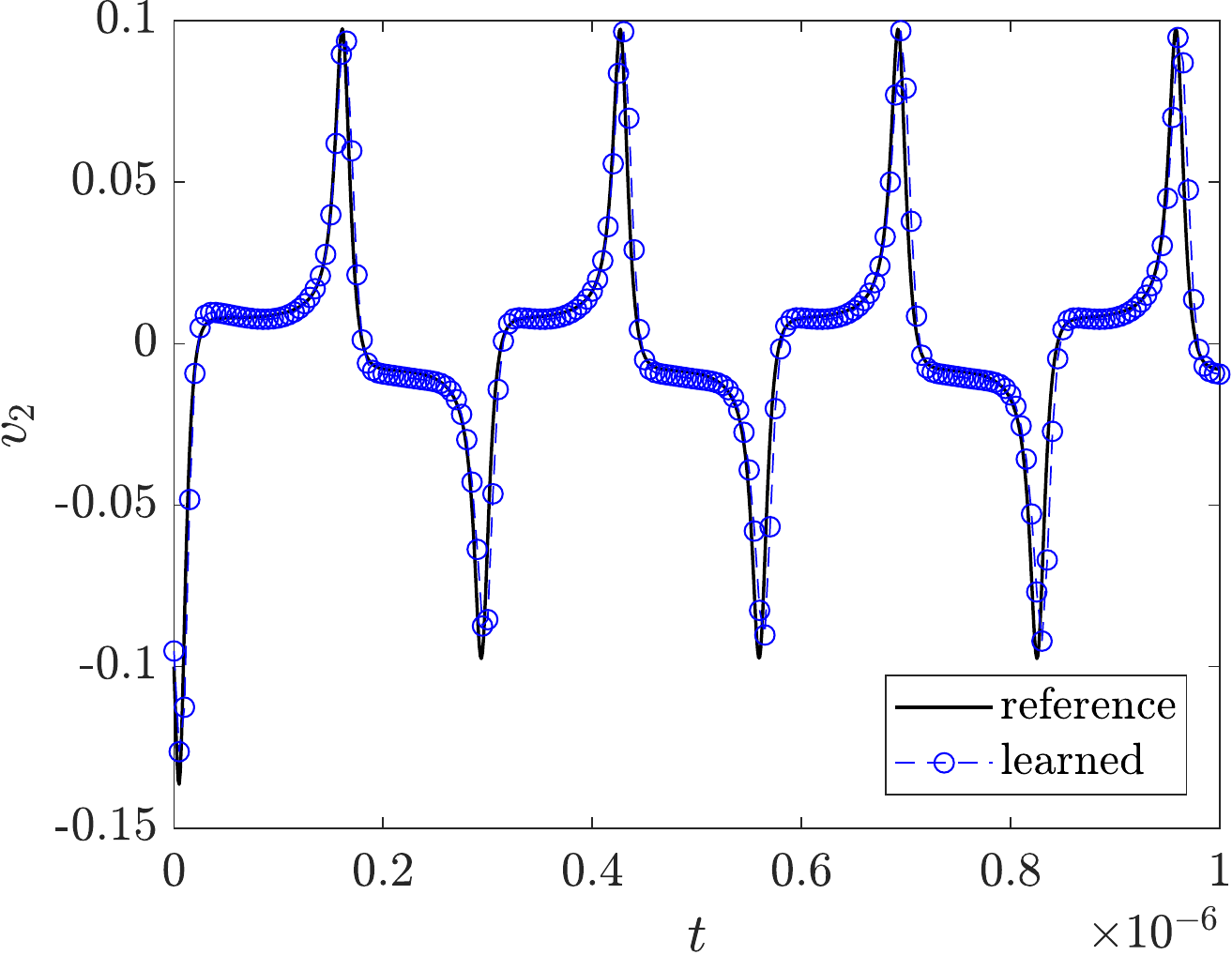}	
	\caption{\small
		Example \ref{ex:DAE1}: Validation of the recovered system  
		with initial state $(0,0.1)$. 
		The solid lines denote the solution of the exact system, 
		and the symbols ``$\circ$'' represent 
		the solution of the learned system. 
		The horizontal axis denotes the time.   
	}\label{fig:DAE1val}
\end{figure}

\begin{figure}[htbp]
	\centering
	\includegraphics[width=0.47\textwidth]{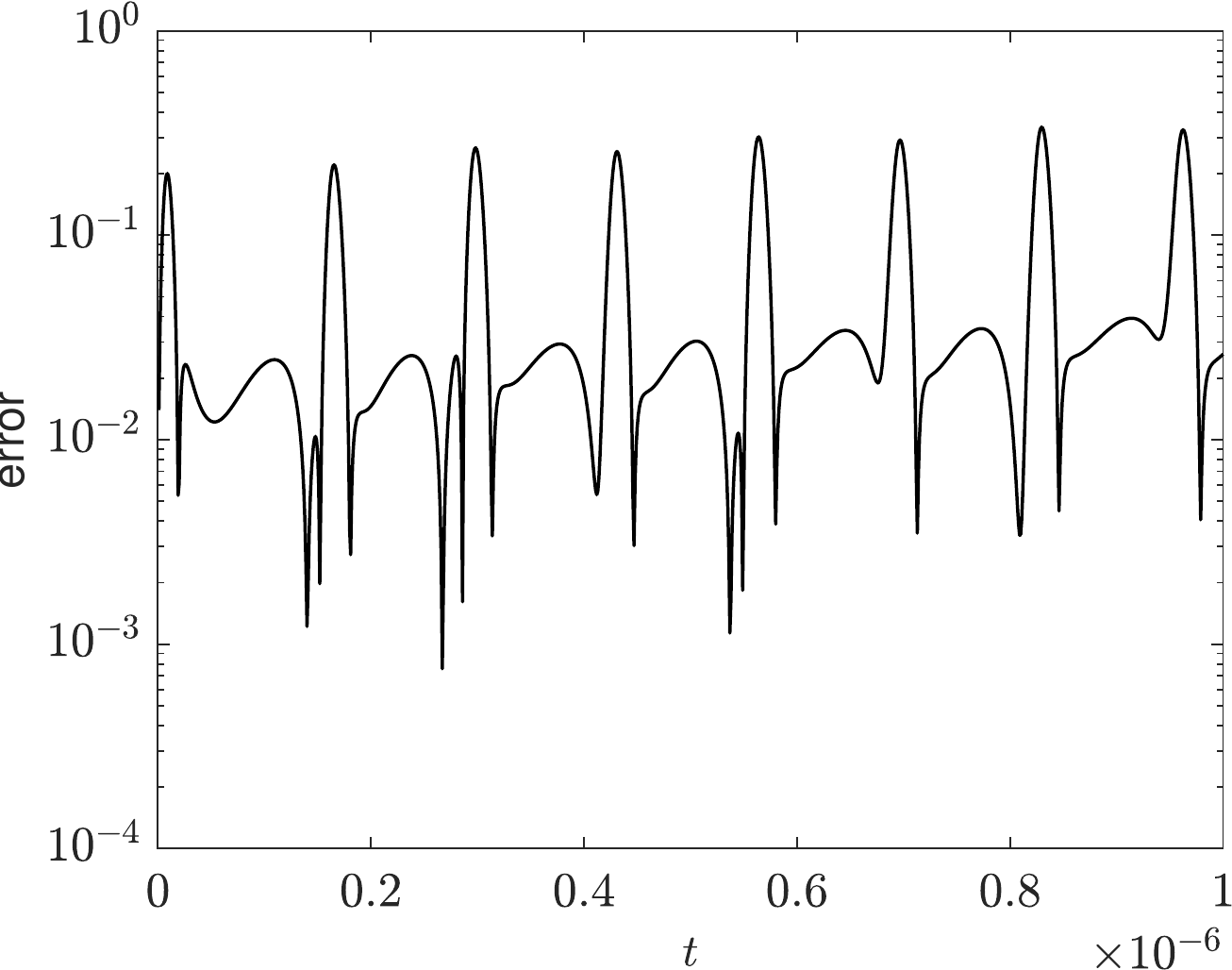}~~~~
	\includegraphics[width=0.47\textwidth]{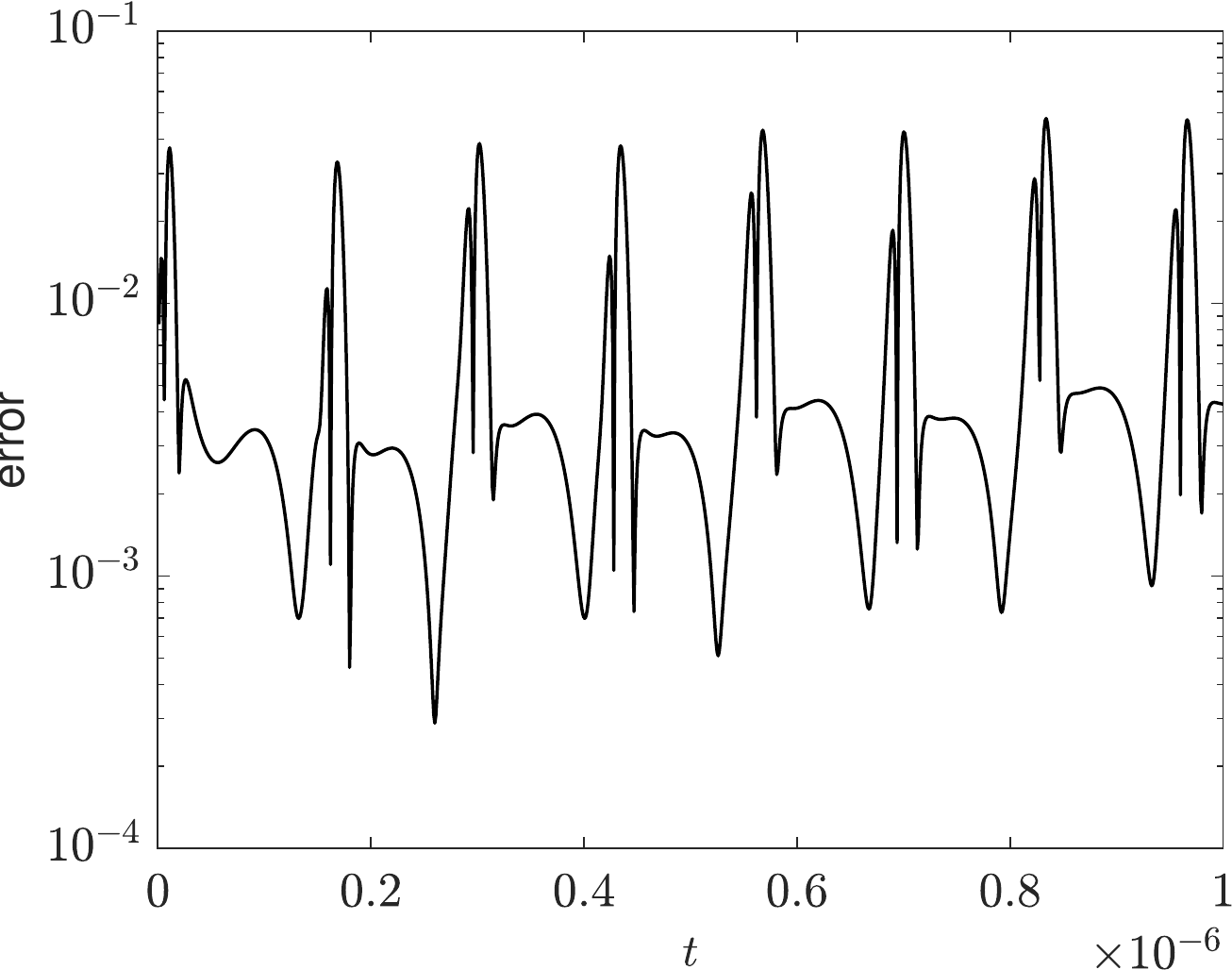}
		\caption{\small
		Example \ref{ex:DAE1}: Same as Figure \ref{fig:DAE1val} except for  
		the errors $\| {\bf x}(t) -  {\bf u} (t) \|_2$ (left) 
		and $\| {\bf y}(t) -  {\bf v} (t) \|_2$ (right). 
	}\label{fig:DAE1err}
\end{figure}
\end{example}

\begin{example}\label{ex:DAE2} \rm
We now
consider a model for a genetic toggle switch in {\em Escherichia
  coli}, which was proposed in \cite{gardner2000construction} 
and numerically studied in \cite{xiu2007efficient}. 
It is composed of two repressors and two constitutive promoters, where
each promoter is inhibited by the repressor that is transcribed by the opposing promoter.
Details of experimental measurement can be found in [10], where the following 
model was also constructed,  
\begin{equation}\label{eq:DAE2}
\begin{cases}
\dot u_1 = \dfr{ \alpha_1}{1+u_2^\beta}-u_1,
\\
\dot u_2  = \dfr{\alpha_2}{1+v_1^\gamma} - u_2,
\\
0 = v_1 + \varepsilon \sin v_1 - \frac{ u_1 }{ ( 1 + {\rm [IPTG]}/K )^\eta}. 
\end{cases}
\end{equation}
In the last equation, we add a small perturbation term $\varepsilon \sin v_1$ to increase the nonlinearity. 
Here, $u_1$ and $u_2$ are the concentration of the two repressors, respectively; $\alpha_1$ and $\alpha_2$ 
are the effective rates of synthesis of the repressors; 
$\gamma$ and $\beta$ represent cooperativity of 
repression of the two promoters, respectively; and $\rm [IPTG]$ 
is the concentration of IPTG, the chemical compound that 
induces the switch.  
The values of these parameters are taken as 
$\alpha_1=156.25$, $\alpha_2=15.6$, $\beta=2.5$, $\gamma=1$, $\eta = 2.0015$. 
In the test, $D=[0,20]^2$ and we generate the measurement data 
by numerically solving \eqref{eq:DAE2} with  
$\rm [IPTG]=10^{-5}$ and $\varepsilon=0.01$. 
Suppose we have $M=500$ bursts of trajectory data with $J_m=2$, 
and five ($5$) arbitrarily chosen bursts of data contain corruption errors following  
normal distribution $N(1,1)$. 
We employ
the $\ell_1$-regression approach with $n=16$ order of
nominalized Legendre polynomial basis.  
The phase portraits for  the system $\dot{\bf u} = {\bf f}({\bf u})$ and the learned system $\dot{\bf x}=\tilde{\bf f}({\bf x})$ 
are plotted in Figure \ref{fig:DAE2}. As we can see, the learned structures agree well with 
the ideal structures everywhere in $D$. 
This demonstrates the ability of our method in finding 
the underlying governing equations. 
The accuracy of the learned model is further examined in  
Figures \ref{fig:DAE2val} and \ref{fig:DAE2err}.

 \begin{figure}[htbp]
	\centering
	\subfigure{\includegraphics[width=0.47\textwidth]{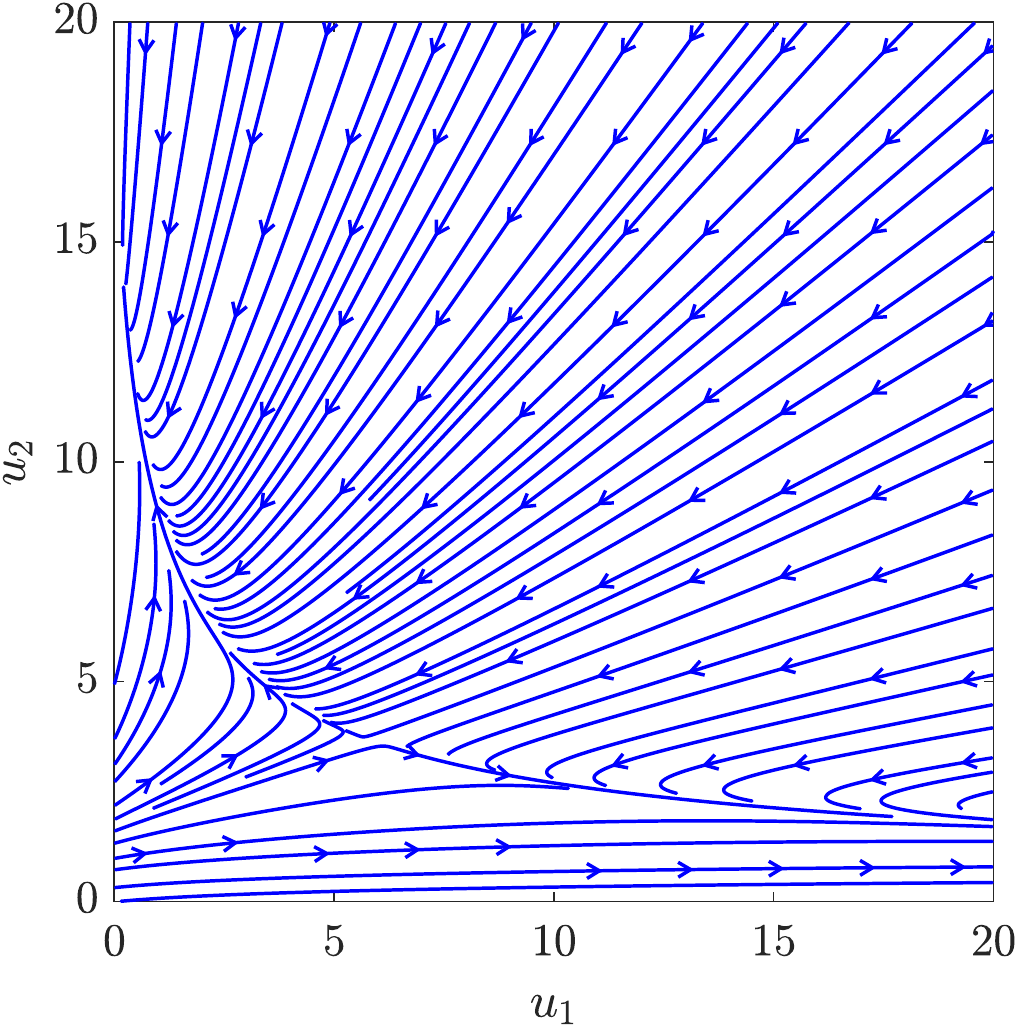}}~~~~
	\subfigure{\includegraphics[width=0.47\textwidth]{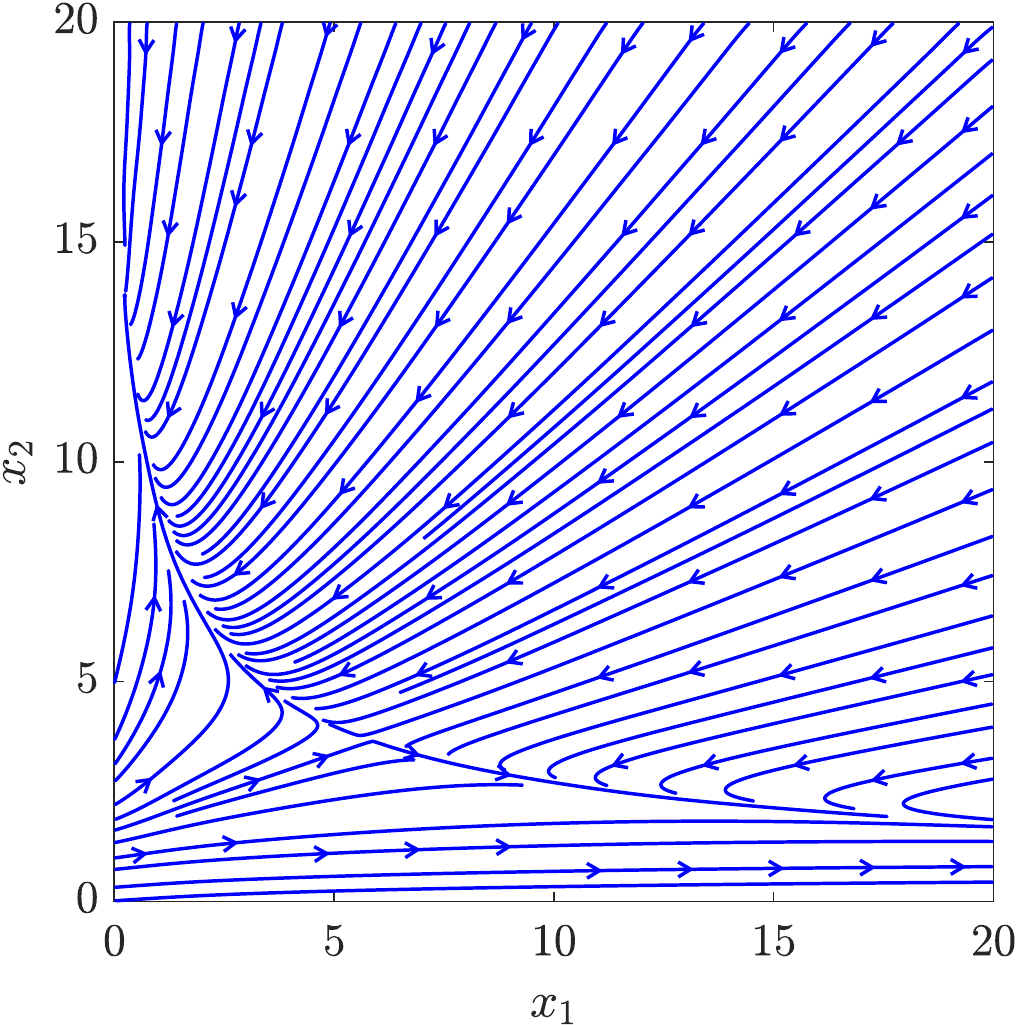}}
	\caption{\small
		Example \ref{ex:DAE2}: Phase portraits for  the system $\dot{\bf u} = {\bf f}({\bf u})$ (left) and the learned system $\dot{\bf x}=\tilde{\bf f}({\bf x})$ (right). 
	}\label{fig:DAE2}
\end{figure}

\begin{figure}[htbp]
	\centering
	\includegraphics[width=0.47\textwidth]{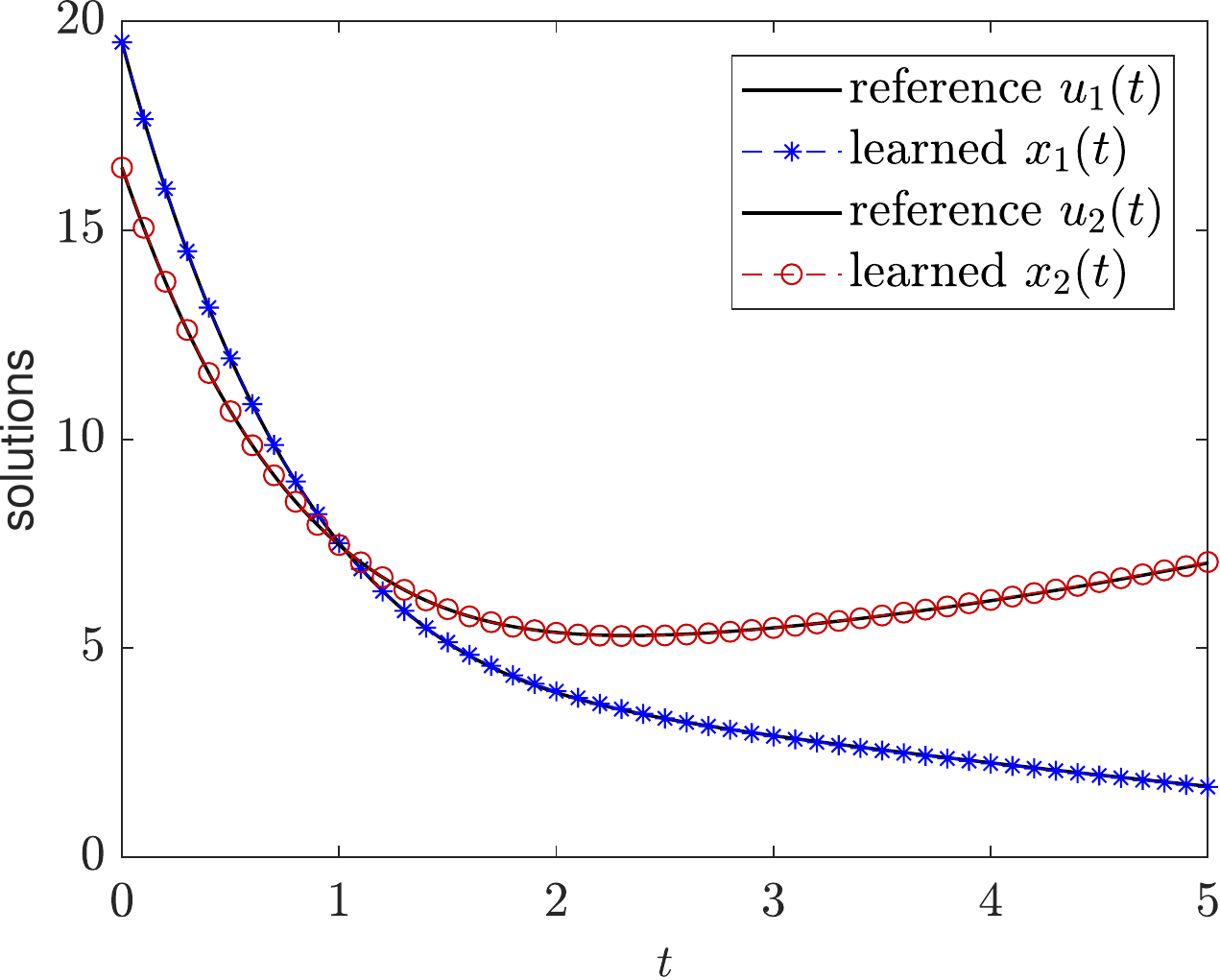}~~~~
	\includegraphics[width=0.47\textwidth]{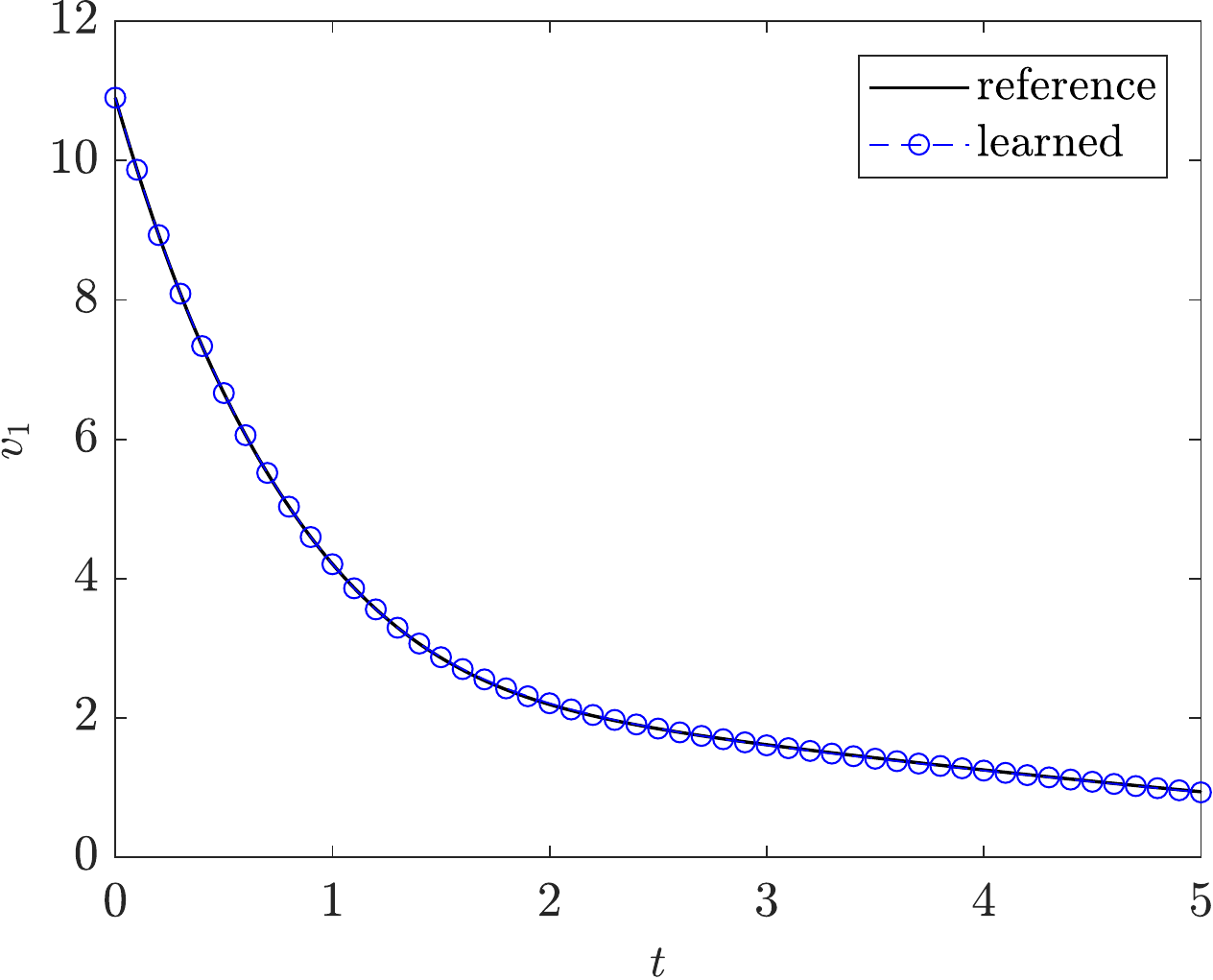}	
	\caption{\small
		Example \ref{ex:DAE2}: Validation of the recovered system  
		with initial state $(19.5,16.5)$. The horizontal axis denotes the time.   
	}\label{fig:DAE2val}
\end{figure}

\begin{figure}[htbp]
	\centering
	\includegraphics[width=0.47\textwidth]{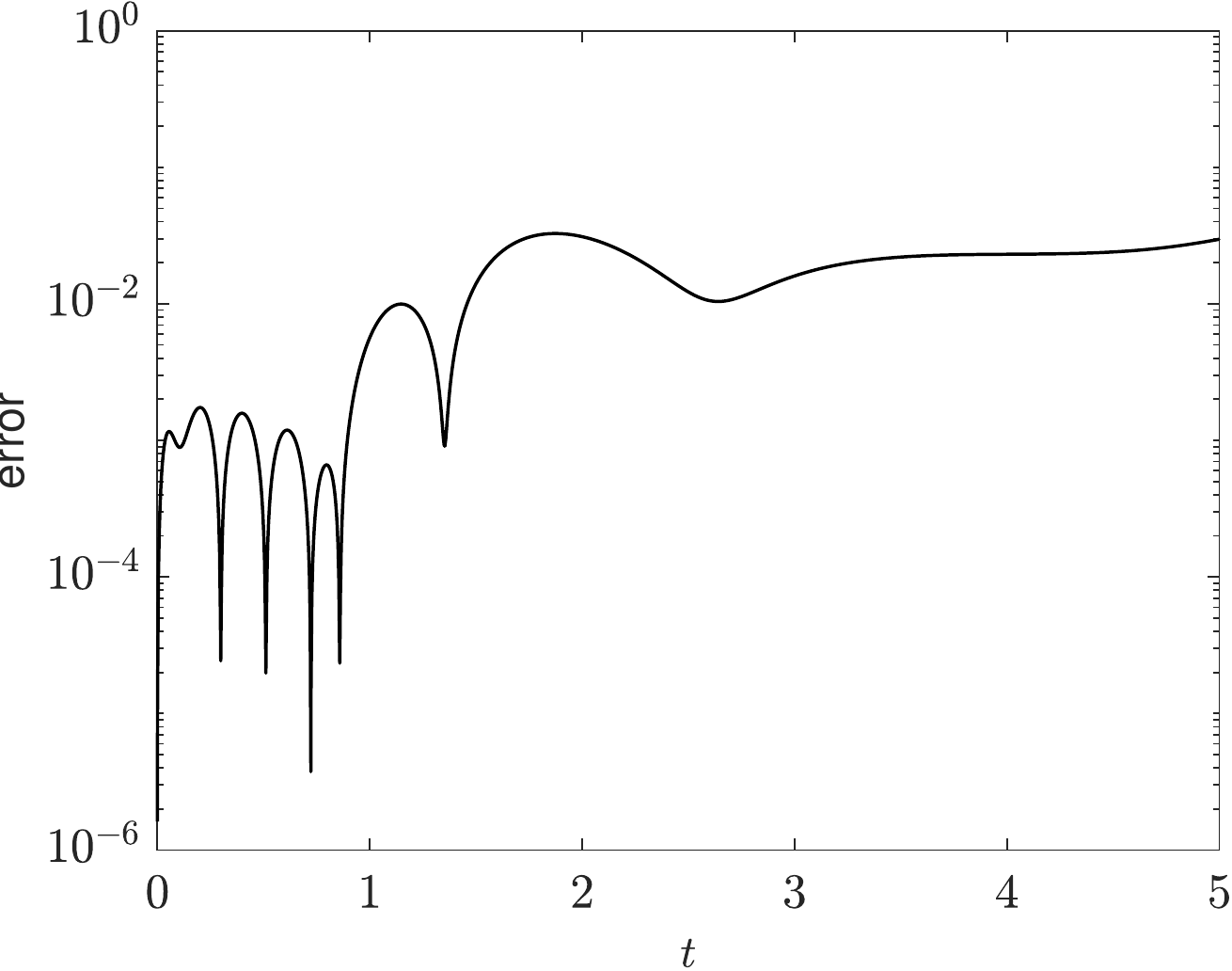}~~~~
	\includegraphics[width=0.47\textwidth]{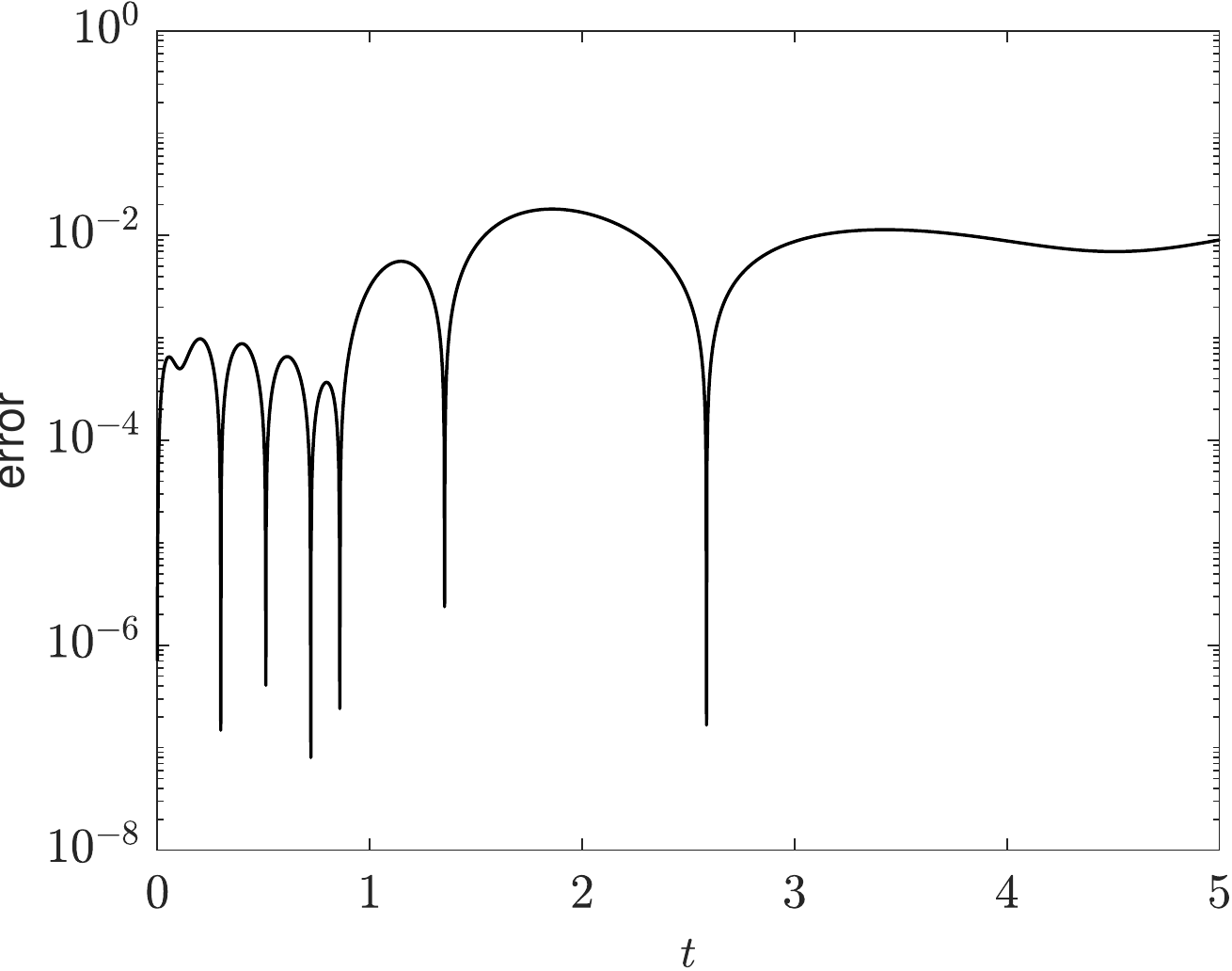}
	\caption{\small
		Example \ref{ex:DAE2}: Same as Figure \ref{fig:DAE2val} except for  
		the errors $\| {\bf x}(t) -  {\bf u} (t) \|_2$ (left) 
		and $\| {\bf y}(t) -  {\bf v} (t) \|_2$ (right). 
	}\label{fig:DAE2err}
\end{figure}
\end{example}

\begin{example}\label{ex:DAE3} \rm
Our last example is a kinetic model of an isothermal
batch reactor system, 
which was proposed in \cite{biegler1986nonlinear} 
and studied in \cite{becerra2001applying,kadu2010modified}, among
other literature. 
It was mentioned in \cite{becerra2001applying} 
that the example in its original form was given by
the Dow Chemical Company as a challenge test for
parameter estimation methods. 
More details about this model can be found in 
\cite{becerra2001applying}. 
The model consists of six differential mass balance
equations, an electroneutrality condition and 
three equilibrium conditions,  
\begin{equation}\label{eq:DAE3}
\begin{cases}
\dot u_1 =  -k_2 u_2 v_2,
\\
\dot u_2  = -k_1 u_2 u_6 + k_{-1} u_4 - k_2 u_2 v_2,
\\
\dot u_3 = k_2 u_2 v_2 + k_3 u_4 u_6 - k_{-3} v_3,
\\
\dot u_4 = -k_3 u_4 u_6 + k_{-3} v_3,
\\
\dot u_5 = _1 u_2 u_6 -k_{-1} v_4,
\\
\dot u_6 = -k_1 u_2 u_6 + k_{-1} v_4 - k_3 v_4 v_6 + k_{-3} v_3,
\\
{\rm [Q^+]} - u_6 + v_1 - v_2 - v_3 - v_4 = 0,
\\
v_2 - K_2 u_1 / (K_2 + v_1) = 0,
\\
v_3 - K_3 u_3 / (K_3 + v_1) = 0,\\
v_4 - K_1 u_5 / (K_1 + v_1) = 0.
\end{cases}
\end{equation}
In our test, the values of the parameters are 
specified as $k_1=25.1911$, $k_2=43.1042$, $k_{-1}=1.1904 \times 10^5$, $k_{-3}= \frac12 k_{-1}$, 
$K_1 = 2.575 \times 10^{-2}$, 
$K_2 = 4.876$, and $K_3 = 1.7884 \times 10^{-2}$, 
and ${\rm [Q^+]} = 0.0131$. 
Assume that one can collect many sets of the time-series measurement data of 
${\bf u}$ and ${\bf v}$ starting from any initial states ${\bf u}_0^{(m)} \in D=[0.6,1.6]\times [6.5,8.5] \times [0,0.7] 
\times [0,0.3] \times [0,0.3] \times [0,0.02] $. 
We collect $500,000$ bursts of trajectory data with $\Delta t=10^{-4}$,  $J_m=2$ and 
${\bf u}_0^{(m)}$ 
sampled from the tensorized Chebyshev distribution on $D$. 
The data set is big as to enable possible accurate discovery of 
the 
complicated reaction process. Handling such a (relatively) large-scale problem 
with the standard regression algorithms can be challenging, primarily
due to the large size of 
matrices. Without using more sophisticated regression methods to deal
with large matrices (which is a topic outside our scope), we 
employ the sequential approximation (SA) algorithm to demonstrate its
capability of handling large data sets. We employ 
the nominalized Legendre polynomial basis of up to
order $n=12$. This induces $N=18,564$ basis functions.
The convergence result is examined by the errors in the function approximation $\tilde {\bf f}^{(k)}, \tilde {\bf g}^{(k)} $ in the $k$-th step iteration 
compared with the exact functions ${\bf f},{\bf g}$. 
To evaluate the errors we independently
draw a set of 20,000 samples uniformly 
and compute the difference between $\tilde {\bf f}^{(k)}, \tilde {\bf g}^{(k)} $ and ${\bf f},{\bf g}$ at these sampling
points. The vector 2-norm is then used to report the difference. 
The relative errors in each iteration of 
SA are displayed in Figure \ref{fig:DAE3err} 
with respect to the number of iterations $k$.  
The SA method uses one piece of data at a time, and 
exhibits good effectiveness in this case. 
As the iteration continues, more data are used and more accurate 
approximations are obtained, yielding 
more accurate models. The exponential
type convergence of the errors is observed. This is consistent with the theoretical analysis in \cite{ShinXiu_SISC17,WuXiu_JCP18,ShinWuXiu_JCP18}. 
We also validate the finial system with the initial state 
${\bf u}_0 = (1.5776, 8.32, 0, 0, 0, 0.0142)^\top$, 
and observe good agreements between the solutions of 
the recovered system and the exact system in Figure \ref{fig:DAE3val}.

\begin{figure}[htbp]
	\centering
	\includegraphics[width=0.48\textwidth]{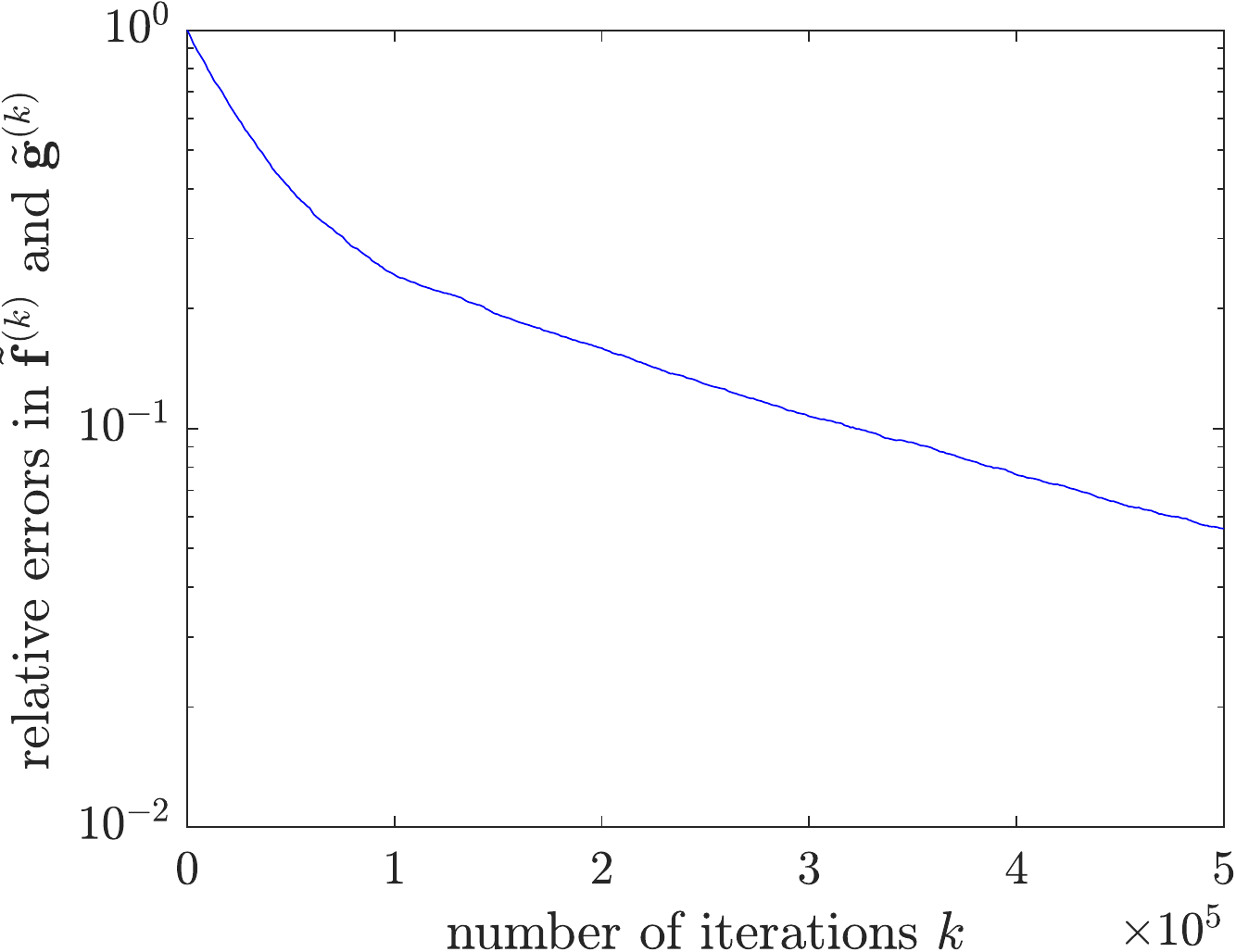}
	\caption{\small
		Example \ref{ex:DAE3}: Convergence history of the sequential approximation. 
	}\label{fig:DAE3err}
\end{figure}

\begin{figure}[htbp]
	\centering
	{\includegraphics[width=0.47\textwidth]{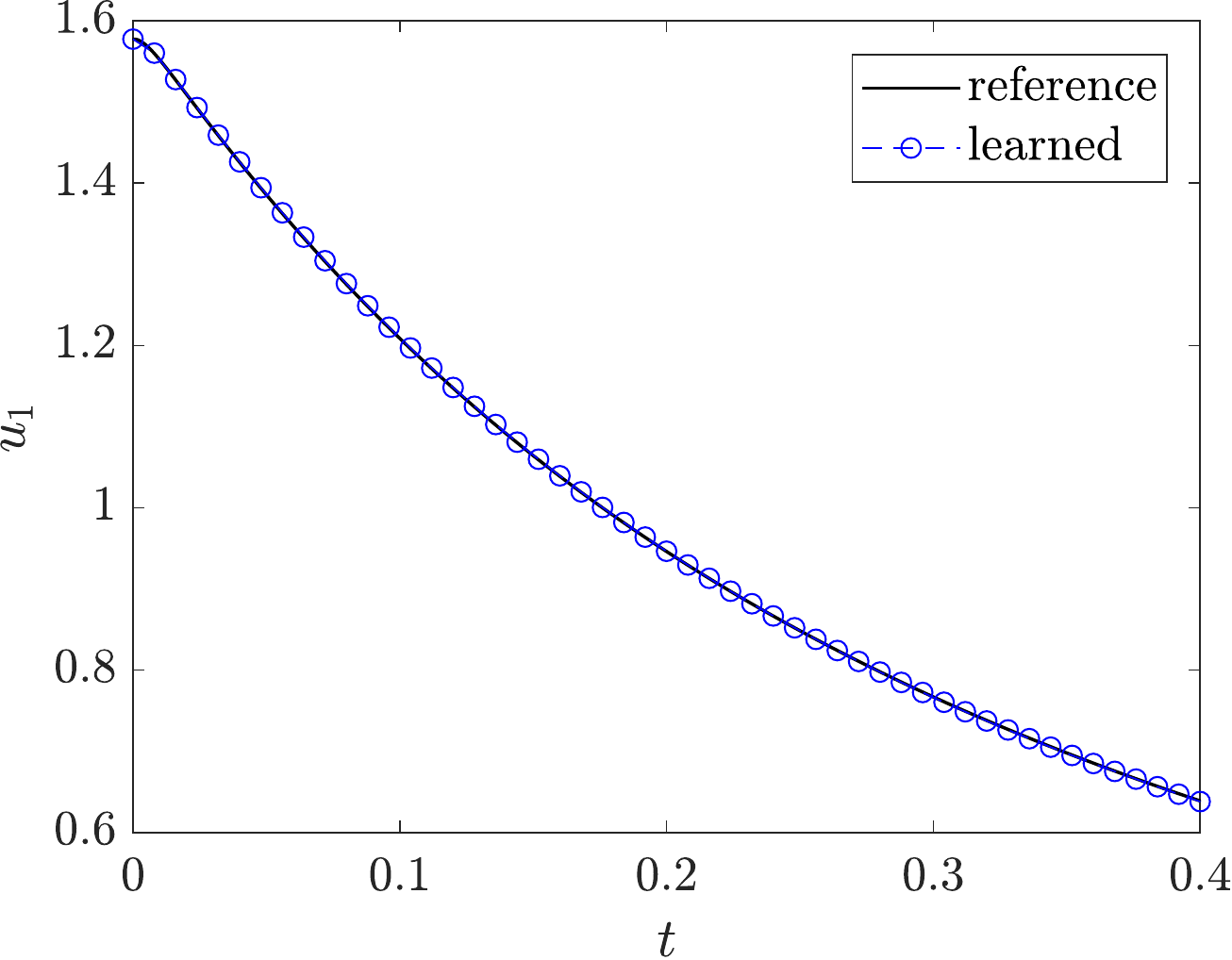}}~~~~
	{\includegraphics[width=0.47\textwidth]{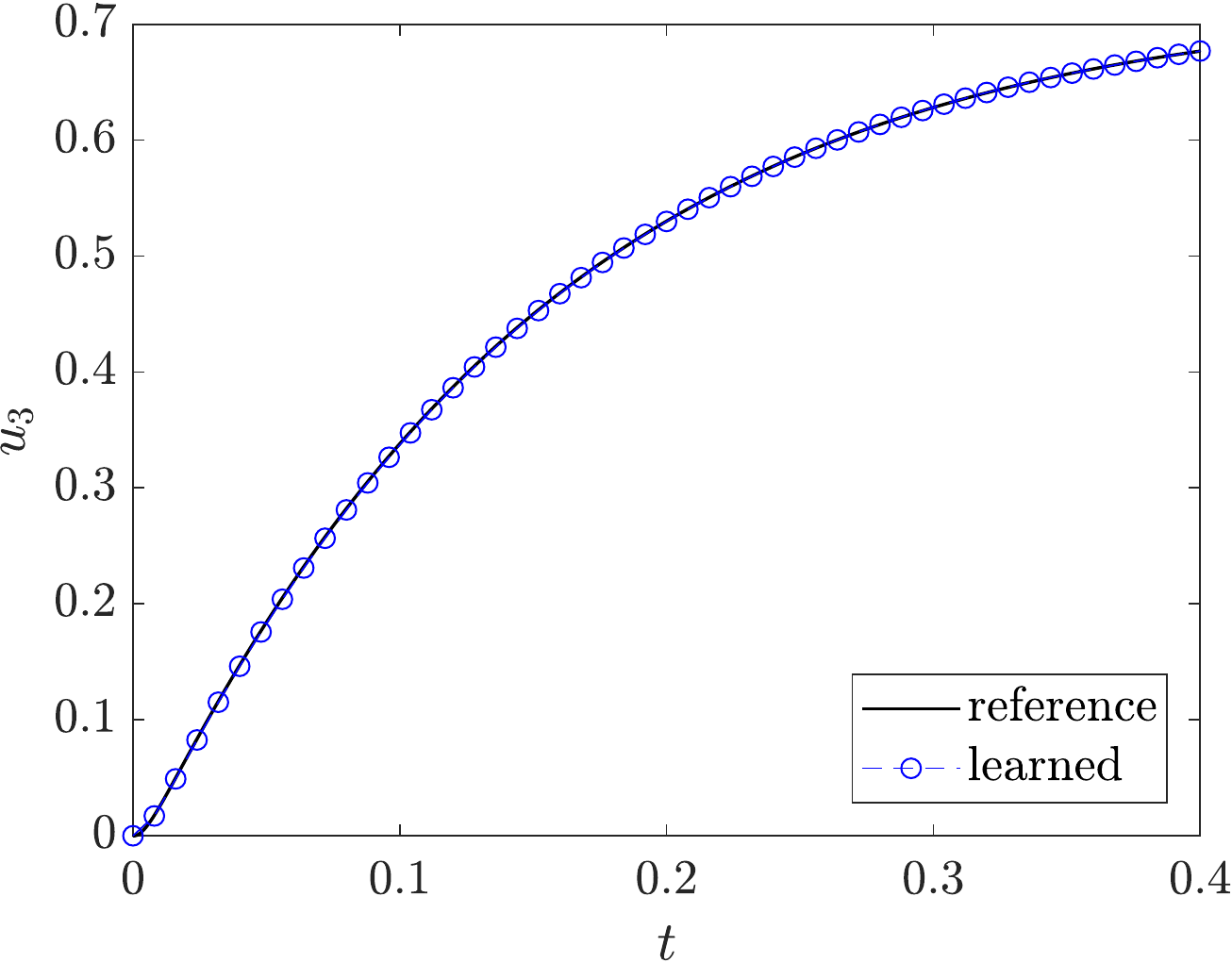}}
	{\includegraphics[width=0.47\textwidth]{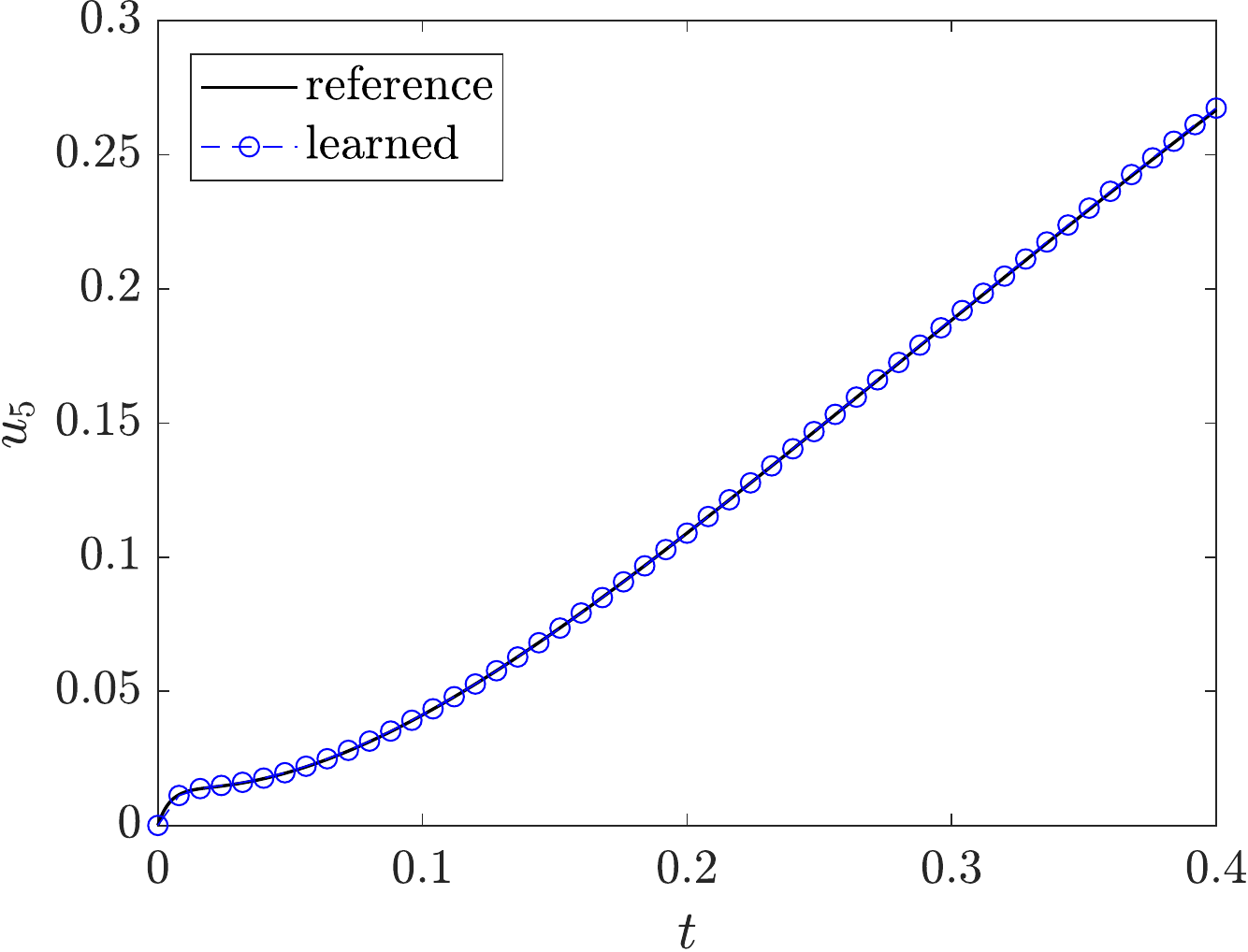}}~~~~
	{\includegraphics[width=0.47\textwidth]{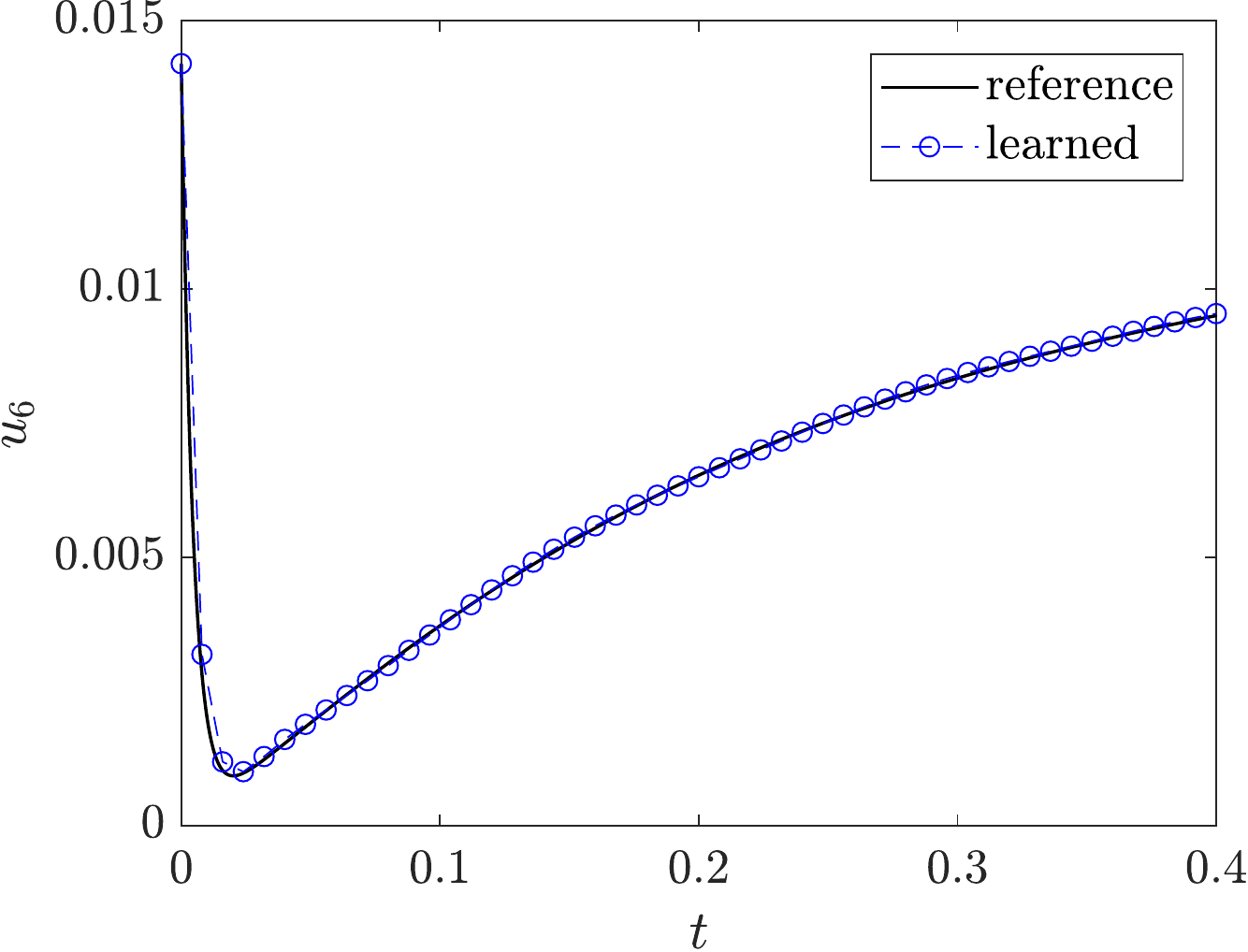}}
	{\includegraphics[width=0.47\textwidth]{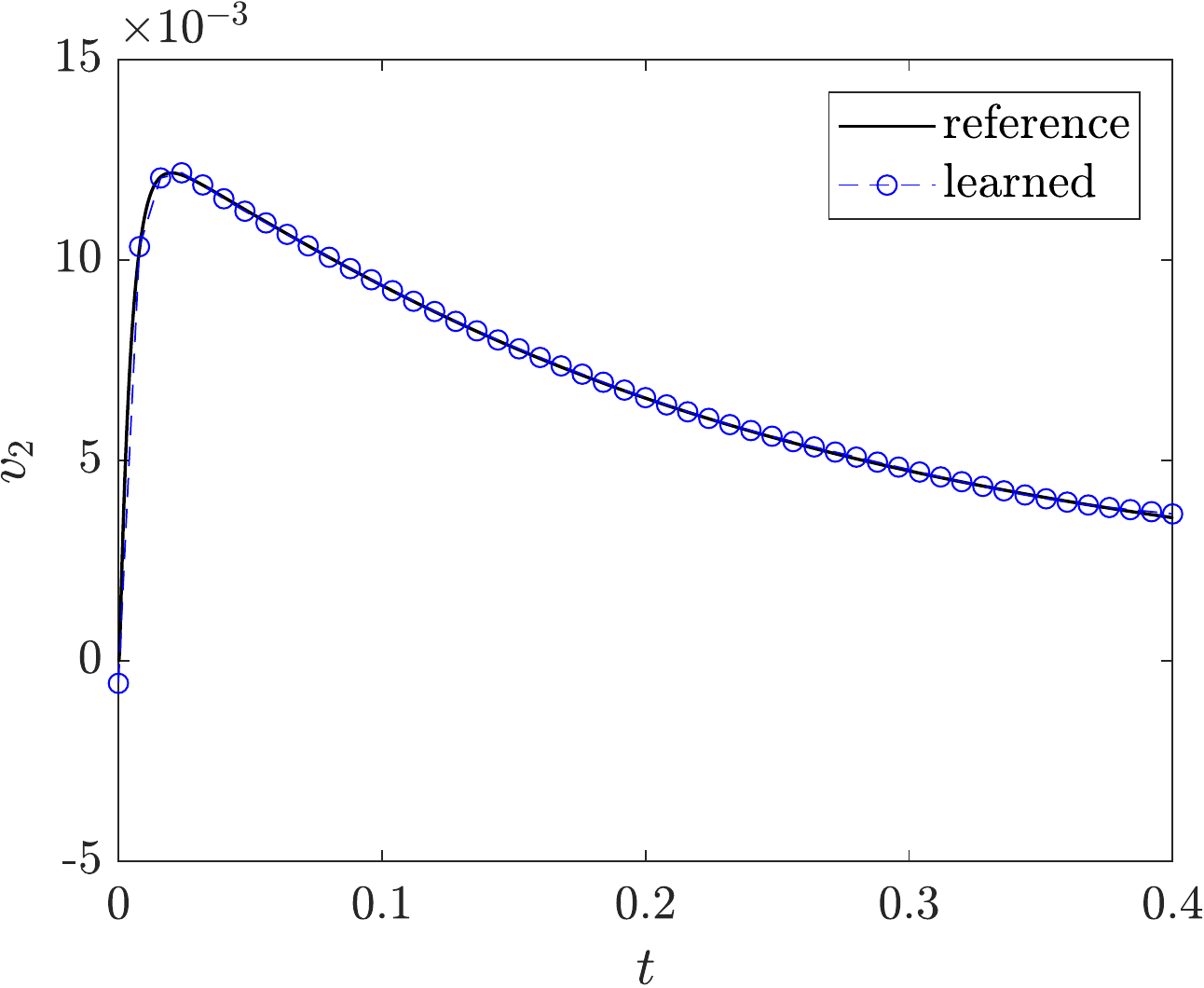}}~~~~
	{\includegraphics[width=0.47\textwidth]{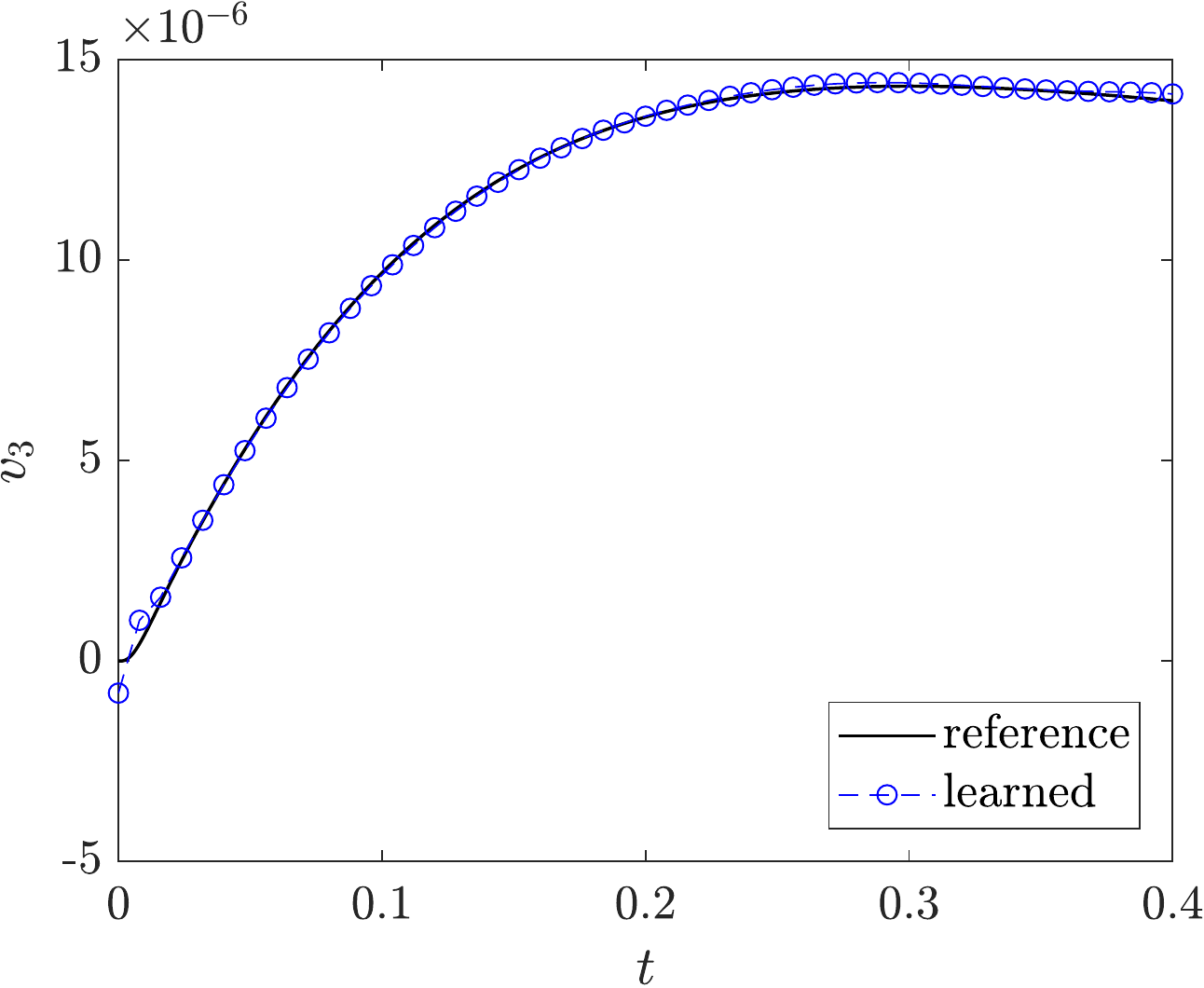}}		
	\caption{\small
		Example \ref{ex:DAE3}: Validation of the learned system. The solid lines denote the solution of the exact system, 
		and the symbols ``$\circ$'' represent 
		the solution of the learned system. The horizontal axis denotes the time.   
	}\label{fig:DAE3val}
\end{figure}

\end{example}

\section{Conclusion} \label{sec:conclusions}

In this paper we studied several effective 
numerical algorithms for data-driven discovery of nonlinear
differential/algebraic equations. We proposed to use the standard basis
functions such as polynomials to construct accurate approximation of
the true governing equations, rather than the exact recovery. We also
discussed the importance of using a (large) number of short bursts of
trajectory data, rather than data from a single (long) trajectory. In
conjunction with some standard approximation algorithms, e.g., least
squares method, the overall method can produce highly accurate
approximation of the true governing equations.
Using an extensive set of test examples, ranging from textbook
examples to more complicated practical examples, we demonstrated that the
proposed method can be highly effective in many situations.

\appendix

\section{Proof of Theorem \ref{thm:errorx2}} \label{app:bound}

We first introduce the following lemma. 

\begin{lemma}\label{thm:funLinfL2}
	If ${\bf f} \in L^\infty(D) $, we have
	\begin{equation} \label{eq:B2infF}
	\| \tilde{\bf f} -  {\bf f} \|_{2,L^\infty} \le 
	\|{\bf f}-{\mathcal P}_V {\bf f}\|_{2,L^\infty} + \| \tilde {\bf f}  - {\mathcal P}_V  {\bf f}  \|_{2,L^2_\omega}  \| \sqrt{K} \|_{L^\infty}. 
	\end{equation}
\end{lemma}

\begin{proof}
	The triangular inequality gives
	\begin{equation}\label{eq:trieq2}
	\| \tilde{\bf f} - {\bf f} \|_{2,L^\infty} \le \|{\bf f}-{\mathcal P}_V {\bf f}\|_{2,L^\infty} + \| \tilde{\bf f} - {\mathcal P}_V {\bf f} \|_{2,L^\infty}.
	\end{equation} 
	Let ${\bf \Psi}=(\psi_1,\cdots,\psi_N)^\top$ be any {\em orthogonal} basis of $V$, and 
	$$
	{\mathcal P}_V  f_\ell=: \langle \hat{\bf c }_{\ell} , {\bf \Psi}({\bf x})\rangle, \qquad
	\tilde f_\ell=: \langle \tilde{\bf c }_{\ell} , {\bf \Psi}({\bf x})\rangle. 
	$$
	Then $K(\x) = \| {\bf \Psi} (\x)\|_2^2 $. Note for any ${\bf x} \in D$ that
	\begin{eqnarray*}
		\| \tilde{\bf f} ({\bf x}) - {\mathcal P}_V {\bf f} ({\bf x}) \|_2^2 & = & \sum_{\ell=1}^d
		\big| \tilde f_\ell ({\bf x}) - {\mathcal P}_V  f_\ell ({\bf x}) \big|^2
		= \sum_{\ell=1}^d
		\big| \langle \tilde{\bf c }_{\ell} -  \hat{\bf {c} }_{\ell}, {\bf \Psi}({\bf x})\rangle  \big|^2
		\\
		&\le &  \sum_{\ell=1}^d
		\|  \tilde {\bf c }_{\ell} -  \hat{\bf {c} }_{\ell} \|_2^2 \| {\bf \Psi}({\bf x}) \|_2^2
		=   \sum_{\ell=1}^d \| \tilde f_\ell  - {\mathcal P}_V  f_\ell  \|_{L^2_\omega}^2   K(\x) \\
		& = & \| \tilde {\bf f}  - {\mathcal P}_V  {\bf f}  \|_{2,L^2_\omega}^2     K(\x) ,
	\end{eqnarray*}
	where the Cauchy-Schwarz inequality and the orthogonality of basis $\{\psi_j (\x) \}$ have been used.
	Taking $L^\infty$-norm in
	the above inequality and using \eqref{eq:trieq2} give \eqref{eq:B2infF}.
\end{proof}

Based on the above lemma, the proof of Theorem \ref{thm:errorx2} is given as follows. 

\begin{proof}
	Integrating the systems \eqref{eq:ODE} and \eqref{eq:appODE} respectively gives
	\begin{align}
	&{\bf u}(t) = {\bf u}(t_0) + \int_{t_0}^t {\bf f} ( {\bf u} (s) ) ds, \label{eq:xex1}
	\\
	&{\bf x}(t) = {\bf x}(t_0) + \int_{t_0}^t  \tilde{\bf f} ( {\bf x} (s) ) ds. \label{eq:xapp1}
	\end{align}
	Subtracting \eqref{eq:xex1} from \eqref{eq:xapp1}  with ${\bf u}(t_0)={\bf x}(t_0) $, we obtain
	\begin{eqnarray*}
		\big\| {\bf x}(t) - {\bf u}(t)  \big\|_2  &=&  \bigg\| \int_{t_0}^t \tilde{\bf f} ( {\bf x} (s) )  - {\bf f} ( {\bf u} (s) ) ds \bigg\|_2
		\\
		&=& (t-t_0) \eta \bigg(\frac{1}{t-t_0} \int_{t_0}^t \tilde{\bf f} ( {\bf x} (s) )  - {\bf f} ( {\bf u} (s) )  ds  \bigg) ,
	\end{eqnarray*}
	where $\eta ({\bf x}) = \| {\bf x} \|_2 $ is a convex function satisfying the Jensen's inequality
	$$
	\eta \bigg(\frac{1}{t-t_0} \int_{t_0}^t \tilde{\bf f} ( {\bf x} (s) )  - {\bf f} ( {\bf u} (s) )  ds  \bigg)
	\le \frac{1}{t-t_0} \int_{t_0}^t  \eta \big(\tilde{\bf f} ( {\bf x} (s) )  - {\bf f} ( {\bf u} (s) ) \big) ds .
	$$
	This together with the triangular inequality yield
	\begin{eqnarray*}
		\big\| {\bf x}(t) - {\bf u}(t)  \big\|_2
		& \le &
		\int_{t_0}^t \big\| \tilde{\bf f} ( {\bf x} (s) )  - {\bf f} ( {\bf u} (s) ) \big\|_2 ds
		\\
		& \le & \int_{t_0}^t
		\big\| \tilde{\bf f} ( {\bf x} (s) )  - {\bf f} ( {\bf x} (s) ) \big\|_2
		+ \big\| {\bf f} ( {\bf x} (s) )  - {\bf f} ( {\bf u} (s) ) \big\|_2 ds
		\\
		& \le &  \int_{t_0}^t   \| \tilde{\bf f}   - {\bf f}  \|_{2,L^\infty} + L_f \big\|  {\bf x} (s)   -  {\bf u} (s)  \big\|_2 ds,
	\end{eqnarray*}
	which implies \eqref{eq:L2Bx}. Applying Gr\"onwall's inequality to \eqref{eq:L2Bx} gives
	\begin{eqnarray*}
		\big\| {\bf x}(t) - {\bf u}(t)  \big\|_2 & \le &  \bigg( t-t_0+ L_f \int_{t_0}^t (s-t_0) \exp\big( L_f(t-s) \big) ds \bigg) \| \tilde{\bf f} - {\bf f} \|_{2,L^\infty}
		\\
		& = & L_f^{-1} \big( e^{L_f (t-t_0)} - 1 \big) \| \tilde{\bf f} - {\bf f} \|_{2,L^\infty},
	\end{eqnarray*}
	Then utilizing \eqref{eq:B2infF}, we get \eqref{eq:L2BxFinal} and complete the proof.
\end{proof}

\section*{Acknowledgment}
This work was partially supported by AFOSR FA95501410022 and NSF DMS 1418771.

\section*{References}

\bibliographystyle{plain}
\bibliography{Kaczmarz,least-squares,collocation,random,learningEqs}

\end{document}